\pdfoutput=1

\documentclass{amsart}
\providecommand{\noopsort}[1]{} 
\usepackage[mathscr]{eucal}
\usepackage{amssymb}
\usepackage[usenames,dvipsnames]{xcolor} 
\usepackage[normalem]{ulem}
\usepackage{amsthm}
\usepackage{bbold}
\usepackage{comment}
\usepackage{enumitem}
\usepackage{amsmath}
\usepackage{wasysym}
\usepackage{tikz-cd}
\usepackage{calc}
\usepackage[backend=biber,style=alphabetic,sorting=nyt,doi=false,isbn=false,url=false,maxalphanames=5,maxcitenames=50,maxnames=50]{biblatex}
\usepackage{xpatch}
\DeclareLabelalphaNameTemplate{
  \namepart[use=true, pre=true, strwidth=1, compound=true]{prefix}
  \namepart[compound=true]{family}
}
\DeclareFieldFormat[article,unpublished,incollection,inproceedings]{title}{#1}
\renewbibmacro{in:}{}
\xpatchbibdriver{misc}
{\printfield{note}}
{\printfield{note}\nopunct,}
 {}{}
\xpatchbibdriver{unpublished}
{\printfield{note}}
{\printfield{note}\nopunct,}
 {}{}
\xpatchbibdriver{book}
{\printfield{volume}}
{\nopunct, \printfield{volume}}
 {}{}
\xpatchbibdriver{book}
{\usebibmacro{chapter+pages}}
 {}
 {}{}
\DeclareFieldFormat[book]{volume}{volume #1 of\nopunct}
\DeclareFieldFormat{pages}{#1} 
\DeclareFieldFormat[book,incollection,inproceedings]{series}{\emph{#1}}
\DeclareBibliographyDriver{inproceedings}{%
  \usebibmacro{bibindex}%
  \usebibmacro{begentry}%
  \usebibmacro{author/translator+others}%
  \setunit{\printdelim{nametitledelim}}\newblock
  \usebibmacro{title}%
  \newunit
  \printlist{language}%
  \newunit\newblock
  \usebibmacro{byauthor}%
  \newunit\newblock
  \addperiod\addspace
  In\nopunct
  \usebibmacro{maintitle+booktitle}%
  \usebibmacro{chapter+pages}%
  \usebibmacro{publisher+location+date}%
  \usebibmacro{finentry}}
\DeclareBibliographyDriver{incollection}{%
  \usebibmacro{bibindex}%
  \usebibmacro{begentry}%
  \usebibmacro{author/translator+others}%
  \setunit{\printdelim{nametitledelim}}\newblock
  \usebibmacro{title}%
  \newunit
  \printlist{language}%
  \newunit\newblock
  \usebibmacro{byauthor}%
  \newunit\newblock
  \addperiod\addspace
  In\nopunct
  \usebibmacro{maintitle+booktitle}%
  \usebibmacro{chapter+pages}%
  \usebibmacro{publisher+location+date}%
  \usebibmacro{finentry}}
\renewbibmacro*{journal+issuetitle}{%
\usebibmacro{journal}%
  \setunit*{\addspace}%
  \iffieldundef{series}
    {}
    {\newunit
     \printfield{series}%
     \setunit{\addspace}}%
  \printfield{volume}%
  \printfield[parens]{number}%
  \iffieldundef{pages}
   {}
   {\addcolon}
  \printfield{pages}%
  \addcomma%
  \addspace%
  \usebibmacro{date}%
  \newunit}
\renewbibmacro*{note+pages}{%
	\printfield{note}%
\newunit}

\usepackage{etoolbox} 

\usepackage[unicode]{hyperref} 

\definecolor{dark-red}{rgb}{0.5,0.15,0.15}
\definecolor{dark-blue}{rgb}{0.15,0.15,0.6}
\definecolor{dark-green}{rgb}{0.15,0.6,0.15}
\hypersetup{
    colorlinks, linkcolor=Blue,
    citecolor=Blue, urlcolor=Blue
}

\usepackage[nameinlink,capitalise,noabbrev]{cleveref}

\newtheorem{thmx}{Theorem}

\numberwithin{equation}{section}
\usepackage[all]{xy}
\xyoption{line}
\usepackage{graphicx}
\usepackage{mathtools}

\swapnumbers 

\newtheorem{Thm}[equation]{Theorem}
\newtheorem*{Thm*}{Theorem}
\newtheorem*{MainThm*}{Main Theorem}
\newtheorem{Prop}[equation]{Proposition}
\newtheorem{Lem}[equation]{Lemma}
\newtheorem{Cor}[equation]{Corollary}

\newtheorem{Conj}[equation]{Conjecture}

\theoremstyle{remark}
\newtheorem{Def}[equation]{Definition}
\newtheorem{Ter}[equation]{Terminology}
\newtheorem{Not}[equation]{Notation}
\newtheorem{Exa}[equation]{Example}

\newtheorem{Rem}[equation]{Remark}

\tikzset{
    labelrotatebelow/.style={anchor=north, rotate=90, inner sep=1.0mm}
}
\tikzset{
    labelrotateabove/.style={anchor=south, rotate=90, inner sep=1.0mm}
}

\newcommand{\nc}{\newcommand}
\nc{\dmo}{\DeclareMathOperator}

\renewcommand{\emptyset}{\varnothing}

\nc{\Beren}[1]{{\color{MidnightBlue}#1}}
\nc{\Drew}[1]{{\color{OliveGreen}#1}}
\nc{\Tobi}[1]{{\color{Orange}#1}}
\nc{\Dout}[1]{\Drew{\sout{#1}}}
\nc{\Bout}[1]{\Beren{\sout{#1}}}
\nc{\Tout}[1]{\Tobi{\sout{#1}}}

\usepackage{todonotes}

\nc{\overbar}[1]{\mkern 1.5mu\overline{\mkern-1.5mu#1\mkern-1.5mu}\mkern 1.5mu}
\nc{\Fbar}{\mkern 2.5mu\overline{\mkern-2.5mu F}}

\nc{\SuppBIK}{\Supp_{\text{BIK}}}

\nc{\kappaaux}{g}
\nc{\kappaM}{{\kappaaux({\frak M})}}
\nc{\kappam}{{\kappaaux({\frak m})}}
\nc{\kappaP}{{\kappaaux(\cat P)}}
\nc{\kappaQ}{{\kappaaux(\cat Q)}}
\nc{\kappaCP}{{\kappaaux_{\cat C}(\cat P)}}
\nc{\kappaDP}{{\kappaaux_{\cat D}(\cat P)}}
\nc{\kappaCQ}{{\kappaaux_{\cat C}(\cat Q)}}
\nc{\kappaDQ}{{\kappaaux_{\cat D}(\cat Q)}}
\nc{\kappaTQ}{{\kappaaux_{\cat T}(\cat Q)}}
\nc{\kappaTViQ}{{\kappaaux_{\cat T(V_i)}(\cat Q)}}
\nc{\kappaphiB}{{\kappaaux(\phi(\cat B))}}
\nc{\kappaphiQ}{{\kappaaux(\varphi(\cat Q))}}
\nc{\kappaphiiQ}{{\kappaaux(\varphi_i(\cat Q))}}
\nc{\kappaCphiQ}{{\kappaaux_{\cat C}(\varphi(\cat Q))}}

\dmo{\Sub}{Sub}
\dmo{\stmod}{stmod}
\dmo{\StMod}{StMod}
\dmo{\Proj}{Proj}
\nc{\SpEn}{\cat S_{E(n)}}
\nc{\SpEnf}{\cat S_n}
\nc{\Lcomp}{L^{\mathrm{com}}} 
\nc{\Ucomp}{U^{\mathrm{com}}}
\nc{\Loco}[1]{\Loc_{\otimes}\hspace{-0.3ex}\langle #1 \rangle}
\nc{\bbullet}{{\scriptscriptstyle\hspace{-1pt}\bullet}}
\nc{\bullett}{{\scriptscriptstyle\bullet}\hspace{-1pt}}
\nc{\LF}{L\hspace{-0.2ex}F}
\nc{\SpG}{\Sp_G}
\nc{\SpGEn}{\Sp_{G,E(n)}}
\nc{\EG}{\bbE_G}
\nc{\EH}{\bbE_H}
\nc{\DEG}{\Der(\EG)}
\nc{\DEH}{\Der(\EH)}
\nc{\DE}{\Der(\bbE)}
\nc{\Prst}{{\cat P}\mathrm{r^{st}}}
\nc{\Mack}[2]{\mathrm{Mack}_{#1}(#2)}
\nc{\SC}{S\cat C}
\dmo{\fin}{{fin}}
\dmo{\DM}{DM}
\dmo{\fp}{fp}
\nc{\DMQ}{\DM_Q}
\dmo{\DerKal}{DMack}
\dmo{\Der}{D}
\dmo{\Derqc}{D_{qc}}
\dmo{\Derperf}{D_{perf}}
\dmo{\DMot}{DMot}
\dmo{\rmH}{H}
\dmo{\piu}{\underline{\pi}}
\dmo{\Sphere}{\mathbb{S}}
\nc{\HA}{{\rmH \hspace{-0.2em}\bbA}}
\nc{\HZ}{{\rmH \hspace{-0.2em}\bbZ}}
\nc{\HR}{{\rmH \hspace{-0.15em}R}}
\nc{\HZbar}{{\rmH \hspace{-0.2em}\underline{\bbZ}}}
\nc{\Fp}{{\bbF_{\hspace{-0.1em}p}}}
\nc{\HFp}{{\rmH \hspace{-0.15em}\bbF_{\hspace{-0.1em}p}}}
\nc{\DHZG}{\Der(\HZ_G)}
\nc{\DHZH}{\Der(\HZ_H)}
\nc{\DHZK}{\Der(\HZ_K)}
\nc{\DHZGN}{\Der(\HZ_{G/N})}
\nc{\DHZGG}{\Der(\HZ_{G/G})}
\nc{\DHZCp}{\Der(\HZ_{C_p})}
\nc{\DHZGprime}{\Der(\HZ_{G'})}
\nc{\DHZ}{\Der(\HZ)}
\nc{\frakp}{\mathfrak{p}}
\nc{\frakq}{\mathfrak{q}}
\nc{\Z}{\mathbb{Z}}
\nc{\SSG}{\text{sSet}_*^G}
\nc{\sSet}{\text{sSet}}

\dmo{\csupp}{csupp}
\dmo{\Con}{Conj}
\dmo{\Id}{Id}
\dmo{\Loc}{Loc}
\dmo{\rmK}{\textrm{\rm K}}
\dmo{\Spc}{Spc}
\dmo{\thick}{thick}
\nc{\thickt}[1]{\thick_\otimes\langle #1 \rangle}
\dmo{\cone}{cone}
\dmo{\End}{End}
\dmo{\Mor}{Mor}
\dmo{\Hom}{Hom}
\dmo{\id}{id}
\dmo{\incl}{incl}
\dmo{\Img}{Im}
\dmo{\im}{im}
\dmo{\Ker}{Ker}
\dmo{\ind}{ind}
\dmo{\CoInd}{coind}
\dmo{\res}{res}
\dmo{\infl}{infl}
\dmo{\triv}{triv}
\dmo{\Tel}{Tel} 
\dmo{\grMod}{grMod}%
\dmo{\Mod}{Mod}%
\dmo{\opname}{op}
\dmo{\SH}{SH}
\nc{\SHfin}{\SH^{\fin}}
\dmo{\smallb}{b}
\dmo{\Spec}{Spec}
\dmo{\supp}{supp}
\dmo{\Supp}{Supp}
\nc{\SHc}{{\SH^c}}
\nc{\SHp}{{\SH_{(p)}}}
\nc{\SHcp}{{\SH^c_{(p)}}}
\nc{\SHG}{\SH_G}
\nc{\SHGEn}{\SH_{G,E(n)}}
\nc{\SHGp}{\SH_{G,(p)}}
\nc{\SHGc}{\SHG^c}
\nc{\quadtext}[1]{\quad\textrm{#1}\quad}
\nc{\qquadtext}[1]{\qquad\textrm{#1}\qquad}
\nc{\adj}{\dashv}
\nc{\adjto}{\rightleftarrows}
\nc{\bbL}{\mathbb{L}}
\nc{\bbA}{\mathbb{A}}
\nc{\bbE}{\mathbb{E}}
\nc{\bbN}{\mathbb{N}}
\nc{\bbQ}{\mathbb{Q}}
\nc{\bbZ}{\mathbb{Z}}
\nc{\bbF}{\mathbb{F}}
\nc{\cat}[1]{\mathscr{#1}}
\nc{\ie}{{\sl i.e.}, }
\nc{\into}{\mathop{\rightarrowtail}}
\nc{\inv}{^{-1}}
\nc{\isoto}{\mathop{\overset{\sim}\to}}
\nc{\isotoo}{\mathop{\overset{\sim}\too}}
\nc{\onto}{\mathop{\twoheadrightarrow}}
\nc{\too}{\mathop{\longrightarrow}\limits}
\nc{\mapstoo}{\longmapsto}
\nc{\adh}[1]{\overline{#1}}
\nc{\adhpt}[1]{\adh{\{#1\}}}
\nc{\aka}{{a.\,k.\,a.}\ }
\nc{\calF}{\mathcal{F}}
\nc{\eg}{{\sl e.\,g.}}
\nc{\Homcat}[1]{\Hom_{\cat #1}}
\nc{\hook}{\hookrightarrow}
\nc{\ideal}[1]{\langle #1\rangle}
\nc{\ihom}{{\underline{\hom}}}
\nc{\Mid}{\,\big|\,}
\nc{\MMod}{\,\text{-}\Mod}%
\nc{\op}{^{\opname}}
\nc{\oto}[1]{\overset{#1}\to}
\nc{\otoo}[1]{\overset{#1}{\,\too\,}}
\nc{\sminus}{\!\smallsetminus\!}
\nc{\poplus}[1]{^{\oplus #1}}%
\nc{\potimes}[1]{^{\otimes #1}}
\nc{\sbull}{{\scriptscriptstyle\bullet}}
\nc{\SET}[2]{\big\{\,#1\Mid#2\,\big\}}
\nc{\SpcK}{\Spc(\cat K)}
\nc{\then}{\Rightarrow}
\nc{\unit}{\mathbb{1}}
\nc{\xra}{\xrightarrow}
\nc{\phigeom}[1]{\widetilde{\Phi}^{#1}}
\nc{\phigeomb}[1]{\Phi^{#1}}
\nc{\PhiHx}{\Phi^H\hspace{-0.25ex}(x)}
\nc{\PhiHX}{\Phi^H\hspace{-0.25ex}(X)}
\dmo{\Oname}{O}
\dmo{\proper}{proper}
\dmo{\lenormal}{\unlhd}
\dmo{\lnormal}{\lhd}
\nc{\normal}{\trianglelefteq}
\nc{\Op}{\Oname^p}
\nc{\Oq}{\Oname^q}
\dmo{\Sp}{Sp}
\dmo{\Ho}{Ho}
\dmo{\Fin}{Fin}
\dmo{\add}{add}
\dmo{\Fun}{Fun}
\dmo{\Ext}{Ext}
\dmo{\CAlg}{CAlg}
\dmo{\CMon}{CMon}
\dmo{\CC}{\cat C} 
\dmo{\DD}{\cat D}
\dmo{\OO}{\mathcal{O}}
\dmo{\Map}{Map}
\dmo{\Span}{Span}
\dmo{\N}{N}
\dmo{\Cat}{Cat}
\dmo{\colim}{colim}
\dmo{\hocolim}{hocolim}
\dmo{\Ch}{Ch}
\dmo{\A}{\mathbb{A}^{eff}}
\nc{\AGeff}{\mathbb{A}_G^{\mathrm{eff}}}
\nc{\BGeff}{\mathcal{B}_G^{\mathrm{eff}}}
\nc{\BG}{{\mathcal{B}_G}}
\nc{\NBGeff}{{\N}{\BGeff}}
\dmo{\Ab}{Ab}
\dmo{\Set}{Set}
\dmo{\ev}{ev}
\dmo{\Spcl}{Spcl}
\nc{\Funadd}{\Fun_{\add}}
\dmo{\proj}{proj}
\dmo{\cof}{cof}

\dmo{\Coideal}{Coideal}
\dmo{\gen}{gen}

\Crefname{Exa}{Example}{Examples}
\Crefname{Rem}{Remark}{Remarks}
\Crefname{thmx}{Theorem}{Theorems}

\nc{\mT}{\kern-0.5em\mod\kern-0.1em\text{-}\cat{T}^c}
\nc{\MT}{\Mod\kern-0.1em\text{-}\cat{T}}
\newcounter{enum-resume-hack}

\begin{document}


\title[Stratification in tensor triangular geometry]{Stratification in tensor triangular geometry with applications to spectral Mackey functors}

\author{Tobias Barthel}
\author{Drew Heard}
\author{Beren Sanders}
\date{\today}

\makeatletter
\patchcmd{\@setaddresses}{\indent}{\noindent}{}{}
\patchcmd{\@setaddresses}{\indent}{\noindent}{}{}
\patchcmd{\@setaddresses}{\indent}{\noindent}{}{}
\patchcmd{\@setaddresses}{\indent}{\noindent}{}{}
\makeatother

\address{Tobias Barthel, Max Planck Institute for Mathematics, Vivatsgasse 7, 53111 Bonn, Germany}
\email{tbarthel@mpim-bonn.mpg.de}
\urladdr{\href{https://sites.google.com/view/tobiasbarthel/home}{https://sites.google.com/view/tobiasbarthel/home}}

\address{Drew Heard, Department of Mathematical Sciences, Norwegian University of Science and Technology, Trondheim}
\email{drew.k.heard@ntnu.no}
\urladdr{\href{https://folk.ntnu.no/drewkh/}{https://folk.ntnu.no/drewkh/}}

\address{Beren Sanders, Mathematics Department, UC Santa Cruz, 95064 CA, USA}
\email{beren@ucsc.edu}
\urladdr{\href{http://people.ucsc.edu/~beren/}{http://people.ucsc.edu/$\sim$beren/}}

\begin{abstract}
We systematically develop a theory of stratification in the context of tensor triangular geometry and apply it to classify the localizing tensor-ideals of certain categories of spectral $G$-Mackey functors for all finite groups~$G$. Our theory of stratification is based on the approach of Stevenson which uses the Balmer--Favi notion of big support for tensor-triangulated categories whose Balmer spectrum is weakly noetherian. We clarify the role of the local-to-global principle and establish that the Balmer--Favi notion of support provides the universal approach to weakly noetherian stratification. This provides a uniform new perspective on existing  classifications in the literature and clarifies the relation with the theory of Benson--Iyengar--Krause. Our systematic development of this approach to stratification, involving a reduction to local categories and the ability to pass through finite \'{e}tale extensions, may be of independent interest. Moreover, we strengthen the relationship between stratification and the telescope conjecture. The starting point for our equivariant applications is the recent computation by Patchkoria--Sanders--Wimmer of the Balmer spectrum of the category of derived Mackey functors, which was found to capture precisely the  height $0$ and height $\infty$ chromatic layers of the spectrum of the equivariant stable homotopy category. We similarly study the Balmer spectrum of the category of $E(n)$-local spectral Mackey functors noting that it bijects onto the height $\le n$ chromatic layers of the  spectrum of the  equivariant stable homotopy category; conjecturally the topologies coincide. Despite our incomplete knowledge of the topology of the Balmer spectrum, we are able to completely classify the localizing tensor-ideals of these categories of spectral Mackey functors.
\end{abstract}


\thanks{The first-named author would like to thank the Max Planck Institute for Mathematics for its hospitality. The second-named author is supported by grant number TMS2020TMT02 from the Trond Mohn Foundation. The third-named author is supported by NSF grant~DMS-1903429.}


\maketitle


{
\hypersetup{linkcolor=black}
\tableofcontents
}

\section*{Introduction}\label{sec:introduction} 

Mathematics is confronted with wild classification problems. For example, there is no hope of classifying the finite-dimensional representations of most finite groups in positive characteristic. Similarly, there is no hope of classifying finite CW-complexes up to homotopy equivalence. With this in mind, an important development of the last thirty years was the realization that it is often possible to classify such mathematical objects up to a weaker notion of equivalence: One works stably, \mbox{that is,} in a suitable stable homotopy category~$\cat K$ of such objects, and considers two objects to be equivalent if they can be built from each other using the natural tensor-triangulated structure of the stable category. Classifying the objects up to this weaker notion of equivalence amounts to classifying the thick $\otimes$-ideals of the category~$\cat K$.

Historical examples of such classification theorems include stable module categories of finite-dimensional representations in modular representation theory, the stable homotopy category of finite spectra in algebraic topology, and derived categories of perfect complexes in algebraic geometry \cite{BensonCarlsonRickard97,HopkinsSmith98,Thomason97}. All of these results are effected by a suitable theory of ``support'' for the objects of the category. For example, the classification theorem in modular representation theory uses the theory of support varieties (developed by Quillen, Carlson and others) which assigns to each finite-dimensional $k$-linear representation of $G$ a certain closed subset of the projective scheme $\Proj(H^*(G;k))$ associated to the group cohomology ring. Two such representations are equivalent in the above sense precisely when their support varieties coincide. Similarly, the Hopkins--Smith classification of finite spectra uses a notion of support defined using the Morava $K$-theories.

These theorems were greatly clarified by Balmer \cite{Balmer05a} who introduced, for any essentially small tensor-triangulated category $\cat K$, the universal notion of support for the objects of $\cat K$. It assigns to each object $x \in \cat K$ a closed subset $\supp(x)$ of a certain topological space~$\Spc(\cat K)$ called the Balmer spectrum. Two objects are equivalent in the above sense precisely when they have the same support. This leads to a bijection between the (radical) thick $\otimes$-ideals of~$\cat K$ and certain subsets (called Thomason subsets) of the space $\Spc(\cat K)$. In this way, the classification problem becomes the question of giving an explicit description of the Balmer spectrum and its universal notion of support.

While these results are satisfying, they only apply to essentially small categories of ``compact'' objects. For example, they only classify the finite-dimensional representations, the finite spectra, and the perfect complexes. Typically the category~$\cat K$ is the subcategory of compact objects $\cat K=\cat T^c \subset \cat T$ in a larger rigidly-compactly generated tensor-triangulated category $\cat T$. Many of the most interesting objects are non-compact objects in $\cat T$ and we would thus like to classify the localizing $\otimes$-ideals of the big category $\cat T$. For example, the stable homotopy category of finite spectra is the subcategory of compact objects $\SHfin = \SHc \subset \SH$ in the stable homotopy category of all spectra, and a classification of the localizing $\otimes$-ideals of $\SH$ would provide a coarse classification of cohomology theories in algebraic topology.

We cannot apply Balmer's construction to these large categories $\cat T$ because they are not essentially small, but such set-theoretic issues are not the main point: The problem is that the axioms for Balmer's universal notion of support are not appropriate for big objects. For example, morally the support of a compact object should be closed (and this is one of Balmer's axioms), but this should not be expected for the support of an arbitrary object (like an infinite-dimensional representation). We can still attempt to classify the localizing $\otimes$-ideals of $\cat T$ by defining a notion of support for arbitrary objects in $\cat T$ but such notions of support are less well-behaved in general and the properties they should satisfy are less clear-cut. Ideally, we would like a ``big'' version of the Balmer spectrum --- a universal notion of support for big tensor-triangulated categories which classifies their localizing $\otimes$-ideals ---
but no one has succeeded in constructing such a support theory and there is evidence that in general it cannot exist \cite{BalmerKrauseStevenson20}. Nevertheless, the search for useful theories of support continues and significant positive results have been attained in certain noetherian situations.

For example, in a series of papers \cite{BensonIyengarKrause08,BensonIyengarKrause11b,BensonIyengarKrause11a}, Benson, Iyengar, and Krause (BIK) have developed a theory of support for compactly generated triangulated categories $\cat T$ equipped with an action of a graded noetherian commutative ring~$R$. If $\cat T$ is \emph{tensor}-triangulated, one can take the canonical action of the graded endomorphism ring of the unit object $R\coloneqq \End_{\cat T}^*(\unit)$ provided this ring is noetherian. In general, the BIK notion of support does not classify the localizing $\otimes$-ideals of $\cat T$ but, inspired by ideas of \cite{HoveyPalmieriStrickland97}, BIK develop a powerful condition called \emph{stratification} which is sufficient to obtain such a classification. Indeed this work of BIK culminates in the celebrated classification of localizing $\otimes$-ideals for the big stable module categories of finite groups \cite{BensonIyengarKrause11a} and has since seen many applications \cite{Shamir12,DellAmbrogioStanley16,BarthelCastellanaHeardValenzuela19,BensonIyengarKrausePevtsova18}.

In parallel, Balmer and Favi \cite{BalmerFavi11} have introduced a notion of support for rigidly-compactly generated tensor-triangulated categories $\cat T$ which takes values in the Balmer spectrum $\Spc(\cat T^c)$ of the subcategory of compact objects. This notion of support has been applied to the classification of localizing $\otimes$-ideals by Stevenson~\cite{Stevenson13}. The overarching goal of the present paper is twofold: On the one hand, we will systematically develop the notion of stratification based on the Balmer--Favi notion of support and the Balmer spectrum, and clarify its relation to the theory of BIK. On the other hand, we will demonstrate that this theory of stratification provides a uniform perspective on many known classification problems --- both those that are and are not amenable to the techniques of BIK --- and also apply our methods to new examples, most notably in equivariant homotopy theory.
\[\ast \ast \ast\]

We will develop these ideas in the context of a rigidly-compactly generated tensor-triangulated category $\cat T$ whose Balmer spectrum of compact objects $\Spc(\cat T^c)$ is weakly noetherian, a point-set topological condition that simultaneously generalizes noetherian spaces and profinite spaces. We will take the geometric perspective of tensor triangular geometry. To guide intuition, it is helpful to picture~$\cat T$ as a bundle of tensor-triangulated categories over the topological space $\Spc(\cat T^c)$. The purpose of the topological condition on $\Spc(\cat T^c)$ is to guarantee that there exists a ``stalk'' category
$\kappaP \otimes \cat T = \Loco{\kappaP}$ for every $\cat P \in \Spc(\cat T^c)$ which captures the information of $\cat T$ supported on the point~$\cat P$. The Balmer--Favi support of an object $t \in \cat T$ is then defined as
	\[
		\Supp(t) \coloneqq \SET{ \cat P \in \Spc(\cat T^c) }{ \kappaP \otimes t \neq 0 },
	\]
consisting of the stalks where the object is nontrivial. Under our topological hypothesis, this notion of support extends the universal support for compact objects: $\Supp(x) = \supp(x)$ for $x \in \cat T^c$.

Following BIK, we say that~$\cat T$ is stratified if the following two conditions hold:
	\begin{itemize}
		\item (\emph{The local-to-global principle}) Any object $t\in\cat T$ can be reconstructed from its germs $\kappaP \otimes t$. In other words, $t \in \Loco{\kappaP \otimes t\mid \cat P \in \Spc(\cat T^c)}$.
		\item (\emph{Minimality at all points}) For each $\cat P \in \Spc(\cat T^c)$, the localizing $\otimes$-ideal $\kappaP \otimes \cat T=\Loco{\kappaP}$ is a minimal localizing $\otimes$-ideal. In other words, it is generated by any non-zero object it contains.
	\end{itemize}
These two conditions guarantee that all localizing $\otimes$-ideals in $\cat T$ can be built from the collection $\SET{\kappaP\otimes\cat T}{\cat P \in \Spc(\cat T^c)}$ and that each of these stalk $\otimes$-ideals cannot be refined further. It follows that the Balmer--Favi notion of support provides a bijection between the localizing $\otimes$-ideals of $\cat T$ and the subsets of $\Spc(\cat T^c)$. Our first contribution is to clarify that these conditions are actually equivalent to the classification of localizing $\otimes$-ideals:

\begin{thmx}\label{thmx:A}
	Let $\cat T$ be a rigidly-compactly generated tensor-triangulated category with $\Spc(\cat T^c)$ weakly noetherian. The following are equivalent:
	\begin{enumerate}
		\item The local-to-global principle holds for $\cat T$ and for each $\cat P \in \Spc(\cat T^c)$, $\kappaP \otimes \cat T = \Loco{\kappaP}$ is a minimal localizing $\otimes$-ideal of $\cat T$.
		\item The map
            \[ 
                \big\{ \text{localizing $\otimes$-ideals of $\cat T$} \big\} \xra{\Supp} \big\{ \text{subsets of $\Spc(\cat T^c)$}\big\}
            \]
			defined by $\cat L \mapsto \bigcup_{t \in \cat L} \Supp(t)$
            is a bijection.
	\end{enumerate}
\end{thmx}

For a more detailed statement see \cref{thm:equiv-strat}. This basic observation --- that stratification \emph{is} the classification of localizing $\otimes$-ideals --- does not appear to be in the literature. It is not particularly deep, but it does involve some set-theoretic subtleties related to the fact that \emph{a priori} we don't know that the localizing $\otimes$-ideals of $\cat T$ form a set. Another of our basic contributions is a strengthening of known results on the local-to-global principle. For example, we show in \cref{thm:local-to-global} that it always holds when the spectrum is noetherian:

\begin{thmx}\label{thmx:local-to-global}
	Let $\cat T$ be a rigidly-compactly generated tensor-triangulated category with $\Spc(\cat T^c)$ noetherian. Then $\cat T$ satisfies the local-to-global principle.
\end{thmx}

We now highlight three key features of our theory of stratification which will be expanded upon below:

	\begin{itemize}
		\item (\emph{Universality}) The Balmer--Favi notion of support provides the universal approach to stratification in weakly noetherian contexts. Any support theory for $\cat T$ lying in a weakly noetherian space that classifies localizing $\otimes$-ideals (in a way compatible with the usual classification of thick $\otimes$-ideals of compact objects) is equivalent to the Balmer--Favi notion of support. As a consequence, if $\cat T$ is stratified in the sense of BIK, then it is stratified in our sense (i.e., by the Balmer--Favi notion of support).
		\item (\emph{Permanence}) Stratification exhibits good permanence properties under base-change functors. In particular, it satisfies versions of Zariski and \'etale descent. This often allows classification problems to be reduced to more regular tt-categories for which stratification is known to hold.
		\item (\emph{Generality}) Fundamental categories in diverse areas of mathematics are stratified in our sense and thus admit a classification of localizing $\otimes$-ideals in terms of their Balmer spectra. This in particular applies to examples for which the BIK support is not applicable due to the lack of a suitable ring action, such as categories of local spectra, categories of equivariant spectra, or derived categories of non-affine schemes. 
	\end{itemize}

Our approach to understanding the global structure of a big tt-category divides into two parts: Firstly, it is a property of the category whether it is stratified by the Balmer--Favi support. If so, then the classification of localizing $\otimes$-ideals of~$\cat T$ is given in terms of the \emph{set} underlying the Balmer spectrum $\Spc(\cat T^c)$. Secondly, one can then try to determine the \emph{topology} of $\Spc(\cat T^c)$, which would result in the classification of thick $\otimes$-ideals of $\cat T^c$. This reverses the historically more common direction of first classifying thick $\otimes$-ideals of compact objects, and then attempting to control big objects. As a consequence, stratification is a rather flexible notion with good permanence properties even in cases where we have only partial knowledge of $\Spc(\cat T^c)$; we illustrate this point in \cref{rem:topology_ideals,rem:no-thick-classification}. Additionally, this two-step process distinguishes our notion of stratification from that of Benson--Iyengar--Krause: the latter simultaneously provides the classification of localizing $\otimes$-ideals \emph{and} determines the Balmer spectrum in terms of $\Spec(R)$. In this sense, BIK-stratification is both stronger and weaker than our notion of stratification.
\[ \ast \ast \ast \]

We will now discuss the above features of our theory of stratification in more detail and provide a coarse summary of our main results.

\subsection*{Universality of stratification}
\sloppy 
The Balmer--Favi notion of support provides the universal approach to stratification in weakly noetherian contexts:

\begin{thmx}\label{thmx:universal}
	Let $\cat T$ be a rigidly-compactly generated tt-category. Let $\sigma\colon \cat T \to \cat P(X)$ be a support function for $\cat T$ lying in a weakly noetherian space~$X$. If this notion of support stratifies $\cat T$ in a way compatible with the usual classification of the thick $\otimes$-ideals of $\cat T^c$, then there is a unique identification $(X,\sigma)\cong (\Spc(\cat T^c),\Supp)$ with the Balmer--Favi notion of support.
\end{thmx}

A more precise statement is \cref{thm:universal-support}. Note that \cref{thmx:universal} is a uniqueness result, rather than a universal property, but philosophically it means that if the localizing $\otimes$-ideals can be classified by a reasonable notion of support, then the Balmer--Favi notion of support will do the job; so the Balmer--Favi notion of support is the ``universal choice'' that should be considered whenever we are confronted with a new category. In any case, the following consequence is \cref{cor:BIK}:

\begin{thmx}\label{thm:intro_bik}
	Let $\cat T$ be a rigidly-compactly generated tensor-triangulated category which, in the terminology of \cite{BensonIyengarKrause11b}, is noetherian and stratified by the action of a graded-noetherian ring $R$. Then the BIK space of supports $\supp_R(\cat T)$ is canonically homeomorphic to~$\Spc(\cat T^c)$ and the BIK notion of support coincide with the Balmer--Favi notion of support.
\end{thmx}

In the companion paper \cite{bhs2}, we compare the Balmer--Favi notion of support with the homological support introduced by Balmer \cite{Balmer20_bigsupport}. In particular, we will prove that they coincide in a strong sense whenever $\cat T$ is stratified. 

\subsection*{Permanence of stratification}

A standard approach to stratify a given tt-category $\cat T$ is to find a good cover of $\cat T$, that is, a collection of tt-functors $F_i\colon \cat T \to \cat U_i$ with $\cat U_i$ stratified, and then verify that the stratification descends along the $F_i$. This technique is the key to the stratification of the stable module category of a finite group in the work of BIK \cite{BensonIyengarKrause11a}, who thereby reduce to the case of elementary abelian groups. 

Since our notion of stratification is based on the Balmer spectrum $\Spc(\cat T^c)$, which by virtue of its very construction behaves well under change of categories and satisfies descent properties, we are able to control the behaviour of stratification under base change. We isolate two general principles, which geometrically speaking amount to Zariski descent and a weak form of \'etale descent for stratification. In \cref{cor:cover}, we prove:

\begin{thmx}\label{thm:zariski_intro}
	Let $\cat T$ be a rigidly-compactly generated tensor-triangulated category with $\Spc(\cat T^c)$ weakly noetherian and satisfying the local-to-global principle. Suppose $\Spc(\cat T^c) = \bigcup_{i \in I} V_i$ is a cover by complements of Thomason subsets $V_i$. Then $\cat T$ is stratified if and only if each of the tt-categories $\cat T(V_i)$ is stratified.
\end{thmx}

In tt-geometry, the quotient map $\cat T \to \cat T(V_i)$ corresponds to a finite localization, and the theorem can be applied to any Zariski cover by quasi-compact open subsets. The theorem thus expresses a tt-geometric form of Zariski descent for stratification. It is somewhat more general, however, as complements of Thomason subsets include arbitrary intersections of quasi-compact open sets. Thus, \cref{thm:zariski_intro} allows not just restriction to a Zariski open cover (as has already been observed in \cite[Theorem~8.11]{Stevenson13}) but also reduces the problem of stratification to the \emph{local categories} $\cat T_{\cat P}\coloneqq \cat T/\Loco{\cat P}$ at each point $\cat P \in \Spc(\cat T^c)$. Combined with \cref{thmx:local-to-global}, we have a powerful method to reduce the problem of stratification to local tt-categories.

Likewise, there is an analogue of \'etale extensions in tt-geometry, introduced by Balmer \cite{Balmer15}. In \cref{thm:finite-etale}, we prove a form of \'etale descent for stratification:

\begin{thmx}\label{thm_etale_intro}
	Let $F\colon \cat C \to \cat D$ be a finite \'{e}tale morphism of rigidly-compactly generated tensor-triangulated categories. Assume that both categories have noetherian spectrum, and let
		\[
			\varphi\colon \Spc(\cat D^c) \to \Spc(\cat C^c)
		\]
	denote the induced map. If $\cat P \in \Spc(\cat D^c)$ is a point such that $\varphi^{-1}(\{\varphi(\cat P)\}) = \{\cat P\}$ then minimality at $\cat P$ in $\cat D$ implies minimality at $\varphi(\cat P)$ in $\cat C$.
\end{thmx}

This result will play an important role in our proof of stratification for categories of spectral Mackey functors (see \cref{thm:mackey_intro} below). We will return to the problem of proving more general kinds of descent in future work with Castellana.

These permanence properties allow us to establish stratification even in circumstances where we are unable to determine the Balmer spectrum. As noted before, the putative classification of localizing $\otimes$-ideals involves only the underlying set of~$\Spc(\cat T^c)$. This leads to the --- at first sight surprising --- observation that (in suitably noetherian situations) the classification of localizing $\otimes$-ideals in~$\cat T$ via stratification is in fact \emph{simpler} than the classification of thick $\otimes$-ideals of compact objects in~$\cat T$; indeed, the latter is tantamount to the computation of the topology on~$\Spc(\cat T^c)$. For the above reasons, we regard the fact that our approach to stratification is based on the Balmer spectrum as a strength rather than a weakness of the theory.

\subsection*{Generality of stratification}

We complement the formal development of the theory with an ample collection of examples from chromatic homotopy theory, equivariant homotopy theory, and algebraic geometry. A uniform perspective on the classification of localizing $\otimes$-ideals emerges: Our framework unifies the modular representation theoretic instances of stratification studied by BIK (via \cref{thm:intro_bik}) with known and new examples. The following theorem collects four examples of tt-categories for which the classification of localizing $\otimes$-ideals was already known, which we revisit with our techniques.

\begin{thmx}\label{thm:old-examples}
	The following rigidly-compactly generated tt-categories are stratified:
	\begin{enumerate}
		\item The stable module category $\StMod(kG)$ of a finite group $G$ over a field $k$. In this case, the Balmer spectrum is homeomorphic to $\Proj(H^*(G,k))$. This example is due to \cite[Theorem~10.3]{BensonIyengarKrause11a}. See \cref{exa:BIK-stratified}.
		\item The derived category $\Derqc(X)$ of a noetherian scheme $X$. In this case, the Balmer spectrum is homeomorphic to the underlying space of $X$. This example is due to \cite[Corollary~8.13]{Stevenson13}. The resulting classification of localizing $\otimes$-ideals is originally due to \cite[Corollary 4.13]{AlonsoJeremiasSouto04}. See \cref{cor:noetherian-schemes-strat}.	
		\item The category of $E(n)$-local spectra for any $n \ge 0$. In this case, the Balmer spectrum is homeomorphic to the poset $[0,n]$. The resulting classification of localizing $\otimes$-ideals is originally due to Hovey--Strickland \cite[Theorem~6.14]{HoveyStrickland99}. See \cref{thm:En-stratified}.
		\item The category of rational $G$-spectra for any compact Lie group $G$. In this case, the Balmer spectrum is homeomorphic to the conjugacy classes of closed subgroups of $G$ equipped with the zf-topology of Greenlees \cite[Section 10]{Greenlees19_rational}. The resulting classification of localizing $\otimes$-ideals is originally due to Greenlees \cite[Theorem 1.6]{Greenlees19_rational}. See \cref{thm:rational_stratification}.
	\end{enumerate}
\end{thmx}

Although the theorem could be regarded as just a repackaging of known results, it does contain new mathematical content. For example, it establishes that the notions of support used by the above authors to classify localizing $\otimes$-ideals coincide with the Balmer--Favi notion of support. For example, the Balmer--Favi support for the category of rational $G$-spectra coincides with the notion of support provided by geometric isotropy used by Greenlees. Similarly, the Balmer--Favi support for the category of $E(n)$-local spectra coincides with the notion of support used by Hovey--Strickland defined using Morava $K$-theories.

We remark that none of these examples are canonically stratified in the sense of BIK. The reader might be surprised by this statement since Example (a) is the paragon example of BIK stratification. However, the endomorphism ring of the unit in the stable module category --- given by Tate cohomology --- is not noetherian, hence one cannot apply the BIK machinery to the canonical action. In order to classify the localizing $\otimes$-ideals, BIK has to let the ordinary group cohomology ring (which is noetherian) act on the category. \Cref{thm:old-examples} thus provides a more canonical classification of the localizing $\otimes$-ideals of the stable module category.

We then turn to a new class of examples. For any finite group $G$, the equivariant stable homotopy category $\SH_G$ can be identified with the category of spectral Mackey functors (in the sense of Barkwick~\cite{Barwick17}) valued in spectra. Similarly, the category of derived Mackey functors (in the sense of Kaledin \cite{Kaledin11}) is the category of spectral Mackey functors valued in $\HZ$-modules \cite[Section~4]{PatchkoriaSandersWimmer22}. More generally, we may study the category of spectral Mackey functors with coefficients in a commutative ring spectrum~$\bbE$. This category can in turn be identified with the derived category $\DEG$ of a certain equivariant ring spectrum $\bbE_G$ obtained by giving $\bbE$ the ``trivial'' $G$-action.

Building on the work of \cite{PatchkoriaSandersWimmer22}, we obtain a complete description of $\Spc(\DEG^c)$, as a set, for all finite groups $G$, in terms of the spectrum of the nonequivariant category of $\bbE$-modules $\Spc(\Der(\bbE)^c)$. See \cref{thm:spcDEG}. Although we obtain some results concerning the topology of $\Spc(\DEG^c)$, a full description of the topology is presently out of reach. Nevertheless, in \cref{thm:mackeystrat} we prove:

\begin{thmx}\label{thm:mackey_intro}
	Let $G$ be a finite group and let $\bbE \in \CAlg(\Sp)$ be a commutative ring spectrum such that $\Spc(\DE^c)$ is noetherian. If $\DE$ is stratified, then so is the category $\DEG$ of spectral $G$-Mackey functors with coefficients in $\bbE$.
\end{thmx}

The proof of this theorem demonstrates the power of our geometric theory of stratification based on the Balmer spectrum. \Cref{thm:zariski_intro,thm_etale_intro} allow us to pass through finite \'{e}tale extensions (restriction to a subgroup) and then to the relevant local categories, which we can identify via geometric fixed point functors. As a special case, the theorem applies to $\bbE=\HZ$ and we thereby obtain a classification of the localizing $\otimes$-ideals of the category of derived $G$-Mackey functors (\cref{cor:dermackstratified}). In this example, the Balmer spectrum $\Spc(\DHZG^c)$ was computed in \cite{PatchkoriaSandersWimmer22} and captures precisely the height 0 and height~$\infty$ parts of the spectrum of the equivariant stable homotopy category.

\Cref{thm:mackey_intro} also applies to the category of $E(n)$-local spectral Mackey functors by taking $\bbE=L_nS^0$ to be the $E(n)$-local sphere spectrum and invoking part (c) of \cref{thm:old-examples}. In this example, we do not have a complete description of the topology of the Balmer spectrum. We establish a continuous bijection onto the height $\le n$ part of the spectrum of the $G$-equivariant stable homotopy category (which is known for some, but not all, finite groups \cite{BalmerSanders17,BHNNNS19,KuhnLloyd20pp}) but we are not able to prove that this continuous bijection is a homeomorphism onto its image. We leave this as \cref{conj:homeomorphism}. It follows that we do not obtain a classification of the thick $\otimes$-ideals of compact objects for the category of $E(n)$-local spectral $G$-Mackey functors. Nevertheless, since we compute the underlying set of the Balmer spectrum, we still obtain a complete classification of localizing $\otimes$-ideals (\cref{cor:En-local-mack-strat}). This result is particularly striking because the classification of localizing $\otimes$-ideals for the $G$-equivariant stable homotopy category is at present completely open (even for $G=1$) and one does not expect such a simple answer.

\subsection*{Consequences of stratification}
One of the consequences of being stratified in the sense of BIK is that a form of the telescope conjecture holds: Every smashing localization is a finite localization. This connection between stratification and the telescope conjecture has also been studied by Stevenson \cite[Theorem 7.15]{Stevenson13}, who showed that, in the language of this paper, if $\cat T$ is a stratified rigidly-compactly generated category which occurs as the homotopy category of a monoidal model category and $\Spc(\cat T^c)$ is noetherian, then $\cat T$ satisfies the telescope conjecture. Based on Stevenson's work, we extend this to a class of spaces we call generically noetherian (\cref{defn:gennoetherian}). These lie in between noetherian spaces and weakly noetherian spaces, and also include all profinite spaces. In \cref{thm:gennoethtelescopeconj}, we prove:

\begin{thmx}\label{thm:intro_telescope}
	Let $\cat T$ be a stratified rigidly-compactly generated $tt$-category with generically noetherian spectrum. The telescope conjecture holds for $\cat T$. 
\end{thmx}

In a forthcoming paper \cite{BalchinBarthelGreenlees21ip}, it will be shown that for any compact Lie group $G$ the Balmer spectrum of the category of rational $G$-spectra is generically noetherian. Together with Example $(d)$ of \cref{thm:old-examples}, this shows that the category of rational $G$-spectra satisfies the telescope conjecture.

We also establish some additional consequences of stratification, such as a tensor-product formula for the Balmer--Favi support (\cref{thm:tensor_product}) and a complete description of the Bousfield lattice (\cref{thm:bousfield-lattice}). In particular, for a stratified category, there is no difference between localizing $\otimes$-ideals and Bousfield classes.

Finally, we turn to the original telescope conjecture for the stable homotopy category of spectra. Using the stratification of the $E(n)$-local category (Example~(c) of \cref{thm:old-examples}), we are able to reformulate the classical telescope conjecture in terms of stratification and the Balmer--Favi notion of support. See \cref{cor:telescope-strat} and \cref{cor:support-approach-to-tc}.

\subsection*{Outline of the document}
The paper consists of five parts. 

In \cref{part:I}, we review the basic theory of tensor triangular geometry, including the Balmer spectrum and finite localizations (\cref{sec:TTG}). We then introduce the Balmer--Favi notion of support in the context of a rigidly-compactly generated category whose spectrum is weakly noetherian (\cref{sec:BF-support}). We study the local-to-global principle in \cref{sec:LTG}, showing in particular that it can be checked on an open cover and that it holds automatically for categories with noetherian spectrum (\cref{thmx:local-to-global}).

In \cref{part:stratification1} we begin our study of stratification. We define the concept and establish conditions and criteria for when it holds (such as \cref{thmx:A}) in \cref{sec:strat}. We establish the forms of Zariski and \'etale descent for stratification mentioned above (\cref{thm:zariski_intro,thm_etale_intro}) in \cref{sec:reducing-to-local} and \cref{sec:etale}, respectively.

In \cref{part:consequences}, we first consider the universal nature of the Balmer--Favi support (\cref{sec:universal}) and prove \cref{thmx:universal,thm:intro_bik}. We then move on to the consequences of a category being stratified. For example, we prove a tensor-product theorem for support and determine the Bousfield lattice in \cref{sec:tensor-and-bousfield}. We then consider the telescope conjecture in \cref{sec:gentelescope} and prove \cref{thm:intro_telescope}. 

In \cref{part:chromatic}, we consider two examples, the $E(n)$-local stable homotopy category (\cref{sec:En-local}) and the rational $G$-equivariant stable homotopy category (\cref{sec:rational-G-spectra}) whose localizing $\otimes$-ideals have already been classified in the literature. We reinterpret these results as establishing that these categories are stratified. In \cref{sec:telescopeconjecture}, we apply the stratification of the $E(n)$-local category to reformulate the classical telescope conjecture in chromatic homotopy theory.

In \cref{part:spectral-mackey}, we study stratification for categories of spectral Mackey functors. We study the Balmer spectrum of these categories in \cref{sec:spectrum-E-mackey}, and discuss different models for these categories in \cref{sec:equivalence-spectral-Mackey}. We then prove \cref{thm:mackey_intro} and its corollaries in \cref{sec:strat-for-spectral-mack}.

\subsection*{Acknowledgements}
We thank Scott Balchin, Paul Balmer, and Greg Stevenson for useful conversations.

\renewcommand{\thepart}{\Roman{part}}

\part{Support and the local-to-global principle}\label{part:I}

We start by laying the foundations for the approach to stratification developed (in the more general relative context) by Stevenson \cite{Stevenson13}, which uses the Balmer spectrum \cite{Balmer05a} and the Balmer--Favi notion of big support \cite{BalmerFavi11}. We will develop these ideas in the context of a rigidly-compactly generated tensor-triangulated category whose Balmer spectrum is a ``weakly noetherian'' space. This ensures that every prime is ``weakly visible'' which allows us to consider the Balmer--Favi support at every point. We will investigate the basic properties of this notion of support, consider the local-to-global principle, and pave the way for the discussion of stratification in \cref{part:stratification1}. Although we use the Balmer--Favi notion of support throughout, the theory also depends crucially on the contributions of Benson--Iyengar--Krause \cite{BensonIyengarKrause08, BensonIyengarKrause11b, BensonIyengarKrause11a} and Hovey--Palmieri--Strickland \cite{HoveyPalmieriStrickland97}.

\section{Tensor triangular geometry}\label{sec:TTG}

We begin by reviewing some basic facts from tensor triangular geometry. For a more detailed discussion, we refer the reader to \cite{Balmer05a, Balmer10b, BalmerICM}.

\begin{Ter}\label{ter:tt-category}
    By a \emph{tensor-triangulated category} we mean a triangulated category equipped with a compatible closed symmetric monoidal structure in the sense of \cite[App.~A]{HoveyPalmieriStrickland97}.  Such a category is \emph{rigidly-compactly generated} if it is compactly generated as a triangulated category and if its compact objects coincide with its rigid objects (a.k.a.~dualizable objects).  In particular, its unit object $\unit$ is compact.  A rigidly-compactly generated tensor-triangulated category is precisely the same thing as a ``unital algebraic stable homotopy category'' in the language of \cite{HoveyPalmieriStrickland97}.  Sometimes we'll drop the ``tensor-triangulated'' and just speak of rigidly-compactly generated categories.  By a \emph{tensor-triangulated functor}, we mean a triangulated functor which is a \emph{strong} monoidal functor.
	We sometimes abbreviate ``tensor-triangulated'' by ``tt''.
\end{Ter}

\begin{Exa}\label{ex:rcg-cats}
    The stable homotopy category of spectra $\SH$ and the derived category of a commutative ring $\Der(R)$ are basic examples of rigidly-compactly generated tensor-triangulated categories. Further examples are mentioned in  \cite[Example~1.2.3]{HoveyPalmieriStrickland97}. Any smashing localization of a rigidly-compactly generated category is again rigidly-compactly generated (see \cite[Section~3.3]{HoveyPalmieriStrickland97}). For example, the stable homotopy category of $E(n)$-local spectra is rigidly-compactly generated. This is not true of an arbitrary Bousfield localization. For example,  the stable homotopy category of $K(n)$-local spectra is not rigidly-compactly generated if $n>0$ (as the unit $\unit$ is dualizable but not compact).
\end{Exa}

\begin{Rem}\label{rem:suffice-tensor-compact}
    If $\cat T$ is rigidly-compactly generated, then to check that a localizing subcategory $\cat L \subset \cat T$ is a $\otimes$-ideal ($\cat L \otimes \cat T \subseteq \cat L$) it suffices to check that it is closed under tensoring with compact objects ($\cat L\otimes \cat T^c \subseteq \cat L$) or even just a set of compact generators. (This is a straightforward exercise.) In particular, if~$\cat T$ is monogenic, meaning that it is generated by the unit $\unit$, then every localizing subcategory is a $\otimes$-ideal. Nevertheless, even in monogenic examples, we will favour language that is appropriate in more general situations and speak, for example, of classifying the localizing $\otimes$-ideals.
\end{Rem}

\begin{Not}
    The localizing $\otimes$-ideal generated by a class of objects $\cat E$ will be denoted $\Loco{\cat E}$. Sometimes we'll simply write $\langle \cat E \rangle \coloneqq \Loco{\cat E}$.
\end{Not}

\begin{Rem}\label{rem:compactly-generated}
    A localizing $\otimes$-ideal $\cat L$ in a rigidly-compactly generated category $\cat T$ is compactly generated as a localizing $\otimes$-ideal if and only if it is compactly generated as a localizing subcategory. (This is a straightforward exercise. The key is to establish that the localizing subcategory $\Loc\langle \cat L \cap \cat T^c\rangle$ is a $\otimes$-ideal.) Thus, there is no potential ambiguity when we speak of ``compactly generated'' localizing $\otimes$-ideals.
\end{Rem}

\begin{Rem}\label{rem:Balmer_spectrum}
    The Balmer spectrum \cite{Balmer05a} of an essentially small tt-category $\cat K$ is a topological space $\Spc(\cat K)$ whose points are the prime thick $\otimes$-ideals of $\cat K$.  Every object $x \in \cat K$ has an associated closed subset
        \begin{equation}\label{eq:supp}
            \supp(x) \coloneqq \SET{\cat P \subset \cat K}{x \not\in \cat P} \subseteq \Spc(\cat K),
        \end{equation}
    and these sets form a basis of closed sets for the topology on $\Spc(\cat K)$. Conceptually, $\Spc(\cat K)$ is the universal space equipped with a ``good'' notion of closed support for each object of~$\cat K$ (see \cite[Theorem~3.2]{Balmer05a}).  If $\cat K$ is rigid (meaning that each object is dualizable), then the thick $\otimes$-ideals of~$\cat K$ are in one-to-one correspondence with the Thomason subsets of $\Spc(\cat K)$ --- the unions of closed sets, each of which has quasi-compact complement.  The bijection sends a thick $\otimes$-ideal $\cat I \subseteq \cat K$ to the Thomason subset $\bigcup_{x \in \cat I} \supp(x)$, while a Thomason subset $Y \subseteq \Spc(\cat K)$ is sent to the  thick $\otimes$-ideal $\cat K_Y \coloneqq \SET{x \in \cat K}{\supp(x) \subseteq Y}$.  This is the abstract Thick Subcategory Classification Theorem (see \cite[Theorem~4.10 and Remark~4.3]{Balmer05a}) which translates the problem of classifying the thick $\otimes$-ideals of~$\cat K$ to the problem of understanding the Balmer spectrum $\Spc(\cat K)$ and its universal notion of support.
\end{Rem}

\begin{Exa}
    If $\cat T$ is a rigidly-compactly generated tensor-triangulated category then its subcategory $\cat T^c$ of compact (=rigid) objects is an essentially small rigid tensor-triangulated category and we can consider its spectrum $\Spc(\cat T^c)$.
\end{Exa}

\begin{Rem}\label{rem:closure-of-prime}
    If $\cat P \in \Spc(\cat K)$ then 
        \[
            \overbar{\{ \cat P \}} = \SET{\cat Q \in \Spc(\cat K)}{ \cat Q \subseteq \cat P}.
        \]
    Thus, the closed points of $\Spc(\cat K)$ are the \emph{minimal} prime $\otimes$-ideals of $\cat K$.
\end{Rem}

\begin{Rem}
    An essentially small rigid tensor-triangulated category $\cat K$ is said to be \emph{local} if $\Spc(\cat K)$ has a unique closed point.  This is equivalent to the zero $\otimes$-ideal~$(0)$ being prime, in which case it is the unique closed point. See \cite[Section 4]{Balmer10b} for further discussion.
\end{Rem}

\begin{Ter}\label{Ter:local}
    We'll say that a rigidly-compactly generated tensor-triangulated category $\cat T$ is \emph{local} if $\cat T^c$ is local, that is, if $\Spc(\cat T^c)$ has a unique closed point.
\end{Ter}

\begin{Rem}\label{rem:quasi-separated}
    The Balmer spectrum $\Spc(\cat K)$ is a spectral space in the sense of Hochster~\cite{Hochster69,DickmannSchwartzTressl19}. A spectral space is quasi-separated by definition, meaning that the intersection of any two quasi-compact open subsets is quasi-compact. However, a stronger statement is true: an arbitrary intersection of quasi-compact open subsets is again quasi-compact. Equivalently, the complement of a Thomason subset is always quasi-compact. Indeed, the complement of a Thomason subset forms a spectral subspace (see the proof of \cite[Lemma~3.2]{Stevenson17}, for example) and thus is quasi-compact.
\end{Rem}

\begin{Rem}\label{rem:open-in-dual}
    The Thomason subsets of a spectral space $X$ are precisely the open subsets of the Hochster-dual spectral topology on $X$ (see \cite[Proposition 8]{Hochster69}). We will denote the dual topological space by $X^*$.
\end{Rem}

\begin{Rem}\label{rem:thomason-closed}
    The Thomason \emph{closed} subsets of $\Spc(\cat K)$ can be characterized as the closed subsets whose complement is quasi-compact and these are precisely the subsets of the form $\supp(x)$ for some $x \in \cat K$; see \cite[Lemma 3.3]{Sanders13}, for example.
\end{Rem}

\begin{Rem}
    The spectrum is functorial: Any tt-functor $F\colon {\cat K \to \cat L}$ between essentially small tt-categories induces a continuous map $\varphi\colon \Spc(\cat L)\to \Spc(\cat K)$ given by $\varphi(\cat P) \coloneqq F^{-1}(\cat P)$. It follows from the definitions that $\varphi^{-1}(\supp_{\cat K}(x)) = \supp_{\cat L}(F(x))$ for every object $x \in \cat K$. We then see from \cref{rem:thomason-closed} that the map $\varphi$ is spectral, meaning that the preimage of a quasi-compact open subset is quasi-compact. Hence the preimage of a Thomason subset is Thomason. In other words, the continuous map $\varphi$ is also continuous with respect to the dual spectral topologies.
\end{Rem}

\begin{Rem}\label{rem:hochster-dual}
    Every Thomason closed subset of a spectral space $X$ is a \emph{quasi-compact} open in the Hochster-dual topology on $X$ (but not conversely, as the example $(0) \in \Spc(\SHcp)$ shows). Indeed, if $A \subseteq X$ is Thomason then $A^* \subseteq X^*$ is open. On the other hand, if $A \subseteq X$ is closed then $A^* \subseteq X^*$ is a subset such that $A^{**} \subseteq X^{**}$ is closed (as $A^{**} = A$ and $X^{**} = X$). Hence $A^{*} \subseteq X^{*}$ is the complement of a Thomason, hence quasi-compact by \cref{rem:quasi-separated}. In particular, if $Z \subseteq \Spc(\cat K)$ is a Thomason closed subset then any cover $Z=\bigcup_{i \in I} Y_i$ by Thomason subsets can be reduced to a finite subcover, since in the Hochster-dual topology this is just a cover of a quasi-compact open by a family of opens.
\end{Rem}

\begin{Rem}\label{rem:determined-by-inclusions}
    Suppose $\Spc(\cat K)$ satisfies the topological condition that every closed set can be written as a finite union of irreducible closed sets:
        \begin{equation}\label{eq:topologydetbyspecialization}
            \forall A \subseteq{\Spc(\cat K)} \text{ closed} \colon A = \bigcup_{i=1}^n \overbar{\{\cat P_i\}} \text{ for finitely many } \cat P_1,\ldots,\cat P_n \in \Spc(\cat K).
        \end{equation}
    Equivalently, suppose every closed set is the closure of a finite set of points. If this is the case, then the topology on $\Spc(\cat K)$ is determined by the specialization relation $\cat Q \in \overbar{\{\cat P\}}$. By \cref{rem:closure-of-prime}, this can be reformulated as follows: If $\Spc(\cat K)$ satisfies~\eqref{eq:topologydetbyspecialization}, then the topology on $\Spc(\cat K)$ is determined by the inclusions among primes $\cat Q \subseteq \cat P$.
\end{Rem}

\begin{Exa}
    If the space $\Spc(\cat K)$ is noetherian, then Property~\eqref{eq:topologydetbyspecialization} holds; see \cite[Theorem~8.1.11]{DickmannSchwartzTressl19} or \cite[Proposition~2.38]{PatchkoriaSandersWimmer22}, for instance. A non-noetherian example is the Balmer spectrum $\Spc(\SH_G^c)$ of the category of compact $G$-spectra for~$G$ a finite group; see \cite[Proposition~6.1]{BalmerSanders17}. 
	There are also naturally occurring non-examples. Indeed, a profinite Balmer spectrum satisfies Property~\eqref{eq:topologydetbyspecialization} if and only if it is finite. An explicit non-example is then the Balmer spectrum $\Spc(\SH_{O(2)}^c)$ of the category of compact $O(2)$-spectra; see~\cite{bgh_balmer}.
\end{Exa}

\begin{Rem}
    In general, the topology on $\Spc(\cat K)$ is determined by the inclusions among primes together with the underlying constructible topology on $\Spc(\cat K)$. This is the Priestley space point of view on the Balmer spectrum, which will be studied in more detail in a forthcoming paper \cite{BalchinBarthelGreenlees21ip}. See also \cite[Section~1.5]{DickmannSchwartzTressl19}.
\end{Rem}

\begin{Rem}\label{rem:set-of-generalizations}
    For a point $\cat P$ in a spectral space $X$, we will write
        \[
            \gen(\cat P) = \gen_X(\cat P) = \SET{\cat Q \in X}{\cat P \in \overbar{\{\cat Q\}}}
        \]
    for the generalization closure of $\cat P$ in $X$. We will also write $Y_{\cat P}$ for the union of all Thomason closed subsets which do not contain~$\cat P$. If $X = \Spc(\cat K)$, then
        \[
            Y_{\cat P} = \bigcup_{x \in \cat P}\supp(x)
        \]
    by Remark~\ref{rem:thomason-closed}. It is the Thomason subset whose corresponding thick $\otimes$-ideal is the prime $\cat P$ itself: $\cat K_{Y_{\cat P}} = \cat P$. The complement $\Spc(\cat K) \setminus Y_{\cat P} = \SET{\cat Q}{\cat P \subseteq \cat Q}=\SET{\cat Q}{\cat P \in \overbar{\{\cat Q\}}} = \gen(\cat P)$ consists of all generalizations of~$\cat P$.
\end{Rem}

\begin{Rem}[Smashing Bousfield localizations]\label{rem:smashing_idempotents}
    Let $\cat T$ be a rigidly-compactly generated tensor-triangulated category.  A Bousfield localization $L\colon\cat T \to \cat T$ is said to be \emph{smashing} if the functor $L$ preserves coproducts.  This is equivalent to the natural map $L\unit \otimes X \to LX$ being an isomorphism (see \cite[Section 3.3]{HoveyPalmieriStrickland97} or \cite[Section 5.5]{Krause10}, for example). A smashing localization $L\colon\cat T \to \cat T$ is completely determined by its kernel $\Ker(L)$ which is a localizing $\otimes$-ideal of $\cat T$.  We can thus think of the smashing localizations of $\cat T$ either as a certain class of Bousfield localizations, or as a certain class of localizing $\otimes$-ideals of $\cat T$ (the \emph{smashing $\otimes$-ideals}).  A third perspective is explained in \cite{BalmerFavi11}: The smashing localizations of~$\cat T$ correspond to the \emph{idempotent triangles} in $\cat T$; that is, exact triangles
        \[
            e \to \unit \to f \to \Sigma e
        \]
    with the property that $e \otimes f=0$.  It follows that the objects $e$ and $f$ are tensor-idempotents ($e\otimes e \simeq e$ and $f \otimes f \simeq f$) and that the functor $f \otimes - \colon \cat T \to \cat T$ is a smashing localization.  Its corresponding smashing $\otimes$-ideal is $\Ker(f\otimes -) = e \otimes \cat T= \Loco{e}$ and the subcategory of local objects is $\Loco{e}^\perp = f \otimes \cat T = \Loco{f} \cong \cat T/\Loco{e}$.  See \cite{BalmerFavi11} for a more detailed discussion.
\end{Rem}

\begin{Rem}[Finite localizations]\label{rem:finite-localizations}
    Let $Y \subseteq \Spc(\cat T^c)$ be a Thomason subset, with corresponding thick $\otimes$-ideal $\cat T^c_Y = \SET{ x\in \cat T^c}{\supp(x) \subseteq Y}$, and let $V \coloneqq {\Spc(\cat T^c) \setminus Y}$ denote the complement. There is an associated idempotent triangle (a.k.a.~smashing localization)
        \[ 
            e_Y \to \unit \to f_Y \to \Sigma e_Y
        \]
    in $\cat T$ such that $\Ker(f_Y \otimes -)=e_Y \otimes \cat T = \Loco{e_Y} = \Loc\langle \cat T^c_Y\rangle$. Define 
        \[
            \cat T(V) := \cat T/\Loco{e_Y} \cong \cat T / \Loc\langle \cat T^c_Y\rangle \cong f_Y \otimes \cat T 
        \]
    to be the associated localization. It is again rigidly-compactly generated. The localization functor $\cat T \to \cat T(V)$ preserves compact objects and hence induces a continuous map $\Spc(\cat T(V)^c) \to \Spc(\cat T^c)$. This map induces a homeomorphism
        \begin{equation}\label{eq:SpcTV}
            \Spc(\cat T(V)^c) \cong V \hookrightarrow \Spc(\cat T^c)
        \end{equation}
    which identifies $\Spc(\cat T(V)^c)$ with $V \subseteq \Spc(\cat T^c)$. This is essentially the content of the Neeman--Thomason Localization Theorem \cite[Theorem~2.1]{Neeman92b} and explains the choice of notation $\cat T(V)$. For further discussion of these finite localizations see~\cite[Section 4]{BalmerFavi11}, \cite[Theorem~3.3.3]{HoveyPalmieriStrickland97}, \cite[Section 5]{BalmerSanders17}, and \cite{Miller92}. Note that these finite localizations of $\cat T$ (which correspond to the thick $\otimes$-ideals of $\cat T^c$) are precisely those smashing localizations of $\cat T$ whose kernel is compactly generated (recall \cref{rem:compactly-generated}).
\end{Rem}

\begin{Def}[The local category at a point]\label{def:primeverdierquotient}\label{exa:local-categories}
    Let $\cat P \in \Spc(\cat T^c)$ be a point in the Balmer spectrum and consider the associated Thomason subset $Y_{\cat P}$ from \cref{rem:set-of-generalizations}. The corresponding finite localization $\cat T_{\cat P} := \cat T/\langle \cat P\rangle$ is $\cat T(\gen(\cat P))$. This category is local in the sense of \cref{Ter:local}. (It has a unique closed point.) We call it the local category of $\cat T$ at~$\cat P$.
\end{Def}

\begin{Rem}\label{rem:union-of-thomason}
    If $Y=\bigcup_{i \in I} Y_i$ is a union of Thomason subsets then $\Loco{e_Y} = \Loco{e_{Y_i} \mid i\in I}$. Indeed, the $\supseteq$ inclusion is immediate. For the converse just note that $\cat T^c_Y = \thickt{\cat T^c_{Y_i} \mid i \in I}$ by the classification of thick $\otimes$-ideals.
\end{Rem}

\begin{Lem}\label{lem:Y1-and-Y2}
    For any two Thomason subsets $Y_1, Y_2 \subseteq \Spc(\cat T^c)$, we have
        \begin{enumerate}
            \item $e_{Y_1} \otimes f_{Y_2} = 0$ if and only if $Y_1 \subseteq Y_2$;
            \item $e_{Y_1} \otimes e_{Y_2} = e_{Y_1 \cap Y_2}$;
            \item $f_{Y_1} \otimes f_{Y_2} = f_{Y_1 \cup Y_2}$.
        \end{enumerate}
\end{Lem}

\begin{proof}
    For part (a) just note that
        $e_{Y_1} \otimes f_{Y_2} = 0$
            $\Leftrightarrow$
        $\cat T_{Y_1}^c \otimes f_{Y_2} = 0$
            $\Leftrightarrow$
        $\cat T_{Y_1}^c \subseteq \Loc\langle \cat T_{Y_2}^c\rangle \cap \cat T^c = \cat T_{Y_2}^c$
            (by \cite[Lem.~2.2]{Neeman92b})
            $\Leftrightarrow$ 
        $Y_1 \subseteq Y_2$.
    For part (b), note that~(a) tells us that if $Y_1 \subseteq Y_2$ then $e_{Y_1} \simeq e_{Y_1} \otimes e_{Y_2}$. Applying this observation twice we have $e_{Y_1 \cap Y_2} \simeq e_{Y_1 \cap Y_2} \otimes e_{Y_1} \otimes e_{Y_2}$. We claim that $f_{Y_1 \cap Y_2} \otimes e_{Y_1} \otimes e_{Y_2} = 0$. One observes that $\cat T^c_{Y_1 \cap Y_2} = \cat T^c_{Y_1} \cap \cat T^c_{Y_2} = \thickt{\cat T^c_{Y_1} \otimes \cat T^c_{Y_2}}$ for example by invoking the classification of thick tensor-ideals of $\cat T^c$. Then use $e_{Y_1} \otimes e_{Y_2} \in \Loco{\cat T^c_{Y_1} \otimes \cat T^c_{Y_2}}$. For part (c), writing $f\coloneqq f_{Y_1 \cup Y_2}$, we have $e_{Y_1} \otimes f =0$ and $e_{Y_2} \otimes f=0$, from which it follows that $f \simeq f \otimes f_{Y_1} \otimes f_{Y_2}$. It remains to show that $e_{Y_1 \cup Y_2} \otimes f_{Y_1} \otimes f_{Y_2}=0$ so that $f_{Y_1} \otimes f_{Y_2} \simeq f \otimes f_{Y_1} \otimes f_{Y_2}$. For this just note that $\thickt{\cat T^c_{Y_1}, \cat T^c_{Y_2}} = \cat T^c_{Y_1 \cup Y_2}$ by the classification of thick $\otimes$-ideals. Hence $\Loc\langle\cat T^c_{Y_1 \cup Y_2}\rangle = \Loc\langle\cat T^c_{Y_1},\cat T^c_{Y_2}\rangle$ so that indeed $e_{Y_1 \cup Y_2} \otimes f_{Y_1} \otimes f_{Y_2} = 0$.
\end{proof}

\begin{Rem}\label{rem:geometric-functor}
    A coproduct-preserving tensor-triangulated functor $F\colon\cat C \to \cat D$ between rigidly-compactly generated tensor-triangulated categories enjoys a number of nice properties, as explained in \cite{BalmerDellAmbrogioSanders16}. In particular, it has a right adjoint $U\colon\cat D \to \cat C$ and we have the projection formula $U(F(x)\otimes y) \simeq x \otimes U(y)$ for all $x \in \cat C$ and $y \in \cat D$. Another nice property is that $F$ preserves dualizable (=compact) objects and hence restricts to a tensor triangulated functor $\cat C^c \to \cat D^c$. This in turn induces a continuous map $\varphi\colon\Spc(\cat D^c) \to \Spc(\cat C^c)$. It is proved in \cite[Proposition~5.11]{BalmerSanders17} that if $e_Y \to \unit \to f_Y \to \Sigma e_Y$ is the idempotent triangle in $\cat C$ associated to a Thomason subset $Y \subseteq \Spc(\cat C^c)$, then its image $F(e_Y) \to \unit \to F(f_Y) \to \Sigma F(e_Y)$ is the idempotent triangle in $\cat D$ associated to the Thomason subset $\varphi^{-1}(Y) \subseteq \Spc(\cat D^c)$. That is:
        \begin{equation}\label{eq:FeY}
             F(e_Y) \simeq e_{\varphi^{-1}(Y)} \quad\text{ and }\quad F(f_Y) \simeq f_{\varphi^{-1}(Y)}.
        \end{equation}
    The following is an immediate consequence:
\end{Rem}

\begin{Prop}\label{prop:restrict-in-target}
    Let $F\colon\cat C \to \cat D$ be a coproduct-preserving tensor-triangulated functor between rigidly-compactly generated tensor-triangulated categories and let
        \[ 
            \varphi\colon\Spc(\cat D^c) \to \Spc(\cat C^c)
        \]
    be the induced map on spectra. For any Thomason subset $Y \subseteq \Spc(\cat C^c)$ ($V\coloneqq \Spc(\cat C^c)\setminus Y$), we have an induced coproduct-preserving tensor-triangulated functor
        \[
            \cat C(V) \to \cat D(\varphi^{-1}(V))
        \]
    which on spectra is the restriction
        \[ 
            \varphi|_V \colon \varphi^{-1}(V) \to V.
        \]
\end{Prop}

\begin{proof}
    This follows from \cref{rem:finite-localizations}, \cref{rem:geometric-functor} and the definitions (in particular, \eqref{eq:SpcTV} and \eqref{eq:FeY}).
\end{proof}

\begin{Rem}\label{rem:finite-of-finite}
    Applied to a finite localization $F\colon\cat C \to \cat C(W)\eqqcolon \cat D$, the induced functor $\cat C(V) \to \cat C(W)(\varphi^{-1}(V)) \cong \cat C(V \cap W)$ is itself a finite localization, namely the finite localization associated to the Thomason subset $V \cap W^c \subset V$. Just recall from \eqref{lem:Y1-and-Y2} that $f_{V^c} \otimes f_{W^c} \simeq f_{V^c \cup W^c}$. It follows that we obtain a commutative diagram
        \[
            \begin{tikzcd}
            \cat C \ar[d] \ar[r,"F"] & \cat C(W) \ar[d] \\
            \cat C(V) \ar[r] & \cat C(V \cap W),
            \end{tikzcd}
        \]
    in which every functor is a finite localization.
\end{Rem}

\begin{Prop}\label{prop:localfinitelocalization}
    Let $V \subseteq \Spc(\cat C^c)$ be the complement of a Thomason subset. For any $\cat P \in \Spc(\cat C(V)^c)\cong V$, the localization $F\colon\cat C \to \cat C(V)$ induces an equivalence
        \[
            \cat C/\langle \varphi(\cat P)\rangle \cong \cat C(V)/\langle \cat P \rangle.
        \]
    Here $\varphi\colon\Spc(\cat C(V)^c) \cong V \hookrightarrow \Spc(\cat C^c)$ denotes the map on spectra induced by $F$.
\end{Prop}

\begin{proof}
    Let $W\coloneqq \gen(\varphi(\cat P)) \subseteq \Spc(\cat C^c)$. Since $V$ is the complement of a Thomason subset, it is generalization closed. Hence $\varphi^{-1}(W) = \varphi^{-1}(\gen(\varphi(\cat P))) = \gen(\cat P)$. Thus, applying \cref{prop:restrict-in-target} to the functor $F$ and the localization $\cat C \to \cat C(W)$, we obtain 
        \[
            \begin{tikzcd}
             \cat C \ar[d] \ar[r,"F"] & \cat C(V) \ar[d]\\
             \cat C(W) \ar[r,"F|_W"] & \cat C(V)(\varphi^{-1}(W))
            \end{tikzcd}
        \]
    and, recalling \cref{def:primeverdierquotient}, the bottom functor is 
        \[
            \cat C/\langle \varphi(\cat P)\rangle = \cat C(W) \to \cat C(V)(\varphi^{-1}(W)) = \cat C(V)/\langle \cat P \rangle.
        \]
    Moreover, by \cref{rem:finite-of-finite}, this can be identified with the localization
        \[
            \cat C(W) \to \cat C(V \cap W)
        \]
    and this is an equivalence since $V \cap W = W$.
\end{proof}

\begin{Rem}\label{rem:local-cat-in-open}
    In less pedantic language, \cref{prop:localfinitelocalization} asserts that the local category of $\cat C(V)$ at $\cat P \in V \subseteq \Spc(\cat C^c)$ is the same as the local category of $\cat C$ at~$\cat P$. For example, the local category at $\cat P$ can be computed in any quasi-compact open neighbourhood of $\cat P$.
\end{Rem}

\begin{Exa}[Algebraic localization]\label{exa:algebraic-localizations}
    In any tt-category $\cat T$, the endomorphism ring $\End_{\cat T}(\unit)$ is commutative and there is a natural continuous map
        \[
            \rho \colon \Spc(\cat T^c) \to \Spec(\End_{\cat T}(\unit)) 
        \]
    defined by $\rho(\cat P) \coloneqq \SET{ f \in \End_{\cat T}(\unit) }{\cone(f) \not\in \cat P}$ and introduced in \cite[Section~5]{Balmer10b}. For any multiplicative subset $S \subset \End_{\cat T}(\unit)$,  let $S^{-1}\cat T$ denote the finite localization of $\cat T$ associated to the Thomason subset $\bigcup_{s \in S}\supp(\cone(s)) \subset \Spc(\cat T^c)$. This ``algebraic localization'' $\cat T \to S^{-1}\cat T$ induces the usual ring-theoretic localization on endomorphism rings: $\End_{\cat T}(\unit) \to \End_{S^{-1}\cat T}(\unit) \cong S^{-1}\End_{\cat T}(\unit)$ (see \cite[Theorem~3.3.7]{HoveyPalmieriStrickland97}). Moreover, the diagram
        \[
            \begin{tikzcd}
                \Spc(\cat T^c) \ar[r,"\rho"] & \Spec(\End_{\cat T}(\unit)) \\
                \Spc((S^{-1}\cat T)^c) \ar[u,hook,] \ar[r,"\rho"] & \Spec(S^{-1}\End_{\cat T}(\unit)) \ar[u,hook]
            \end{tikzcd}
        \]
    induced by $\cat T \to S^{-1} \cat T$ is a pullback (see \cite[Corollary~5.6(c)]{Balmer10b}).
\end{Exa}

\begin{Exa}\label{exa:localization-of-DR}
    Let $S \subset R$ be a multiplicative subset in a commutative ring~$R$. The extension-of-scalars $\Der(R) \to \Der(S^{-1}R)$ coincides (up to tensor-triangular equivalence) with the associated algebraic localization $\Der(R) \to S^{-1}\Der(R)$. Indeed, the monad associated with the $\Der(R) \adjto \Der(S^{-1}R)$ adjunction is (by formal nonsense, e.g., \cite[Lemma~2.8]{BalmerDellAmbrogioSanders15}) the monad associated to the idempotent ring ${S^{-1}R \in \Der(R)}$. It follows (by a direct argument or by invoking \cite[Theorem~1.6]{DellAmbrogioSanders18}) that the adjunction is monadic. In other words, extension-of-scalars $\Der(R) \to \Der(S^{-1}R)$ is the smashing localization associated to the idempotent ring~$S^{-1}R$. Then consider the finite localization $L\colon \Der(R) \to S^{-1}\Der(R)$ of \cref{exa:algebraic-localizations}. The kernel $\Ker(L) = \Loco{\cone(s) \mid s \in S}$ is  contained in $\Ker(-\otimes S^{-1}R)$. Hence $S^{-1}R \cong L(S^{-1}R)$. On the other hand, we claim that the map $R \to S^{-1}R$ becomes an isomorphism after applying $L$. Indeed, for any $X \in \Der(R)$, $H_*(L(X)) \cong S^{-1}H_*(X)$. In particular, $S^{-1}H_0(R) \cong H_*(L(R)) \to H_*(L(S^{-1}R)) \cong S^{-1}H_0(S^{-1}R)$ is indeed an isomorphism. It follows that $S^{-1}R \cong L(R)$.  In other words, the smashing localization associated to $S^{-1}R$ coincides with the algebraic localization associated to $S \subset R = \End_{\Der(R)}(\unit)$.
\end{Exa}

\begin{Exa}\label{exa:local-for-DR}
    Consider the derived category $\Der(R)$ of a commutative ring $R$. It follows from the Neeman--Thomason Classification Theorem \cite{Neeman92a,Thomason97}, that $\Spc(\Der(R)^c) \cong \Spec(R)$. More precisely, the comparison map  $\rho \colon \Spc(\Der(R)^c) \to \Spec(R)$ is a homeomorphism. (This is explained in \cite[Proposition~8.1]{Balmer10b}.) If $\cat P \in \Spc(\Der(R)^c)$ is the (unique) prime such that $\rho(\cat P)= \frak p \in \Spec(R)$, one readily checks that $Y_{\cat P} = \bigcup_{x \in \cat P}\supp(x)$ coincides with $\bigcup_{s \not\in \frak p} \supp(\cone(s))$. In other words, algebraic localization at  $\frak p \in \Spec(R)$ in the sense of \cref{exa:algebraic-localizations} coincides with localization at $\cat P$ in the sense of \cref{exa:local-categories}. In particular,  
        \[
            \Der(R)/\langle \cat P \rangle \cong \Der(R)_{\frak p} \cong \Der(R_{\frak p})
        \]
    where the last equivalence is explained in \cref{exa:localization-of-DR}. We remark in passing that for a commutative ring $R$, the derived category $\Der(R)$ is local in the sense of \cref{Ter:local} if and only if $R$ is a local ring. (This is \cite[Example~4.4]{Balmer10b}.) Indeed, $\Spc(\Der(R)^c)\cong \Spec(R)$ and $\Spec(R)$ has a unique closed point if and only if $R$ is local.
\end{Exa}

\section{Balmer--Favi support}\label{sec:BF-support}

We now introduce the Balmer--Favi notion of support from \cite{BalmerFavi11}, weakening their noetherian assumption on the spectrum to a condition we call weakly noetherian. This extension also includes the case of profinite Balmer spectra and forms the natural context for the development of Balmer--Favi support.

\begin{Def}
    A point $\cat P$ in a spectral space $X$ is said to be \emph{weakly visible} if $\{\cat P\}$ is the intersection of a Thomason subset and the complement of a Thomason subset.
\end{Def}

\begin{Rem}\label{rem:visible}
    This notion of ``weakly visible'' point coincides with the notion of ``visible'' point in \cite{Stevenson14,Stevenson17}. In contrast, we follow the terminology of \cite{BalmerFavi11} and say that a point $\cat P$ is \emph{visible} if its closure $\overbar{\{\cat P\}} \subseteq X$ is a Thomason subset. Every visible point is weakly visible since $\{\cat P\} = \overbar{\{\cat P\}} \cap Y_{\cat P}^c$ is then an intersection of a Thomason and the complement of a Thomason (see \cref{rem:set-of-generalizations} and \cite[Lemma~7.8]{BalmerFavi11}). Since a spectral space is noetherian if and only if every point is visible  \cite[Proposition~7.13]{BalmerFavi11}, we are led to the following terminology:
\end{Rem}

\begin{Def}
    A spectral space is \emph{weakly noetherian} if every point is weakly visible.
\end{Def}

\begin{Rem}
    A spectral space is $T_1$ if and only if it is Hausdorff if and only if it is profinite, see for example \cite[\href{https://stacks.math.columbia.edu/tag/0905}{Lemma 0905}]{stacks-project} and \cite[Section~1.3]{DickmannSchwartzTressl19}. In particular, we note that a $T_1$ spectral space (a.k.a.~profinite space) is noetherian if and only if it is finite. On the other hand, every $T_1$ spectral space is weakly noetherian. Indeed, every singleton $\{ \cat P \} = X \cap \gen(\cat P)$ is the intersection of a Thomason (the whole space itself) and the complement of a Thomason (recall \cref{rem:set-of-generalizations}).
\end{Rem}

\begin{Exa}
    For a concrete incarnation in tt-geometry, the spectrum $\Spec(R) \cong \Spc(\Der(R)^c)$ of a non-noetherian absolutely flat ring $R$ (such as an infinite product of fields) is an example of a spectral space that is weakly noetherian but not noetherian (see \cite[Section~3]{Stevenson14}).
\end{Exa}

\begin{Rem}\label{rem:sub-weakly-noetherian}
    Let $\varphi\colon X \to Y$ be an injective spectral map between spectral spaces. If $Y$ is weakly noetherian, then so is $X$. Indeed, preimages under $\varphi$ of Thomason subsets are Thomason. Thus, if $\{\varphi(\cat P)\} = Y_1 \cap Y_2^c$ then $\{\cat P\} = \varphi^{-1}(\{\varphi(\cat P)\}) = \varphi^{-1}(Y_1) \cap \varphi^{-1}(Y_2)^c$.
\end{Rem}

\begin{Def}[Balmer--Favi]
    Let $\cat T$ be a rigidly-compactly generated tt-category. For each weakly visible point $\cat P \in \Spc(\cat T^c)$, we define a $\otimes$-idempotent $\kappaP \in \cat T$ by
        \[
            \kappaP \coloneqq e_{Y_1} \otimes f_{Y_2}
        \]
    for any choice of Thomason subsets $Y_1, Y_2 \subseteq \Spc(\cat T^c)$ such that $\{\cat P \} = Y_1 \cap Y_2^c$. This object $\kappaP \in \cat T$ does not depend on the choice of Thomason subsets $Y_1, Y_2$ by \cite[Corollary~7.5]{BalmerFavi11}.
\end{Def}

\begin{Rem}\label{rem:choice-of-thomason}
    We can always take the second Thomason $Y_2$ to be $Y_{\cat P}$ (\cref{rem:set-of-generalizations}). Indeed, $\cat P \not\in Y_2$ implies that $Y_2$ contains no generalization of $\cat P$, hence $Y_2 \subseteq Y_{\cat P}$ and it follows that $\{\cat P\} = Y_1 \cap Y_{\cat P}^c$. In particular, we always have $\kappaP \otimes f_{Y_{\cat P}} \simeq \kappaP$. On the other hand, the first Thomason $Y_1$ always contains $\overbar{\{\cat P\}}$. We see that a point $\cat P \in \Spc(\cat T^c)$ is weakly visible if and only if there exists a Thomason closed subset~$Z$ which contains~$\cat P$ but does not contain any other generalization of $\cat P$. When the point is visible, the Thomason closed subset $\overbar{\{\cat P\}}$ works, but in general it could be larger.
\end{Rem}

\begin{Rem}\label{rem:closed-weakly-visible}
	An immediate consequence of \cref{rem:choice-of-thomason} is that the closed point of a local category is weakly visible if and only if it is visible. Moreover, combined with \cref{rem:thomason-closed} and \cref{rem:sub-weakly-noetherian}, ones sees that in general a point $\cat P \in \Spc(\cat T^c)$ is weakly visible if and only if the unique closed point $\cat P \in \Spc(\cat T_{\cat P}^c)$ of the local category at $\cat P$ (\cref{exa:local-categories}) is (weakly) visible.
	In particular, for $X$ a spectral space, a point $x \in X$ is weakly visible if and only if $x \in \gen(x)$ is weakly visible; this also follows, for example, from \cite[Corollary 4.5.13]{DickmannSchwartzTressl19} applied to $X^*$.
\end{Rem}

\begin{Rem}\label{rem:kappa-non-zero}
    It follows from \cref{lem:Y1-and-Y2}(a) that for each weakly visible point $\cat P$, the object $\kappaP$ is nonzero.
\end{Rem}

\begin{Def}[Balmer--Favi]\label{def:balmer-favi}
    Let $\cat T$ be a rigidly-compactly generated tt-category with $\Spc(\cat T^c)$ weakly noetherian. Define the \emph{support} of an object $t \in \cat T$ by
        \[ 
            \Supp(t) \coloneqq \SET{\cat P \in \Spc(\cat T^c)}{t \otimes \kappaP \neq 0}.
        \]
    The basic properties of this notion of support are established in \cite[Section~7]{BalmerFavi11} under the assumption that $\Spc(\cat T^c)$ is noetherian. As we shall see below, it is not difficult to generalize these results to the situation where $\Spc(\cat T^c)$ is weakly noetherian.
\end{Def}

\begin{Rem}\label{rem:support-axioms}
    The following basic properties follow directly from the definition and \cref{rem:kappa-non-zero}:
    \begin{enumerate}
        \item $\Supp(0) = \emptyset$ and $\Supp(\unit) = \Spc(\cat T^c)$;
        \item $\Supp(\Sigma t) = \Supp(t)$ for every $t \in\cat T$;
        \item $\Supp(c) \subseteq \Supp(a) \cup \Supp(b)$ for any exact triangle $a \to b \to c \to \Sigma a$ in~$\cat T$;
        \item $\Supp(\coprod_{i\in I} t_i) = \bigcup_{i \in I} \Supp(t_i)$ for any set of objects $t_i$ in $\cat T$;
        \item $\Supp(t_1 \otimes t_2) \subseteq \Supp(t_1) \cap \Supp(t_2)$ for any $t_1,t_2 \in \cat T$.
    \end{enumerate}
    These properties (excluding $\Supp(\unit) = \Spc(\cat T^c)$) are equivalent to the statement that $\cat T_Y \coloneqq \SET{ t \in \cat T}{\Supp(t) \subseteq Y}$ is a localizing $\otimes$-ideal of $\cat T$ for any subset $Y \subseteq \Spc(\cat T^c)$.
\end{Rem}

\begin{Lem}\label{lem:support-compat-with-finite-loc}
    For any object $t \in \cat T$ and Thomason subset $Y \subseteq \Spc(\cat T^c)$, we have
        \[
            \Supp(t \otimes e_Y) = \Supp(t) \cap Y \quad\text{ and }\quad \Supp(t \otimes f_Y) = \Supp(t) \cap Y^c.
        \]
    In particular, $\Supp(\kappaP) = \{\cat P\}$ for any $\cat P \in \Spc(\cat T^c)$.
\end{Lem}

\begin{proof}
    First we establish that $\Supp(e_Y) \subseteq Y$ and $\Supp(f_Y) \subseteq Y^c$. To this end, let $\cat P \in \Spc(\cat T^c)$ and write $\{\cat P\} = Y_1 \cap Y_2^c$. If $\cat P \not\in Y$ then $Y \cap Y_1 \cap Y_2^c = \emptyset$ so that $Y \cap Y_1 = Y \cap Y_1 \cap Y_2 \subseteq Y_2$. It then follows from \cref{lem:Y1-and-Y2} that $e_Y \otimes e_{Y_1} \otimes f_{Y_2} = 0$ so that $\cat P \not\in \Supp(e_Y)$. On the other hand, if $\cat P \in Y$ then $Y_1 \cap Y_2^c \subseteq Y$ so that $Y_1 \subseteq Y_2 \cup Y$. Hence $e_{Y_1} \otimes f_{Y_2} \otimes f_Y = 0$, again by \cref{lem:Y1-and-Y2}, so that $\cat P \not\in \Supp(f_Y)$. Now consider an arbitrary $t \in \cat T$. By property (e) of \cref{rem:support-axioms} we have $\Supp(t \otimes e_Y) \subseteq \Supp(t) \cap Y$ and $\Supp(t \otimes f_Y) \subseteq \Supp(t) \cap Y^c$. Moreover, from the exact triangle $t \otimes e_Y \to t \to t\otimes f_Y \to \Sigma t\otimes e_Y$ we have $\Supp(t) \subseteq \Supp(t \otimes e_Y) \cup \Supp(t \otimes f_Y)$. Intersecting with $Y$ and $Y^c$ we obtain $\Supp(t) \cap Y \subseteq \Supp(t \otimes e_Y)$ and $\Supp(t) \cap Y^c \subseteq \Supp(t \otimes f_Y)$, which completes the proof.
\end{proof}

\begin{Prop}\label{prop:uniqueness-of-support}
    Let $\cat T$ be a rigidly-compactly generated tt-category with $\Spc(\cat T^c)$ weakly noetherian. The Balmer--Favi notion of support is the only assignment of a subset $\sigma(t) \subseteq \cat \Spc(\cat T^c)$ to each object $t \in \cat T$ which can satisfy the following two properties:
    \begin{enumerate}
        \item \label{it:unique-support-a} For every $t \in \cat T$, $\sigma(t) = \emptyset$ implies $t=0$.
        \item \label{it:unique-support-b} For every $t \in \cat T$ and Thomason subset $Y \subseteq \Spc(\cat T^c)$, $\sigma(t \otimes e_Y) = \sigma(t) \cap Y$ and $\sigma(t \otimes f_Y) = \sigma(t) \cap Y^c$.
    \end{enumerate}
\end{Prop}

\begin{proof}
    Indeed, if $\sigma$ satisfies (b) then $\sigma(\kappaP) = \{\cat P\}$ and $\sigma(t\otimes \kappaP) = \sigma(t) \cap \{\cat P\}$. We then have $\cat P \in \sigma(t)$ iff $\cat P \in \sigma(t \otimes \kappaP)$ iff $\sigma(t \otimes \kappaP) \neq \emptyset$ iff $t \otimes \kappaP \neq 0$ iff $\cat P \in \Supp(t)$.
\end{proof}

\begin{Def}[The detection property]\label{def:detectionproperty}
    A rigidly-compactly generated tensor-triangulated category $\cat T$ with $\Spc(\cat T^c)$ weakly noetherian is said to have the \emph{detection property} if $\Supp(t) = \emptyset$ implies $t=0$ for every $t \in \cat T$.
\end{Def}

\begin{Rem}\label{rem:bfsupport-idempotents-compatibility}
    The Balmer--Favi notion of support (\cref{def:balmer-favi}) always satisfies property~\hyperref[it:unique-support-b]{(b)} of \cref{prop:uniqueness-of-support} by \cref{lem:support-compat-with-finite-loc}. However, we do not know in what generality the detection property holds.
\end{Rem}

\begin{Rem}
    We finish this section with the observation, generalizing \cite[Theorem 7.22]{BalmerFavi07}, that the Balmer--Favi support satisfies a half-smash product formula and thus in particular coincides with Balmer's universal support function on compact objects:
\end{Rem}

\begin{Lem}\label{lem:half-tensor}
    For any compact $x\in \cat T^c$ and arbitrary $t \in \cat T$,
        \[ 
            \Supp(x \otimes t) = \supp(x) \cap \Supp(t).
        \]
    In particular, for any compact object $x \in \cat T^c$, the Balmer--Favi notion of support coincides with the usual notion of support: $\Supp(x) = \supp(x)$.
\end{Lem}

\begin{proof}
    Let $Y\coloneqq \supp(x)$. Since $x$ is compact, $\Loco{x} = \Loco{\cat T_Y^c} = \Loco{e_Y} = e_Y \otimes \cat T$ in the notation of \cref{rem:finite-localizations}. In particular, $x \otimes t \simeq e_Y \otimes x \otimes t$, hence $\Supp(x\otimes t) = \Supp(e_Y \otimes x \otimes t) = \Supp(x\otimes t)\cap Y$ by \cref{lem:support-compat-with-finite-loc}. Thus ${\Supp(x\otimes t) \subseteq Y}$. Applying property (e) of \cref{rem:support-axioms}, it follows that $\Supp(x \otimes t) \subseteq \Supp(t) \cap Y$. On the other hand, if $\cat P \not\in \Supp(x \otimes t)$ so that $x\otimes t \otimes \kappaP = 0$ then $\Loco{x} \otimes t \otimes \kappaP = 0$. This implies that $e_Y \otimes t \otimes \kappaP = 0$ so that $\cat P \not\in \Supp(e_Y \otimes t) = \Supp(t) \cap Y$. Therefore $\Supp(x \otimes t) = \Supp(t) \cap Y$. Finally, the last claim follows by specializing to $t = \unit$. 
\end{proof}

\section{The local-to-global principle}\label{sec:LTG}

The Balmer--Favi notion of support extends in a natural way from objects to localizing $\otimes$-ideals. As such, it can be used to study the classification of localizing $\otimes$-ideals. This problem splits into two parts:
    \begin{enumerate}
        \item The local-to-global principle, which formulates how all localizing $\otimes$-ideals may be obtained from $\kappaP \otimes \cat T = \Loco{\kappaP}$ for $\cat P \in \Spc(\cat T^c)$. This is the topic of the present section.
        \item The minimality of the candidate irreducible localizing $\otimes$-ideals $\kappaP \otimes \cat T$. This will be studied in \cref{part:stratification1}. 
    \end{enumerate}

\begin{Def}\label{def:supp-of-localizing}
    For a localizing $\otimes$-ideal $\cat L$ of $\cat T$, define $\Supp(\cat L) \coloneqq \bigcup_{t \in \cat L} \Supp(t)$. 
\end{Def}

\begin{Rem}\label{rem:supp-of-localizing}
    It is a routine exercise using the basic properties of support (\cref{rem:support-axioms}) that if $\cat L=\Loco{\cat E}$ for some collection of objects~$\cat E$ then $\Supp(\cat L) = \bigcup_{t \in \cat E} \Supp(t)$. In particular, if $\Loco{t_1} = \Loco{t_2}$ then $\Supp(t_1) = \Supp(t_2)$.
\end{Rem}

\begin{Rem}\label{rem:set-generated}
    If a localizing $\otimes$-ideal is generated by a set of objects then it is also generated by a single object: $\Loco{\cat E} = \Loco{\coprod_{t \in \cat E}t}$. We'll refer to such localizing \mbox{$\otimes$-ideals} as the ``set-generated'' localizing $\otimes$-ideals. As we shall see below, when stratification holds, every localizing $\otimes$-ideal will necessarily be set-generated.
\end{Rem}

\begin{Lem}\label{lem:surjectivity}
    Let $\cat T$ be a rigidly-compactly generated tensor-triangulated category with $\Spc(\cat T^c)$ weakly noetherian. Then the map
        \[
            \Supp \colon \big\{ \text{set-generated localizing $\otimes$-ideals of $\cat T$}\big\} \to \big\{\text{subsets of $\Spc(\cat T^c)$}\big\}
        \]
    is surjective.
\end{Lem}

\begin{proof}
    For any subset $Y \subseteq \Spc(\cat T^c)$, consider the set-generated localizing $\otimes$-ideal
        \[ 
            \cat L_Y \coloneqq \Loco{\kappaP \mid \cat P \in Y}.
        \]
    Recall from \cref{lem:support-compat-with-finite-loc} that $\Supp(\kappaP) = \big\{\cat P\big\}$. Hence
        \[
            \Supp(\cat L_Y) = \bigcup_{\cat P \in Y} \Supp(\kappaP) = Y
        \]
    by \cref{rem:supp-of-localizing}.
\end{proof}

\begin{Prop}[Krause--Stevenson]\label{prop:only-set-then-all-set}
    Let $\cat T$ be a rigidly-compactly generated tt-category. If the class of set-generated localizing $\otimes$-ideals of $\cat T$ forms a set, then every localizing $\otimes$-ideal of $\cat T$ is generated by a set.
\end{Prop}

\begin{proof}
    The proof of \cite[Lemma~3.3.1]{KrauseStevenson19} goes through by replacing all instances of ``localizing subcategory'' and ``$\Loc$'' with ``localizing $\otimes$-ideal'' and ``$\Loc_{\otimes}$''. If there exists a localizing $\otimes$-ideal of $\cat T$ which is not generated by a set, their argument constructs by transfinite recursion an ordinal-indexed chain of distinct \emph{set-generated} localizing $\otimes$-ideals. This contradicts the hypothesis since the collection of ordinals does not form a set and cannot be embedded in the set of set-generated localizing $\otimes$-ideals. In particular, the hypothesis of \cite[Lemma~3.3.1]{KrauseStevenson19} can be weakened to the assumption that there is only a set of \emph{set-generated} localizing subcategories;
	this was already observed, at least implicitly, in the proof of
	\cite[Lemma~3.3.4]{KrauseStevenson19}.
\end{proof}

\begin{Lem}\label{lem:tensor-in}
    Let $\cat G$ be a class of objects in $\cat T$. If $e \in \Loco{\cat G}$ then $e\otimes f \in \Loco{\cat G\otimes f}$ for any two objects $e,f\in \cat T$. In particular, if $e$ is a $\otimes$-idempotent then
        \[ 
            e\in \Loco{\cat G} \Longrightarrow e \in \Loco{\cat G \otimes e}.
        \]
\end{Lem}

\begin{proof}
    Note that $ - \otimes f$ is a coproduct-preserving triangulated functor $F:\cat T \to \cat T$ with the property that $F(s\otimes t) \simeq s \otimes F(t)$ for all $s,t \in \cat T$. It is a standard lemma that if $F$ is such a functor and $e \in \Loco{\cat G}$, then $F(e) \in \Loco{F(\cat G)}$. Indeed, consider the full subcategory of $\cat T$ consisting of those $t \in \cat T$ such that $F(t) \in \Loco{F(\cat G)}$. This is a localizing $\otimes$-ideal which contains $\cat G$ and hence contains the object $e$.
\end{proof}

\begin{Lem}\label{lem:big-supp-in-small}
    Let $\cat T$ be a rigidly-compactly generated tensor-triangulated category with $\Spc(\cat T^c)$ weakly noetherian and having the detection property (\cref{def:detectionproperty}). Let $t \in \cat T$ and $x \in \cat T^c$ be objects such that
        \[
            \Supp(t) \subseteq \supp(x).
        \]
    Then $t \in \Loco{t\otimes x}$.
\end{Lem}

\begin{proof}
    Consider the finite localization (\cref{rem:finite-localizations}) associated to the Thomason subset $\supp(x)$:
        \[
            e \to \unit \to f\to \Sigma e.
        \]
    Observe that $\Loco{e} = \Loc\langle \cat T^c_{\supp(x)}\rangle = \Loc\langle \thickt{x}\rangle = \Loco{x}$. Then consider  
        \[ 
            e\otimes t \to t \to f\otimes t \to \Sigma e\otimes t.
        \]
    Note that $\Supp(f \otimes t) = \Supp(t) \cap \supp(x)^c =\emptyset$ by \cref{lem:support-compat-with-finite-loc}. Hence $f\otimes t=0$ by the detection property. Therefore $t \simeq e\otimes t$. Since $e \in \Loco{x}$, we conclude that  $t\simeq t\otimes e \in \Loco{t\otimes x}$ by \cref{lem:tensor-in}.
\end{proof}

\begin{Def}[The local-to-global principle]
    Let $\cat T$ be a rigidly-compactly generated tensor-triangulated category with $\Spc(\cat T^c)$ weakly noetherian. We say that $\cat T$ satisfies the \emph{local-to-global principle} if
        \[ 
            \Loco{t} = \Loco{t \otimes \kappaP \mid \cat P \in \Spc(\cat T^c)}
        \]
    for every object $t \in \cat T$. Of course, the right-hand side is the same as $\Loco{{t \otimes \kappaP} \mid \cat P \in \Supp(t)}$. 
\end{Def}

\begin{Rem}\label{rem:ltgimpliesdetection}
    Note that the local-to-global principle for $\cat T$ implies that the detection property holds for $\cat T$, i.e., that $t=0$ if and only if ${\Supp(t) = \emptyset}$ for any $t \in \cat T$.
\end{Rem}

\begin{Cor}
	If $\cat T$ satisfies the local-to-global principle
	then the Balmer--Favi notion of support is the unique assignment of a subset $\Supp(t) \subseteq \Spc(\cat T^c)$ to each object $t \in \cat T$ satisfying the following properties: 
		\begin{enumerate}
			\item For every $t \in \cat T$, $\Supp(t) = \emptyset$ implies $t=0$.
			\item  For every $t \in \cat T$ and Thomason subset $Y \subseteq \Spc(\cat T^c)$,
			$\Supp(t \otimes e_Y) = \Supp(t) \cap Y$ and $\Supp(t \otimes f_Y) = \Supp(t) \cap Y^c$.
		\end{enumerate}
\end{Cor}

\begin{proof}
    This is a consequence of \cref{prop:uniqueness-of-support}. Indeed, part \hyperref[it:unique-support-a]{(a)} holds by 
	\cref{rem:ltgimpliesdetection}, while the Balmer--Favi notion of support satisfies \hyperref[it:unique-support-a]{(b)} by \cref{rem:bfsupport-idempotents-compatibility}.
\end{proof}

\begin{Rem}
    Our next task is to demonstrate that the local-to-global principle passes to finite localizations. We'll state the result in more generality in case it could be useful for future applications.
\end{Rem}

\begin{Prop}\label{prop:ltg-permanence}
    Let $F\colon\cat C \to \cat D$ be a coproduct-preserving tensor-triangulated functor between rigidly-compactly generated tensor-triangulated categories and let~$U$ denote its right adjoint. Suppose that the induced map
        \[
            \varphi\colon\Spc(\cat D^c) \to \Spc(\cat C^c)
        \]
    is injective and that $\Spc(\cat C^c)$ is weakly noetherian. Then $\Spc(\cat D^c)$ is also weakly noetherian and the following hold:
        \begin{enumerate}
            \item For each $\cat P \in \Spc(\cat C^c)$, we have
                \[
                    F(\kappaCP) =
                        \begin{cases}
                            0        & \text{ if $\cat P \not\in \im(\varphi)$}\\
                            \kappaDQ & \text{ if $\cat P=\varphi(\cat Q).$}
                        \end{cases}
                \]
            \item For any $c \in \cat C$, we have $\Supp(F(c)) \subseteq \varphi^{-1}(\Supp(c))$.
            \item For any $d \in \cat D$, $\Supp(U(d)) \subseteq \varphi(\Supp(d))$ with equality when $U$ is conservative.
            \item If $\cat C$ satisfies the detection property and the right adjoint $U$ is conservative then $\cat D$ also satisfies the detection property.
            \item If $\cat C$ satisfies the local-to-global principle then $\cat D$ also satisfies the local-to-global principle.
        \end{enumerate}
\end{Prop}

\begin{proof}
    Note that $\Spc(\cat D^c)$ is also weakly noetherian by \cref{rem:sub-weakly-noetherian}. If $e_Y \to \unit \to f_Y \to \Sigma e_Y$ is the idempotent triangle associated to a Thomason subset $Y \subseteq \Spc(\cat C^c)$ (\cref{rem:finite-localizations}), then $F(e_Y) \simeq e_{\varphi^{-1}(Y)}$ and $F(f_Y) \simeq f_{\varphi^{-1}(Y)}$ (\cref{rem:geometric-functor}). Thus, if $\{ \cat P \} = Y_1 \cap Y_2^c$ then 
        \[
            F(\kappaCP) \simeq F(e_{Y_1}) \otimes F(f_{Y_2}) \simeq e_{\varphi^{-1}(Y_1)} \otimes f_{\varphi^{-1}(Y_2)}.
        \]
    By \cref{lem:Y1-and-Y2}, we then have 
        \begin{align*}
            F(\kappaCP) = 0
            &\Leftrightarrow \varphi^{-1}(Y_1) \subseteq \varphi^{-1}(Y_2)\\
            &\Leftrightarrow \varphi^{-1}(Y_1 \cap Y_2^c) = \emptyset\\
            &\Leftrightarrow \varphi^{-1}(\{\cat P\}) = \emptyset\\
            &\Leftrightarrow \cat P \not\in \im(\varphi).
        \end{align*}
    On the other hand, for a point $\cat P=\varphi(\cat Q)$ in the image, the hypothesis that $\varphi$ is injective ensures that $\{\cat Q\} = \varphi^{-1}(\{\varphi(\cat Q)\}) = \varphi^{-1}(Y_1) \cap \varphi^{-1}(Y_2)^c$. Hence $\kappaDQ \simeq e_{\varphi^{-1}(Y_1)} \otimes f_{\varphi^{-1}(Y_2)} \simeq F(\kappaCP)$. This proves part (a).

    For part (b), suppose that $\cat Q \in \Supp(F(c))$, so that $\kappaDQ\otimes F(c) \neq 0$. By part~(a), this is $F(\kappaCphiQ\otimes c)\neq 0$. Hence $\kappaCphiQ \otimes c \neq 0$, so that $\varphi(\cat Q) \in \Supp(c)$.

    For part (c), observe that $\cat P \in \Supp(U(d))$ if and only if 
        \begin{equation}\label{eq:partb}
            \kappaCP \otimes U(d) \simeq U(F(\kappaCP) \otimes d)
        \end{equation}
    is nonzero. Thus, if $\cat P \in \Supp(U(d))$ then $F(\kappaCP) \otimes d \neq 0$. Hence by part~(a), $\cat P = \varphi(\cat Q)$ for some (unique) $\cat Q \in \Spc(\cat D^c)$ and $\kappaDQ \otimes d \simeq F(\kappaCP)\otimes d \neq 0$. This establishes $\Supp(U(d)) \subseteq \varphi(\Supp(d))$. Conversely, if $\cat Q \in \Supp(d)$ then $0 \neq \kappaDQ \otimes d \simeq F(\kappaCphiQ)\otimes d$ which, if $U$ is conservative, implies that \eqref{eq:partb} is nonzero (with $\cat P \coloneqq \varphi(\cat Q)$) so that $\varphi(\cat Q) \in \Supp(U(d))$.

    For part (d), consider an object $d \in \cat D$ with $\Supp(d) = \emptyset$. By part (c), we see that $\Supp(U(d)) = \emptyset$, so $U(d) = 0$ by the detection property for $\cat C$. If $U$ is conservative, this implies $d=0$, and we conclude that $\cat D$ also has the detection property.

    Part (e) follows from part (a) since if $\unit_{\cat C} \in \Loco{\kappaCP \mid \cat P \in \Spc(\cat C^c)}$ then $\unit_{\cat D} = F(\unit_{\cat C}) \in \Loco{F(\kappaCP) \mid \cat P \in \Spc(\cat C^c)} = \Loco{\kappaDQ \mid \cat Q \in \Spc(\cat D^c)}$ and hence $d \in \Loco{d \otimes \kappaDQ \mid \cat Q \in \Spc(\cat D^c)}$ for any $d \in \cat D$ by \cref{lem:tensor-in}.
\end{proof}

\begin{Cor}\label{cor:ltg-passes-to-finite}
    Let $\cat T$ be a rigidly-compactly generated tt-category with $\Spc(\cat T^c)$ weakly noetherian. If $\cat T$ satisfies the local-to-global principle, then so does each finite localization~$\cat T(V)$.
\end{Cor}

\begin{Rem}\label{rem:LGP-not-local}
    Note that the local-to-global principle cannot be checked one prime at a time. For example, if $R$ is an absolutely flat ring then each local  category $\Der(R)/\langle \cat P \rangle \cong \Der(R_{\frak p})$ (\cref{exa:local-for-DR})	satisfies the local-to-global principle (since each localization $R_{\frak p}$ is a field) and yet $\Der(R)$ need not satisfy the local-to-global principle, as will be discussed in \cref{exa:LGP-for-absolutely-flat}. Nevertheless, it can be checked on an open cover (see \cref{prop:ltg-cover} and \cref{rem:ltg-cover} below). First we need the following lemma:
\end{Rem}

\begin{Lem}\label{lem:pullback-localizing}
    Let $L:\cat T \adjto \cat T(V) : U$ denote finite localization with respect to the Thomason subset $Y \subseteq \Spc(\cat T^c)$. We have an inclusion-preserving bijection
        \[ 
            \big\{ \text{localizing $\otimes$-ideals of $\cat T(V)$} \big\} \xra{L^{-1}} \big\{ \text{localizing $\otimes$-ideals of $\cat T$ which contain $e_Y$} \big\}.
        \]
    whose inverse is given by $U^{-1}$. Explicitly, for a class $\cat E$ of objects of $\cat T(V)$,
        \begin{equation}\label{eq:explicit-pullback}
            L^{-1}(\Loco{\cat E}) = \Loco{e_Y,U(\cat E)}.
        \end{equation}
\end{Lem}

\begin{proof}
    Note that $L(e_Y)=0$ so the map $L^{-1}$ in the statement is well-defined. It is evidently injective since $L$ is essentially surjective. Then consider a localizing \mbox{$\otimes$-ideal} $\cat L$ of $\cat T$. Since the right adjoint $U$ preserves coproducts, $U^{-1}(\cat L)$ is a localizing subcategory of $\cat T(V)$. Moreover, since $L$ is essentially surjective, the projection formula  $U(L(s)\otimes t) \simeq s\otimes U(t)$ (\cref{rem:geometric-functor}) implies that $U^{-1}(\cat L)$ is a localizing $\otimes$-ideal of $\cat T(V)$. Then consider $L^{-1}U^{-1}(\cat L)$. We always have $\cat L \subseteq L^{-1}U^{-1}(\cat L)$ since if $t \in \cat L$ then $U(L(t)) \simeq U(\unit)\otimes t \simeq f_Y \otimes t \in \cat L$. On the other hand, if $e_Y \in \cat L$, then the exact triangle
        \[
            e_Y \otimes t \to t \to f_Y \otimes t \to \Sigma e_Y \otimes t
        \]
    demonstrates the reverse inclusion. So the map in the statement is a bijection with inverse $U^{-1}$. Now for the explicit description \eqref{eq:explicit-pullback}. First note that $\cat E \subseteq U^{-1}(\Loco{e_Y,U(\cat E)})$ and therefore $\Loco{\cat E} \subseteq U^{-1}(\Loco{e_Y,U(\cat E)})$ since the right-hand side is a localizing $\otimes$-ideal. Hence 
        \[
            L^{-1}(\Loco{\cat E}) \subseteq L^{-1}U^{-1}(\Loco{e_Y,U(\cat E)}) = \Loco{e_Y,U(\cat E)}.
        \]
    On the other hand, $LU \simeq \Id$ since the right adjoint $U$ is fully faithful. Hence $U(\cat E) \subseteq L^{-1}(\Loco{\cat E})$ and it follows that $\Loco{e_Y,U(\cat E)} \subseteq L^{-1}(\Loco{\cat E})$.
\end{proof}

\begin{Prop}\label{prop:ltg-cover}
    Let $\cat T$ be a rigidly-compactly generated tt-category with $\Spc(\cat T^c)$ weakly noetherian. Suppose $\Spc(\cat T^c) = V_1 \cup \cdots \cup V_n$ is a finite cover by complements of Thomason sets. If $\cat T(V_i)$ satisfies the local-to-global principle for each $i$ then $\cat T$ satisfies the local-to-global principle.
\end{Prop}

\begin{proof}
    Let $Y_i \coloneqq \Spc(\cat T^c)\setminus V_i$ denote the Thomason complement. If $L_i \colon \cat T \to \cat T(V_i)$ denotes the corresponding localization which induces $\Spc(\cat T(V_i)^c) \cong V_i \subseteq \Spc(\cat T^c)$ on spectra (\cref{rem:finite-localizations}), then \[L_i(\unit) =\unit \in \Loco{\kappaTViQ \mid \cat Q \in V_i} = \Loco{L_i(\kappaQ) \mid \cat Q \in V_i}\] by the assumed local-to-global principle in $\cat T(V_i)$ and \cref{prop:ltg-permanence}. Hence by \cref{lem:pullback-localizing}, $\unit \in L_i^{-1}\Loco{L_i(\kappaQ) \mid \cat Q \in V_i} = \Loco{e_{Y_i}, \{f_{Y_i} \otimes \kappaQ \}_{\cat Q \in V_i}}$. \Cref{lem:tensor-in} then implies that \[f_{Y_i} \in \Loco{\kappaQ \otimes f_{Y_i} \mid \cat Q \in V_i} \subseteq \Loco{\kappaQ \mid \cat Q \in \Spc(\cat T^c)}.\] Thus, $f_{Y_i} \in \Loco{\kappaQ \mid \cat Q \in \Spc(\cat T^c)}$ for all $i$. Inductively applying Abstract Mayer--Vietoris \cite[Theorem~3.13]{BalmerFavi11}, we obtain \begin{equation}\label{eq:fY1Yn} f_{Y_1 \cap \cdots \cap Y_n} \in \Loco{\kappaQ \mid \cat Q \in \Spc(\cat T^c)}. \end{equation} Since $Y_1 \cap \cdots \cap Y_n = \emptyset$, we have $f_{Y_1 \cap \cdots \cap Y_n} = \unit$ and the local-to-global principle for~$\cat T$ then follows from \eqref{eq:fY1Yn} by applying \cref{lem:tensor-in}.
\end{proof}

\begin{Rem}\label{rem:ltg-cover}
    As explained in \cref{rem:LGP-not-local}, the statement of the proposition is no longer true (in general) if we take an arbitrary (not necessarily finite) cover of $\Spc(\cat T^c)$ by complements of Thomason subsets. On the other hand, we can take an arbitrary cover of $\Spc(\cat T^c)$ by quasi-compact opens, since such a cover can be reduced to a finite cover because $\Spc(\cat T^c)$ is quasi-compact.
\end{Rem}

\begin{Rem}
    While we do not have a general characterization for when the local-to-global principle holds, we can show that it is satisfied whenever the topological space $\Spc(\cat T^c)$ is noetherian:
\end{Rem}

\begin{Thm}\label{thm:local-to-global}
    Let $\cat T$ be a rigidly-compactly generated tt-category with $\Spc(\cat T^c)$ noetherian. Then $\cat T$ satisfies the local-to-global principle.
\end{Thm}

\begin{proof}
    The proof of \cite[Theorem 3.6]{BensonIyengarKrause11a} goes through in our setting. Consider $\cat L \coloneqq \Loco{\kappaP \mid \cat P \in \Spc(\cat T^c)}$. It suffices to establish that $\unit \in \cat L$ since \cref{lem:tensor-in} will then imply that $t \in \Loco{t \otimes \kappaP \mid \cat P \in \Spc(\cat T^c)}$ for every $t \in \cat T$. To this end consider 
        \[
            Y \coloneqq \SET{\cat P \in \Spc(\cat T^c)}{ e_{\overbar{\{\cat P\}}} \in \cat L}.
        \]
    This is a specialization closed set since if $\cat Q \subseteq \cat P$ then $e_{\overbar{\{\cat Q\}}} \in \Loco{e_{\overbar{\{\cat P\}}}}$. We claim that $Y = \Spc(\cat T^c)$. If this were not the case then, since $\Spc(\cat T^c)$ is noetherian, we could choose a prime $\cat Q$ which is minimal in $\Spc(\cat T^c)\setminus Y$. Then let $Z := {\overbar{\{\cat Q\}} \setminus \{\cat Q\}}$, which is specialization closed and thus Thomason, as the ambient space is noetherian. By construction, we have $e_Z \simeq e_Z \otimes e_{\overbar{\{\cat Q\}}}$. Moreover, since $\overbar{\{\cat Q\}} \setminus Z = \{\cat Q\}$, \cite[Lemma 7.4]{BalmerFavi11} implies that $f_Z \otimes e_{\overbar{\{\cat Q\}}} \simeq g(\cat Q)$. Thus, tensoring the exact triangle $e_Z \to \unit \to f_Z \to \Sigma e_Z$ with $e_{\overbar{\{\cat Q\}}}$ provides an exact triangle
        \begin{equation}\label{eq:ex-tri}
            e_Z \to e_{\overbar{\{\cat Q\}}} \to g(\cat Q) \to \Sigma e_Z.
        \end{equation}
    Now since $Y$ is specialization closed, we have $Y=\bigcup_{\cat P \in Y} \overbar{\{ \cat P \}}$. Hence
        \[
            \Loco{e_Y} = \Loco{e_{\overbar{\{\cat P\}}} \mid \cat P \in Y} \subseteq \cat L
        \]
    by \cref{rem:union-of-thomason}. Then note that $Z \subseteq Y$ since $\cat Q$ is minimal in the complement of~$Y$. Hence $e_Z \in \Loco{e_Y}$ so that $e_Z \in \cat L$. On the other hand, $\cat L$ contains $g(\cat Q)$ by definition. Hence by the exact triangle \eqref{eq:ex-tri}, $\cat L$ also contains $e_{\overbar{\{\cat Q\}}}$. This is a contradiction since $\cat Q$ is not contained in $Y$ by assumption. We conclude that $Y=\Spc(\cat T^c)$. It follows that $\unit = e_Y \in \Loco{ e_{\overbar{\{\cat P\}}} \mid \cat P \in Y} \subseteq \cat L$ again by \cref{rem:union-of-thomason}.
\end{proof}

\begin{Rem}
    \Cref{thm:local-to-global} demonstrates that the extra hypotheses of \cite[Theorem 6.9]{Stevenson13} and \cite[Theorem~4.7]{Stevenson17} are unnecessary. For a rigidly-compactly generated tt-category $\cat T$ with $\Spc(\cat T^c)$ noetherian, the local-to-global principle \emph{always holds} without needing any additional hypotheses.
\end{Rem}

\begin{Exa}\label{exa:LGP-for-absolutely-flat}
    Stevenson \cite{Stevenson14,Stevenson17} proves that the derived category $\Der(R)$ of an absolutely flat ring satisfies the local-to-global principle if and only if $R$ is semi-artinian. Note that the spectrum $\Spec(R)$ of an absolutely flat ring is always weakly noetherian. On the other hand, for an absolutely flat ring, $\Spec(R)$ is noetherian  $\Leftrightarrow$ the ring $R$ is noetherian $\Leftrightarrow$ $R$ is artinian  $\Leftrightarrow$ $R$ is semisimple $\Leftrightarrow$ $R$ is a finite product of fields $\Leftrightarrow$ $\Spec(R)$ is finite. Moreover, a ring is semi-artinian and noetherian if and only if it is artinian (see \cite[Corollaire~3.1]{NastasescuPopescu68}). In particular, any non-noetherian non-semi-artinian absolutely flat ring (such as an infinite product of fields) provides an example where the spectrum is weakly noetherian and yet the local-to-global principle fails. Thus, the hypothesis in  \cref{thm:local-to-global} that $\Spc(\cat T^c)$ is noetherian cannot be weakened to the hypothesis that $\Spc(\cat T^c)$ is weakly noetherian. On the other hand, any non-noetherian semi-artinian absolutely flat ring  provides an example where the local-to-global principle holds and yet the spectrum is only weakly noetherian. As explained in \cite[Section 4]{NastasescuPopescu68}, an example of such a ring is the subring of $\prod_{\bbN}\Fp$ consisting of those sequences which are eventually constant. 
\end{Exa}

\begin{Rem}
    The local-to-global principle is also studied in slightly different settings than our present one. See \cite[Theorem~3.4]{BensonIyengarKrause11b}, \cite[Proposition~6.8]{Stevenson13}, \cite[Section 4]{Stevenson14}, and \cite{Stevenson17} for positive and negative results in this direction.
\end{Rem}

\part{Stratification: classification and criteria}\label{part:stratification1}

We now begin our study of stratification. Informally, a rigidly-compactly generated tt-category $\cat T$ is stratified if its Balmer spectrum $\Spc(\cat T^c)$ is weakly noetherian and the Balmer--Favi notion of support classifies the localizing $\otimes$-ideals of $\cat T$. It is a property that a given $\cat T$ might or might not possess. Our present goal is to formulate this notion precisely and establish conditions and criteria for when stratification holds. For example, we'll reduce the problem to the local categories at each point and also establish a form of finite \'{e}tale descent for stratification. These results will provide some basic tools for establishing that a given category is stratified.

\section{Stratification via Balmer--Favi supports}\label{sec:strat}

In this section we introduce the notion of stratification for a rigidly-compactly generated tensor-triangulated category with $\Spc(\cat T^c)$ weakly noetherian. The next result is inspired by and closely related to the work of Benson, Iyengar, and Krause, see in particular \cite[Theorem~4.2]{BensonIyengarKrause11b} and \cref{rem:BIK}. 

\begin{Thm}\label{thm:equiv-strat}
    Let $\cat T$ be a rigidly-compactly generated tensor-triangulated category with $\Spc(\cat T^c)$ weakly noetherian. Then the following are equivalent:
\begin{enumerate}
    \item ``Minimality'': The local-to-global principle holds for $\cat T$ and for each $\cat P\in \Spc(\cat T^c)$, $\Loco{\kappaP}$ is a minimal localizing $\otimes$-ideal of $\cat T$.
    \item For all $t \in \cat T$, $\Loco{t} = \Loco{\kappaP\mid \cat P \in \Supp(t)}$.
    \item ``Stratification'': The map
            \[ 
                \big\{ \text{localizing $\otimes$-ideals of $\cat T$} \big\} \xra{\Supp} \big\{ \text{subsets of $\Spc(\cat T^c)$}\big\}
            \]
            is an injection (and hence a bijection).
\end{enumerate}
\end{Thm}

\begin{proof}
    $(a)\Rightarrow(b)$: By the local-to-global principle,
        \[ 
            \Loco{t} = \Loco{\kappaP\otimes t \mid \cat P \in \Supp(t)}
        \]
    and certainly
        \begin{equation}\label{eq:loc-kappa-in}
            \Loco{\kappaP\otimes t\mid \cat P \in \Supp(t)} \subseteq \Loco{\kappaP \mid \cat P \in \Supp(t)}.
        \end{equation}
    On the other hand, if $\cat P \in \Supp(t)$, then $0\neq \Loco{\kappaP\otimes t} \subseteq \Loco{\kappaP}$, hence by the minimality hypothesis, $\Loco{\kappaP\otimes t} = \Loco{\kappaP}$. So $\kappaP \in \Loco{\kappaP \otimes t}$. Thus \eqref{eq:loc-kappa-in} is an equality.

    $(b)\Rightarrow(c)$: We actually show that the map
        \begin{equation}\label{eq:supp-set-gen}
            \big\{\text{set-generated localizing $\otimes$-ideals of $\cat T$}\big\} \xra{\Supp} \big\{\text{subsets of $\Spc(\cat T^c)$}\big\}
        \end{equation}
    is an injection. It will then follow that there is only a set of set-generated localizing $\otimes$-ideals. Hence, by \cref{prop:only-set-then-all-set}, every localizing $\otimes$-ideal is set-generated, so we will have established (c). Note that the map is always surjective by  \cref{lem:surjectivity}. Recall that every set-generated localizing $\otimes$-ideal of $\cat T$ is generated by a single object (\cref{rem:set-generated}) and $\Supp(\Loco{t}) = \Supp(t)$ (\cref{rem:supp-of-localizing}). Thus the injectivity of \eqref{eq:supp-set-gen} is equivalent to
        \[ 
            \forall t_1,t_2 \in \cat T, \Supp(t_1)=\Supp(t_2) \Longrightarrow \Loco{t_1} = \Loco{t_2}.
        \]
    If $\Supp(t_1) = \Supp(t_2)$ then (b) implies that $\Loco{t_1}=\Loco{t_2}$, so we are done.

    $(c)\Rightarrow(a)$: Suppose $0\neq \cat L \subseteq \Loco{\kappaP}$. Then
        \[ 
            \emptyset \neq \Supp(\cat L) \subseteq \Supp(\Loco{\kappaP}) = \big\{\cat P\big\}.
        \]
    Hence $\Supp(\cat L) = \Supp(\Loco{\kappaP})$. So (c) implies that $\cat L=\Loco{\kappaP}$. This establishes the minimality of $\Loco{\kappaP}$. To establish the local-to-global principle, just observe that the support of $\Loco{t \otimes \kappaP \mid \cat P \in \Spc(\cat T^c)}$ is
        \[
            \bigcup_{\cat P \in \Spc(\cat T^c)} \Supp(t \otimes \kappaP)
             = \bigcup_{\cat P \in \Spc(\cat T^c)} 
             \{\cat P\} \cap \Supp(t)
             = \Supp(t)
        \]
    by \cref{rem:supp-of-localizing} and \cref{lem:support-compat-with-finite-loc}. By (c) this implies that we have an equality of localizing $\otimes$-ideals $\Loco{t} = \Loco{t \otimes \kappaP \mid \cat P \in \Spc(\cat T^c)}$.
\end{proof}

\begin{Def}[Stratification]\label{def:stratified}
    A rigidly-compactly generated \mbox{tt-category} $\cat T$ with $\Spc(\cat T^c)$ weakly noetherian is said to be \emph{stratified} when the equivalent conditions (a)--(c) of \cref{thm:equiv-strat} hold. We will also say that $\cat T$ satisfies (or has) ``minimality at $\cat P$'' if $\Loco{\kappaP}$ is a minimal localizing $\otimes$-ideal of $\cat T$ as in part (a). (See \cref{prop:single-prime-local} below.)
\end{Def}

\begin{Rem}\label{rem:min-consequence}
    Note that if $\cat T$ has minimality at $\cat P$ then $\kappaP \in \Loco{t \otimes \kappaP}$ whenever $\cat P \in \Supp(t)$.
\end{Rem}

\begin{Exa}\label{rem:Neeman_example}
    Neeman \cite{Neeman00} provides an example of a (non-noetherian) commutative ring $R$ whose spectrum $\Spec(R)$ consists of a single point and yet whose derived category $\Der(R)$ has many localizing $\otimes$-ideals. As the spectrum is a single point, it is certainly noetherian, hence $\Der(R)$ satisfies the local-to-global principle. However, the category is not stratified. In this example, the Balmer--Favi notion of support merely detects whether an object is zero or not. On the other hand, when~$R$ is noetherian, $\Der(R)$ is always stratified, see \cite{Neeman92a} or \cref{thm:ring-strat} below. 
\end{Exa}

\begin{Rem}
    In analogy with \cref{prop:ltg-permanence}, we will now demonstrate that stratification passes to finite localizations.
\end{Rem}

\begin{Prop}\label{prop:homeo}
    Let $F\colon\cat C \to \cat D$ be a coproduct-preserving tensor-triangulated functor between rigidly-compactly generated tensor-triangulated categories. Suppose that the induced map
        \[
            \varphi\colon\Spc(\cat D^c) \to \Spc(\cat C^c)
        \]
    is injective and that $\Spc(\cat C^c)$ is weakly noetherian. Then $\Spc(\cat D^c)$ is also weakly noetherian and the following hold:
        \begin{enumerate}
            \item If $\cat C$ satisfies minimality at $\varphi(\cat Q)$, then $\cat Q \in \varphi^{-1}(\Supp(c)) \Rightarrow \cat Q \in \Supp(F(c))$ for any $c \in \cat C$.
            \item If $\cat C$ is stratified and $F$ is essentially surjective, then $\cat D$ is also stratified.
        \end{enumerate}
\end{Prop}

\begin{proof}
    Note that \cref{rem:sub-weakly-noetherian} shows that $\Spc(\cat D^c)$ is weakly noetherian. For part~(a), suppose that $\varphi(\cat Q) \in \Supp(c)$. Then 
        \[
            \kappaphiQ \in \Loco{c\otimes \kappaphiQ} \subseteq \Loco{c}
        \]
    since we have minimality at $\varphi(\cat Q)$ (\cref{rem:min-consequence}). Hence, by \cref{prop:ltg-permanence}(a), $\kappaQ \simeq F(\kappaphiQ) \in \Loco{F(c)}$ so that $\cat Q \in \Supp(F(c))$.

    For part (b), first recall from \cref{thm:equiv-strat} that stratification implies minimality at every point. Hence part (a) and \cref{prop:ltg-permanence}(b) imply that 
        \[
            \Supp(F(c)) = \varphi^{-1}(\Supp(c))
        \]
    for any $c \in \cat C$. Now to establish stratification for $\cat D$ it suffices to prove that $d_1 \in \Loco{d_2}$ whenever $\Supp(d_1) \subseteq \Supp(d_2)$. (Recall the proof of (b) $\Rightarrow$ (c) in \cref{thm:equiv-strat}.) Since $F$ is essentially surjective, $d_i = F(c_i)$ for some $c_i \in \cat C$. Then $\varphi^{-1}(\Supp(c_1)) \subseteq \varphi^{-1}(\Supp(c_2))$. Hence $\Supp(c_1) \subseteq \Supp(c_2) \cup (\im(\varphi))^c$ and thus $c_1 \in \Loco{c_2, \{ \kappaP \}_{\cat P \not\in \im(\varphi)}}$ by the local-to-global principle. It follows that 
        \[
            d_1\simeq F(c_1) \in \Loco{F(c_2), \{ F(\kappaP)\}_{\cat P\not\in\im(\varphi)}} = \Loco{F(c_2)} = \Loco{d_2}
        \]
    by \cref{prop:ltg-permanence}(a). This completes the proof.
\end{proof}

\begin{Cor}\label{cor:passes-to-finite}
    Let $\cat T$ be a rigidly-compactly generated tt-category with $\Spc(\cat T^c)$ weakly noetherian. If $\cat T$ is stratified, then so is each finite localization $\cat T(V)$.
\end{Cor}

\begin{Rem}
    In the next section, we will consider the problem of going the other way around --- of deducing stratification for $\cat T$ from its localizations. In particular, we will reduce the problem of stratification to the local categories of \cref{exa:local-categories}.
\end{Rem}

\begin{Rem}
The notion of `stratification' studied in this paper (\cref{def:stratified}) should not be confused with the notion of `stratification' recently considered by Ayala,
Mazel-Gee, and Rozenblyum 
in \cite{AyalaMazelGeeRozenblyum20pp}.
The former provides an approach to the classification of localizing tensor ideals,
while the latter expresses an $\infty$-categorical version of the local-to-global principle 
which is then used to establish reconstruction theorems for presentable stable $\infty$-categories; see also 
\cite{BalchinGreenlees20}.
In particular, the notion of `stratification' in \cite{AyalaMazelGeeRozenblyum20pp} does not require
any minimality of the `stalk categories',
 which for us is a crucial property.
\end{Rem}

\section{Reducing to the local case}\label{sec:reducing-to-local}

Our next task is to reduce the problem of stratification to the local categories.

\begin{Rem}\label{rem:local-kappa}
    Let $\cat T$ be a local rigidly-compactly generated tt-category so that $\Spc(\cat T^c)$ has a unique closed point $\frak m = (0)$. The closed point $\frak m$ is weakly visible if and only if it is visible if and only if the subset $\{\frak m\}$ is Thomason (\cref{rem:closed-weakly-visible}). In this case, 
        \[ 
            \kappam = e_{\{\frak m\}} \otimes f_{Y_{\frak m}}
        \]
    but one readily checks that $Y_{\frak m}=\emptyset$. In other words, $f_{Y_{\frak m}}=\unit$ and we simply have
        \[
            \kappam = e_{\{\frak m\}}.
        \]
\end{Rem}

\begin{Prop}\label{prop:single-prime-local}
    Let $\cat T$ be a rigidly-compactly generated tensor-triangulated category with $\Spc(\cat T^c)$ weakly noetherian and satisfying the local-to-global principle. Then $\cat T$ satisfies minimality at a prime $\cat P \in \Spc(\cat T^c)$ if and only if the local category $\cat T/\langle \cat P\rangle$ satisfies minimality at its unique closed point.
\end{Prop}

\begin{proof}
    First we will prove that $\Loco{\kappaP}$ is minimal in $\cat T$ if $\Loco{\kappam}$ is minimal in the local category $\cat T/\langle \cat P \rangle$ where $\frak m$ denotes  the local category's unique closed point. Let $L\colon \cat T \to \cat T_{\cat P} \coloneqq \cat T/\langle \cat P \rangle$ denote the localization. Then $L(\kappaP) \cong \kappam$ by part (a) of \cref{prop:ltg-permanence}. Our hypothesis is that  $\Loco{\kappam}$ is a minimal localizing $\otimes$-ideal in $\cat T_{\cat P}$. By the relationship between localizing $\otimes$-ideals in $\cat T_{\cat P}$ and~$\cat T$ (\cref{lem:pullback-localizing}), this amounts to
        \[ 
            L^{-1}(0) \subsetneq L^{-1}(\Loco{\kappam})
        \]
    with no localizing  $\otimes$-ideals between them. Using the explicit description of the pullback \eqref{eq:explicit-pullback}, this can be rewritten as
        \[ 
            \Loco{e_{Y_{\cat P}}} \subsetneq \Loco{e_{Y_{\cat P}},U(\kappam)}
        \]
    and since $U(\kappam) \cong UL(\kappaP) \cong f_{Y_{\cat P}} \otimes \kappaP \cong \kappaP$ (\cref{rem:choice-of-thomason}), we have
        \[
			\Loco{e_{Y_{\cat P}}} \subsetneq \Loco{e_{Y_{\cat P}},\kappaP}
        \]
    with no localizing $\otimes$-ideals between them. Now we are ready to establish the minimality of $\Loco{\kappaP}$ in $\cat T$. Consider any nonzero $t \in \Loco{\kappaP}$. We claim that $t \not\in \Loco{e_{Y_{\cat P}}}$. Otherwise, if $t \in \Loco{e_{Y_{\cat P}}}$ then $t \simeq t\otimes e_{Y_{\cat P}}$ and hence $\kappaP \otimes t=0$ (\cref{rem:choice-of-thomason}). But our assumption $t \in \Loco{\kappaP}$ implies $\Supp(t) \subseteq \Supp(\kappaP) = \{ \cat P \}$. So $\kappaP \otimes t=0$ would imply that $\Supp(t)=\emptyset$ and hence that $t=0$ by the local-to-global principle. We conclude that if $0 \neq t \in \Loco{\kappaP}$ then $t \not\in \Loco{e_{Y_{\cat P}}}$. Hence
        \[
            \Loco{e_{Y_{\cat P}}} \subsetneq \Loco{e_{Y_{\cat P}},t} \subseteq \Loco{e_{Y_{\cat P}},\kappaP}
        \]
    so that $\Loco{e_{Y_{\cat P}},t} = \Loco{e_{Y_{\cat P}},\kappaP}$. Therefore $\kappaP \in \Loco{e_{Y_{\cat P}},t}$. It follows that
        \[
            \kappaP\simeq \kappaP\otimes \kappaP \in \Loco{e_{Y_{\cat P}}\otimes \kappaP,t\otimes \kappaP}
        \]
    by \cref{lem:tensor-in}. Since $e_{Y_{\cat P}} \otimes \kappaP=0$, we conclude that 
        \[
            \kappaP \in \Loco{t \otimes \kappaP} = \Loco{t}
        \]
    where the last equality is by the local-to-global principle. Therefore $\Loco{\kappaP} = \Loco{t}$. This establishes that $\Loco{\kappaP}$ is a minimal localizing $\otimes$-ideal in~$\cat T$.

    For the only if part, we need to establish minimality at $\frak m$ in $\cat T_{\cat P}$ assuming minimality at $\cat P$ in $\cat T$. To this end, consider $0 \neq d \in \Loco{\kappam}$. Then $\frak m \in \Supp(d)$ by the local-to-global principle. Write $d=L(t)$ for some $t \in \cat T$. Then $\frak m \in \Supp(L(t))$. By part (b) of \cref{prop:ltg-permanence}, it follows that $\cat P \in \Supp(t)$. Minimality at $\cat P$ then implies that $\kappaP \in \Loco{t}$ (\cref{rem:min-consequence}). Hence $\kappam \cong L(\kappaP) \in \Loco{L(t)}=\Loco{d}$. Therefore $\Loco{\kappam} = \Loco{d}$.
\end{proof}

\begin{Cor}\label{cor:reduced-to-local}
    Let $\cat T$ be a rigidly-compactly generated tensor-triangulated category with $\Spc(\cat T^c)$ weakly noetherian and satisfying the local-to-global principle. Then $\cat T$ is stratified if and only if for each $\cat P \in \Spc(\cat T^c)$ the local category $\cat T/\langle \cat P\rangle$ satisfies minimality at its unique closed point.
\end{Cor}

\begin{proof}
    By \cref{thm:equiv-strat}, stratification of $\cat T$ is equivalent to minimality at each prime $\cat P \in \Spc(\cat T^c)$. The claim thus follows from \cref{prop:single-prime-local}.
\end{proof}

\begin{Rem}\label{rem:cover-suffices}
    Let $V \subseteq \Spc(\cat T^c)$ be the complement of a Thomason subset. For any $\cat P \in V$ we have $\cat T/\langle \cat P\rangle \cong \cat T(V)/\langle \cat P \rangle$ by \cref{prop:localfinitelocalization} and \cref{rem:local-cat-in-open}. That is, the local category of $\cat T$ at $\cat P$ is the same as the local category of $\cat T(V)$ at $\cat P$. Then one consequence of \cref{prop:single-prime-local} is that $\cat T$ satisfies minimality at $\cat P$ if and only if $\cat T(V)$ satisfies minimality at $\cat P \in V$. In order to establish stratification of~$\cat T$ it thus suffices to establish stratification of $\cat T(V)$ for $V$ forming a cover of the spectrum of~$\cat T$:
\end{Rem}

\begin{Cor}\label{cor:cover}
    Let $\cat T$ be a rigidly-compactly generated tensor-triangulated category with $\Spc(\cat T^c)$ weakly noetherian and satisfying the local-to-global principle. Suppose $\Spc(\cat T^c) = \bigcup_{i \in I} V_i$ is a cover by complements of Thomason subsets. Then $\cat T$ is stratified if and only if each finite localization $\cat T(V_i)$ is stratified.
\end{Cor}

\begin{proof}
    The only if part follows from \cref{cor:passes-to-finite}. The rest follows from  \cref{cor:reduced-to-local} keeping in mind \cref{rem:cover-suffices}.
\end{proof}

\begin{Rem}
    If the cover in \cref{cor:cover} is finite, then we may drop the assumption that $\cat T$ satisfies the local-to-global principle. Indeed, in this case it is a consequence of \cref{prop:ltg-cover} and \cref{thm:equiv-strat}.
\end{Rem}

\begin{Exa}
    The first example of a classification of localizing $\otimes$-ideals was for the derived category $\Der(R)$ of a commutative noetherian ring $R$, due to Neeman~\cite{Neeman92a}. He utilizes a notion of support provided by the  residue fields $\kappa(\frak p)$ of the ring~$R$. Since these are not the same as the Balmer--Favi idempotents $\kappaP$ defined using finite localizations, it isn't immediate that $\Der(R)$ is stratified in the sense of \cref{def:stratified}. In fact, in general the two notions of support (defined using the $\kappaP$s and the $\kappa(\frak p)$s, respectively) differ; see \cite{bhs2}. For noetherian rings, however, the universal nature of Balmer--Favi supports (\cref{thm:universal-support}) asserts that stratification for $\Der(R)$ is indeed a consequence of Neeman's result, but it is instructive to give a more direct proof, as explained in \cite[Section~3.4]{Stevenson18}:
\end{Exa}

\begin{Thm}[Neeman]\label{thm:ring-strat}
    The derived category $\Der(R)$ of a commutative noetherian ring is stratified.
\end{Thm}

\begin{proof}
    By \cref{cor:reduced-to-local} and \cref{exa:local-for-DR}, it suffices to assume that $R$ is local and check minimality at the closed point $\frak M \in \Spc(\Der(R)^c)$. By \cref{rem:local-kappa}, $\kappaM = e_{\{\frak M\}}$. Now if the maximal ideal of $R$ is $\frak m=(f_1,\ldots,f_n)$ then 
        \[
            \{\frak M\} = \rho^{-1}(\{ \frak m\}) = \supp(\cone(f_1) \otimes \cdots \otimes \cone(f_n))
        \] 
    is the support of the associated Koszul complex $K_{\frak m}\coloneqq \cone(f_1) \otimes \cdots \otimes \cone(f_n)$. Thus $\Loco{e_{\{\frak M\}}} = \Loco{ K_{\frak m}}$. 
The residue field $\kappa(\frak m)$ is a ``field object'' in $\Der(R)$ in the sense that for any object $X \in \Der(R)$, $\kappa(\frak m)\otimes X$ is a direct-sum of suspensions of copies of $\kappa(\frak m)$.
Since $\kappa(\frak m)\otimes \cone(f) \neq 0$ for any $f\in\frak m$, it follows inductively that
$\kappa(\frak m) \otimes K_{\frak m}$ is a nonzero direct-sum of suspensions of $\kappa(\frak m)$.
Hence $\kappa(\frak m) \in \Loco{K_{\frak m}}$. On the other hand, $\thick\langle\kappa(\frak m)\rangle$ consists precisely of those complexes with finite length homology (see \cite[Example 3.5]{DwyerGreenleesIyengar06} for example). The Koszul complex has finite length homology, so $K_{\frak m} \in \Loco{\kappa(\frak m)}$. We conclude that
        \[ 
            \Loco{\kappa(\frak m)} = \Loco{K_{\frak m}} = \Loco{e_{\{\frak M\}}} = \Loco{\kappaM}.
        \]
    Then consider $X \in \Loco{g(\frak M)}$. If $X \otimes \kappa(\frak m) = 0$ then $\Loco{\kappa(\frak m)}\subseteq \Ker(X\otimes -)$ so that $X \otimes g(\frak M) =0$. Therefore $\Supp(X) = \emptyset$ so that $X=0$ by the local-to-global principle. Thus, if $0 \neq X \in \Loco{g(\frak M)}$ then $X \otimes \kappa(\frak m) \neq 0$
	is a nonzero direct-sum of suspensions of 
	$\kappa(\frak m)$. Hence $\kappa(\frak m) \in \Loco{X}$ so that $\Loco{X} = \Loco{\kappaM}$. We conclude that $\Loco{\kappaM}$ is minimal and this completes the proof.
\end{proof}

\begin{Rem}\label{rem:Dqc(X)}
    For $X$ a quasi-compact and quasi-separated scheme, we have the derived category $\Derqc(X)$ of complexes of $\cat O_X$-modules with quasi-coherent cohomology. It is rigidly-compactly generated and its subcategory of rigid-compact objects $\Derqc(X)^c=\Derperf(X)$ is the derived category of perfect complexes. A fundamental result concerning the Balmer spectrum is that $\Spc(\Derperf(X)) \cong X$. See \cite{Thomason97}, \cite[Theorem~6.3]{Balmer05a}, and \cite[Theorem~9.5]{BuanKrauseSolberg07}. More precisely, the map $X \to \Spc(\Derperf(X))$ induced by the usual notion of cohomological support (via the universal property of the Balmer spectrum) is an isomorphism. Moreover, for any quasi-compact open subset $U \subseteq X \cong \Spc(\Derqc(X)^c)$, the corresponding finite localization of $\Derqc(X)$ can be identified with $\Derqc(U)$. That is, $\Derqc(X)(U) \cong \Derqc(U)$ in the notation of \cref{rem:finite-localizations}.
\end{Rem}

\begin{Cor}[Alonso Tarr\'{\i}o--Jerem\'{\i}as L\'{o}pez--Souto Salorio; Stevenson]\label{cor:noetherian-schemes-strat}
    The derived category $\Derqc(X)$ of a noetherian scheme is stratified.
\end{Cor}

\begin{proof}
    Since $\Spc(\Derqc(X)^c)\cong X$ is noetherian, $\Derqc(X)$ satisfies the local-to-global principle by \cref{thm:local-to-global}. The result then follows from \cref{cor:cover}, \cref{thm:ring-strat} and \cref{rem:Dqc(X)} by taking an open affine cover of $X$.
\end{proof}

\begin{Rem}
    The classification of localizing $\otimes$-ideals was first established by \cite[Corollary~4.13]{AlonsoJeremiasSouto04}. The proof by reduction to the affine case is due to \cite[Corollary~8.13]{Stevenson13}.
\end{Rem}

\section{Extending by finite \'{e}tale morphisms}\label{sec:etale}

Next we want to explain that minimality can be checked (in some situations) after a finite \'{e}tale extension.

\begin{Def}\label{def:finite-etale}
    Let $F\colon \cat C \to \cat D$ be a coproduct-preserving tensor-triangulated functor between rigidly-compactly generated tensor-triangulated categories. We say that $F$ is \emph{finite \'{e}tale} if there is a compact commutative separable algebra~$A$ in~$\cat C$ and an equivalence of tensor-triangulated categories $\cat D\cong A\MMod_{\cat C}$ such that $F\colon\cat C \to \cat D$ becomes the extension-of-scalars functor $F_A\colon \cat C \to A\MMod_{\cat C}$. For background on this notion of ``finite \'{e}tale morphism'' in tensor triangular geometry see \cite{Sanders22,Balmer16a,Balmer16b}.
\end{Def}

\begin{Exa}
	In equivariant examples, restriction to a finite-index subgroup is often  a finite \'{e}tale morphism in the sense of \cref{def:finite-etale}, as explained in~\cite{BalmerDellAmbrogioSanders15}. This is the case, for example, for the restriction functor $\res^G_H \colon \SH_G\to \SH_H$ on the equivariant stable homotopy category. The same is true for categories of spectral Mackey functors considered in \cref{part:spectral-mackey} (see \cref{lem:resfiniteetale}). Further examples are discussed in \cite[Section 5]{Sanders22}. In particular, \cite[Example~5.12]{Sanders22} shows that finite \'{e}tale morphisms are ``local in the target'':
\end{Exa}

\begin{Lem}\label{lem:restriction-is-finite-etale}
    Let $F\colon \cat C \to \cat D$ be a finite \'{e}tale morphism of rigidly-compactly generated tensor-triangulated categories and let
        \[
            \varphi\colon \Spc(\cat D^c) \to \Spc(\cat C^c)
        \]
    denote the induced map. For any Thomason subset $Y \subseteq \Spc(\cat C^c)$ with complement $V\coloneqq \Spc(\cat C^c)\setminus Y$, the induced functor (\cref{prop:restrict-in-target})
        \[
            \cat C(V) \to \cat D(\varphi^{-1}(V)) 
        \]
    is also finite \'{e}tale.
\end{Lem}

\begin{Thm}\label{thm:finite-etale}
    Let $F\colon \cat C \to \cat D$ be a finite \'{e}tale morphism of rigidly-compactly generated tensor-triangulated categories. Assume that both categories have noetherian spectrum, and let
        \[
            \varphi\colon \Spc(\cat D^c) \to \Spc(\cat C^c)
        \]
    denote the induced map. If $\cat P \in \Spc(\cat D^c)$ is a point such that $\varphi^{-1}(\{\varphi(\cat P)\}) = \{\cat P\}$ then minimality at $\cat P$ in $\cat D$ implies minimality at $\varphi(\cat P)$ in $\cat C$.
\end{Thm}

\begin{proof}
	First we prove the theorem under the additional assumption that $\cat C$ is local and $\varphi(\cat P)$ is the closed point in $\Spc(\cat C^c)$. It follows that the unique point~$\cat P$ in its preimage is a closed point of $\Spc(\cat D^c)$ (although it need not be the only one). We then proceed with the following observation. Recall that the closed point in the local category $\cat C^c$ is $(0)$. We claim that $F$ is conservative on big objects supported on the closed point $\frak m \in \Spc(\cat C^c)$. Indeed, let $0 \neq t \in \cat C$ with $\Supp(t) \subseteq \{\frak m\}$. If $F(t)=0$ then $A \otimes t \simeq UF(t) = 0$ where $U$ denotes the right adjoint of $F$ and $A=U(\unit)$ is the compact commutative separable algebra providing the finite \'{e}tale morphism. The fact that $A$ is compact means that $\emptyset = \Supp(A \otimes t) = \supp(A) \cap \Supp(t)$ by the half-tensor theorem for Balmer--Favi supports (see \cref{lem:half-tensor}). Moreover, compactness of $A$ implies that $\supp(A)$ is specialization closed and thus contains~$\frak m$. Hence $\Supp(t) =\emptyset$ and therefore $t=0$.

    Recall from \cref{rem:local-kappa} that in a local category, $g(\frak m) = e_{\{\frak m\}}$ is the left idempotent associated to the Thomason subset $Y\coloneqq \{\frak m\}$. By our hypothesis, the preimage $\varphi^{-1}(Y) = \{\cat P\}$ is a single closed point. It follows that $F(e_{\{\frak m\}}) = e_{\{\cat P\}}$. With this in hand, we can prove that minimality at $\frak m$ follows from minimality of $\cat P$. Consider any $0\neq t \in \Loco{e_{\{\frak m\}}}$. Then $0 \neq F(t)$ by the conservativity claim above. Also $F(t) \in \Loco{F(e_{\{\frak m\}})} = \Loco{e_{\{\cat P\}}}$. Now $\Supp(e_{\{\cat P\}}) = \{\cat P\}$ so
        \[
            0 \neq \Loco{e_{\{\cat P\}}} = \Loco{e_{\{\cat P\}}\otimes \kappaP} \subseteq \Loco{\kappaP}
        \]
    by the local-to-global principle in $\cat D$. Therefore, minimality at $\cat P$ implies that $\Loco{e_{\{\cat P\}}} = \Loco{\kappaP}$. Minimality of this then implies that
        \[
            \Loco{F(t)} = \Loco{e_{\{\cat P\}}} = \Loco{\kappaP}.
        \]
    Therefore $F(e_{\{\frak m\}}) = e_{\{\cat P\}} \in \Loco{F(t)}$. Then consider $U^{-1}(\Loco{UF(t)})$.
	It is a localizing subcategory of $\cat D$ since $U$ preserves coproducts.
	The separable Neeman--Thomason theorem due to Balmer \cite[Theorem~4.2]{Balmer16a} establishes that every compact object of $\cat D$ is a direct summand of the image $F(c)$ of a compact object $c \in \cat C^c$. It follows, using the projection formula (\cref{rem:geometric-functor}) that  $U^{-1}(\Loco{UF(t)})$ is closed under tensoring with any compact object, and hence is closed under tensoring with any object (\cref{rem:suffice-tensor-compact}), so that it is a localizing $\otimes$-ideal. It contains
$F(t)$ and thus also contains
	$F(e_{\{\frak m\}})$. Hence we conclude that $UF(e_{\{\frak m\}}) \in \Loco{UF(t)}$. That is, $A \otimes e_{\{\frak m\}} \in \Loco{A\otimes t}$. Then
        \[
            \Loco{e_{\{\frak m\}}} = \Loco{A\otimes e_{\{\frak m\}}} \subseteq \Loco{A\otimes t} \subseteq \Loco{t}
        \]
    where the first equality follows from  \cref{lem:big-supp-in-small} since $\frak m \in \supp(A)$. This establishes the special case of the proposition when $\cat C$ is local and $\varphi(\cat P)$ is the unique closed point.

    Now we turn to the general statement. Recall from \cref{exa:local-categories} that the localization $\cat C \to \cat C/\langle \cat \varphi(\cat P)\rangle$ is the finite localization associated to the Thomason subset $Y_{\varphi(\cat P)}$ whose complement $Y_{\varphi(\cat P)}^c$ consists of all generalizations of $\varphi(\cat P)$. The preimage $\varphi^{-1}(Y_{\varphi(\cat P)})$ need not coincide with $Y_{\cat P}$. Indeed, $Y_{\cat P}^c = \SET{\cat Q}{\cat Q \supseteq \cat P}$ while $V\coloneqq \varphi^{-1}(Y_{\varphi(\cat P)})^c = \SET{\cat Q}{\varphi(\cat Q) \supseteq \varphi(\cat P)}$. In any case, we have the induced functor as from \cref{prop:restrict-in-target}:
        \[
            \begin{tikzcd}
            \cat C \ar[d] \ar[r,"F"] & \cat D \ar[d]  \ar[dr,bend left, "(-)_{\cat P}"]\\
            \cat C/\langle \varphi(\cat P)\rangle \ar[r,"F"] & \cat D(V) \ar[r,"(-)_{\cat P}"] & \cat D/\langle \cat P \rangle
            \end{tikzcd}
        \]
    and the right triangle of functors are all finite localizations, while the induced functor $F$ is finite \'{e}tale by \cref{lem:restriction-is-finite-etale}. On spectra the bottom line is (reversing directions)
        \[
            \SET{\cat Q \subset \cat D^c}{\cat Q\supseteq \cat P} \hookrightarrow \SET{\cat Q \subset \cat D^c}{\varphi(\cat Q) \supseteq \varphi(\cat P)}  \xra{\varphi} \SET{\cat Q' \subset \cat C^c}{\cat Q' \supseteq \varphi(\cat P)}. 
        \]
    If $\cat P$ is minimal in $\cat D$ then it is minimal in the local category and hence in $\cat D(V)$ (recall \cref{prop:single-prime-local}). The lower $F$ satisfies the hypotheses of the local case of the proposition we have already proved, so we have established that $\varphi(\cat P)$ is minimal in the  local category $\cat C/\langle \varphi(\cat P)\rangle$, and hence $\varphi(\cat P)$ is minimal in $\cat C$ (again invoking \cref{prop:single-prime-local}).
\end{proof}

\begin{Rem}
    We'll use \cref{thm:finite-etale} in \cref{sec:strat-for-spectral-mack} to establish stratification for certain equivariant examples by passing to smaller subgroups.
\end{Rem}

\part{Stratification: universality and consequences}\label{part:consequences}

\section{Universal weakly noetherian stratification}\label{sec:universal}
We now establish that the approach to stratification described above, which uses the Balmer--Favi notion of support, is in a certain sense the universal approach to stratification, at least in weakly noetherian contexts. Throughout this section, we will assume that $\cat T$ is a rigidly-compactly generated tensor-triangulated category. We begin with an axiomatization of the properties of the Balmer--Favi notion of support recorded in \cref{rem:support-axioms}.

\begin{Def}\label{def:axiomaticsupp}
    Let $X$ be a topological space and let $\sigma\colon \cat T \to \mathcal{P}(X)$ be a function defined on the objects of $\cat T$. The pair $(X,\sigma)$ is called a \emph{theory of support} if it satisfies the following axioms:
    \begin{enumerate}
        \item $\sigma(0) = \emptyset$ and $\sigma(\unit)=X$;
        \item $\sigma(\Sigma t) = \sigma(t)$ for every $t \in\cat T$;
        \item $\sigma(c) \subseteq \sigma(a) \cup \sigma(b)$ for any exact triangle $a \to b \to c \to \Sigma a$ in~$\cat T$;
        \item $\sigma(\coprod_{i\in I} t_i) = \bigcup_{i \in I} \sigma(t_i)$ for any set of objects $t_i$ in $\cat T$;
        \item $\sigma(t_1 \otimes t_2) \subseteq \sigma(t_1) \cap \sigma(t_2)$ for any $t_1,t_2 \in \cat T$.
    \end{enumerate}
    We also refer to $\sigma$ as a \emph{support function}.
\end{Def}

\begin{Rem}
    Note that we have not yet made use of the topology on $X$. Nevertheless, it will play an important role below.
\end{Rem}

\begin{Def}\label{def:axiomaticstratifying}
    We extend a support function $\sigma$ to a collection of objects $\cat S \subset \cat T$ by defining $\sigma(\cat S) \coloneqq \bigcup_{s \in \cat S}\sigma(s)$. The axioms ensure that $\sigma(\Loco{\cat E}) = \bigcup_{t \in \cat E} \sigma(t)$ for any set of objects $\cat E \subset \cat T$. In particular, $\sigma$ induces a map 
    \begin{equation}\label{eq:axiomaticidealclassification}
        \big\{ \text{localizing $\otimes$-ideals of $\cat T$} \big\} \longrightarrow \big\{ \text{subsets of $X$}\big\}.
    \end{equation}
    If this map is a bijection, then the support theory $(X,\sigma)$ is said to be \emph{stratifying}.
\end{Def}

\begin{Exa}
    If the Balmer spectrum $\Spc(\cat T^c)$ is weakly noetherian then the Balmer--Favi notion of support provides a support theory $(\Spc(\cat T^c),\Supp)$ for $\cat T$.
\end{Exa}

\begin{Thm}\label{thm:universal-support}
    Let $\cat T$ be a rigidly-compactly generated tensor-triangulated category. Let $X$ be a spectral space equipped with a notion of support $\sigma(t) \subseteq X$ for the objects~$t \in \cat T$ satisfying axioms (a)--(e) of \cref{def:axiomaticsupp}. Suppose that this notion of support stratifies $\cat T$ in the sense that the map
        \begin{align}\label{eq:sigma-bijection}
            \begin{split}
                \big\{ \text{localizing $\otimes$-ideals of $\cat T$} \big\} &\to \big\{ \text{subsets of $X$} \big\}
            \end{split}
        \end{align}
    defined by $\cat L \mapsto \sigma(\cat L) \coloneqq \bigcup_{t \in \cat L} \sigma(t)$ is a bijection. Then the following are equivalent:
        \begin{enumerate}
            \item The subset $\sigma(x) \subseteq X$ is Thomason closed for each compact object $x \in \cat T^c$ and every Thomason closed subset of $X$ arises in this way.
            \item Under the bijection, the compactly generated localizing $\otimes$-ideals of $\cat T$ correspond to the Thomason subsets of $X$.
            \item There is a unique homeomorphism $f\colon X\xra{\sim} \Spc(\cat T^c)$ such that $\sigma(x) = f^{-1}(\supp(x))$ for every compact $x \in \Spc(\cat T^c)$.
        \end{enumerate}
    If these conditions hold and the spectral space $X$ is weakly noetherian then  under the homeomorphism $X\xra{\sim} \Spc(\cat T^c)$ the notion of support $\sigma$ coincides with the Balmer--Favi notion of support $\Supp$.
\end{Thm}

\begin{proof}
    The axioms on~$\sigma$ imply that for each subset $Y \subseteq X$, 
        \[
            \Psi(Y) \coloneqq \SET{t \in \cat T}{{\sigma(t) \subseteq Y}}
        \]
    is a localizing $\otimes$-ideal of $\cat T$. One readily checks that $\cat L \subseteq \Psi(\sigma(\cat L))$ for any localizing $\otimes$-ideal $\cat L$. Hence $\sigma(\cat L) \subseteq \sigma(\Psi(\sigma(\cat L))) \subseteq \sigma(\cat L)$ so that $\sigma(\cat L) = \sigma(\Psi(\sigma(\cat L)))$ and therefore $\cat L = \Psi(\sigma(\cat L))$. In other words, $\Psi$ provides the inverse to the bijection~\eqref{eq:sigma-bijection}. With these comments in hand, we proceed with the proof. First we establish that (a) and (b) are equivalent.

    (a)$\Rightarrow$(b): The first part of the hypothesis implies that if $\cat L$ is compactly generated then $\sigma(\cat L)$ is a Thomason subset. On the other hand, if $Y = \bigcup_{i\in I} Z_i$ is a Thomason subset (a.k.a.~union of Thomason closed sets) then $Z_i = \sigma(x_i)$ for compact $x_i \in \cat T^c$ by the second part of the hypothesis. Hence $Y=\sigma(\Loco{x_i \mid i \in I})$ is the support of a compactly generated $\otimes$-ideal. This establishes (b).

    (b)$\Rightarrow$(a): First we prove that every Thomason closed subset $Z$ is of the form~$\sigma(x)$ for a compact object $x \in \cat T^c$. Well, $Z = \sigma(\Psi(Z))$ and by hypothesis $\Psi(Z)$ is compactly generated. Hence
        \[ 
            Z = \sigma(\Psi(Z)) = \bigcup_{t \in \cat T^c \cap \Psi(Z)} \sigma(t).
        \]
    In the Hochster-dual topology, this is an open cover of a quasi-compact open subset (see \cref{rem:hochster-dual}). Hence $Z=\sigma(x_1)\cup \cdots \cup\sigma(x_n)$ for finitely many $x_1,\ldots, x_n \in \cat T^c \cap \Psi(Z)$. That is, $Z = \sigma(x_1 \oplus \cdots \oplus x_n)$ is the support of a compact object.

    To establish (a) it remains to show that $\sigma(x)$ is closed for every compact object~$x$. By hypothesis, $\sigma(x)$ is a Thomason subset, that is, a union of Thomason closed subsets. By what we have just proved, we can write each Thomason closed subset as the support of a compact object. Thus $\sigma(x) = \bigcup_{i \in I} \sigma(x_i)$ for some compact $x_i$. This implies by the stratification hypothesis that we have an equality of (compactly generated) localizing $\otimes$-ideals $\Loco{x} = \Loco{x_i \mid i \in I}$ which in turn implies that
        \[
            \thickt{x} = \cat T^c \cap \Loco{x} = \cat T^c \cap \Loco{x_i \mid i \in I} = \thickt{x_i \mid i\in I}.
        \]
    It follows that $x \in \thickt{x_1,\ldots,x_n}$ for a finite collection of the $x_i$. Hence $\thickt{x} = \thickt{x_1 \oplus \cdots \oplus x_n}$ and $\sigma(x) = \sigma(x_1) \cup \cdots \cup \sigma(x_n)$ is closed. This establishes (a).

    We have established that (a) and (b) are equivalent. Note that (c) trivially implies (a) since the statements in (a) are basic properties of the universal support~\eqref{eq:supp}. We now prove that (a) and (b) imply (c).

    To prove (c) we first establish that $\sigma$ satisfies the tensor product formula for compact objects: $\sigma(x \otimes y) = \sigma(x) \cap \sigma(y)$. First note that if $x$ is a compact object then $\Loco{x} = \Psi(\sigma(x))$ and $\thickt{x} = \Loco{x} \cap \cat T^c = \Psi(\sigma(x)) \cap \cat T^c$. Then for any two compact objects $x$ and $y$, we have
        \begin{align*}
            \Psi(\sigma(x)\cap\sigma(y)) \cap \cat T^c &= \SET{a \in \cat T^c}{\sigma(a) \subseteq \sigma(x) \cap \sigma(y)}\\
            &= \Psi(\sigma(x)) \cap \Psi(\sigma(y)) \cap \cat T^c \\
            &= \thickt{x} \cap \thickt{y}.
        \end{align*}
    Now the universal notion of support in $\Spc(\cat T^c)$ has the tensor product property: $\supp(x\otimes y) = \supp(x) \cap \supp(y)$. Hence, by the classification of thick $\otimes$-ideals, 
        \[ 
            \thickt{x\otimes y} = \thickt{x} \cap \thickt{y} = \Psi(\sigma(x) \cap \sigma(y)) \cap \cat T^c.
        \]
    Therefore
        \[
            \Loco{x\otimes y} = \Loc\langle \thickt{x\otimes y}\rangle = \Loc\langle \Psi(\sigma(x) \cap \sigma(y)) \cap \cat T^c\rangle = \Psi(\sigma(x) \cap \sigma(y))
        \]
    since $\Psi(\sigma(x) \cap \sigma(y))$ is compactly generated by hypothesis (since the intersection of the Thomason subsets $\sigma(x)$ and $\sigma(y)$ is Thomason). This establishes that 
        \[
            \sigma(x \otimes y) = \sigma(x) \cap \sigma(y)
        \]
    for any two compact objects $x,y \in \cat T^c$.

    This tensor product formula, combined with the axioms for $\sigma$ and the fact that~$\sigma(x)$ is Thomason closed for compact $x$, implies that the restriction of $\sigma$ to $\cat T^c$ is a support datum in the sense of \cite[Definition 3.1]{Balmer05a}. Hence by the universal property of the Balmer spectrum \cite[Theorem~3.2]{Balmer05a}, there is a unique continuous map $f\colon X \to \Spc(\cat T^c)$ such that $f^{-1}(\supp(x)) = \sigma(x)$ for each compact $x$. Moreover, by our hypotheses we have an induced bijection between the thick $\otimes$-ideals of~$\cat T^c$ and the Thomason subsets of $X$ given by $\cat C \mapsto \sigma(\Loco{\cat C}) = \sigma(\cat C)$ with inverse $Y \mapsto \Psi(Y) \cap \cat T^c$. Therefore, $\sigma$ defines a classifying support datum in the sense of \cite[Section 6]{BuanKrauseSolberg07} for the lattice of thick $\otimes$-ideals of $\cat T^c$ (cf.~\cite[Lemma 7.10 and Definition~5.2]{BuanKrauseSolberg07}). Hence by \cite[Corollary~6.2]{BuanKrauseSolberg07}, the map $f\colon X\to \Spc(\cat T^c)$ is a homeomorphism. This establishes (c).

    Under the homeomorphism $X \xra{\sim} \Spc(\cat T^c)$ of part (c), our support function~$\sigma$ induces a support function in $\Spc(\cat T^c)$ which coincides with the universal support on compact objects. It remains to establish that this induced notion of support in $\Spc(\cat T^c)$ coincides with the Balmer--Favi notion of support when $X$ is weakly noetherian. This will follow from \cref{prop:uniqueness-of-support}  if we can establish that the induced notion of support satisfies properties  \hyperref[it:unique-support-a]{(a)} and \hyperref[it:unique-support-b]{(b)} of that proposition. Property~\hyperref[it:unique-support-a]{(a)} is immediate since the stratification hypothesis ensures that $t=0$ if $\sigma(t) = \emptyset$. To invoke the proposition it thus remains to verify property~\hyperref[it:unique-support-b]{(b)}. By hypothesis, a Thomason subset $Y$ corresponds to  the compactly generated $\Psi(Y) = \Loc\langle \cat T^c \cap \Psi(Y)\rangle =\Loc\langle x \in \cat T^c \mid \sigma(x) \subseteq Y \rangle$. But, since $\sigma(x) = \supp(x)$ for compact~$x$, this is $\Loc\langle x \in \cat T^c \mid \supp(x) \subseteq Y\rangle = \Loco{e_Y}$. Thus $\sigma(e_Y) = \sigma(\Psi(Y))=Y$. Now under the stratification, $\sigma(e_Y) \cap \sigma(f_Y) = \emptyset$ corresponds to $\Loco{e_Y} \cap \Loco{f_Y} = 0$. Hence $Y \cap \sigma(f_Y) = \emptyset$ so that $\sigma(f_Y) \subseteq Y^c$. Moreover, from the exact triangle $e_Y \to \unit \to f_Y \to \Sigma e_Y$ and the axiom that $\sigma(\unit)$ is the whole space we conclude that $\sigma(f_Y) = Y^c$. Now consider an arbitrary object $t \in \cat T$ and the exact triangle $e_Y \otimes t \to t \to f_Y \otimes t\to\Sigma e_Y \otimes t$. Our basic axioms for $\sigma$ (\cref{def:axiomaticsupp}) imply that $\sigma(e_Y \otimes t) \subseteq \sigma(e_Y) \cap \sigma(t) = Y \cap \sigma(t)$ and $\sigma(f_Y \otimes t) \subseteq \sigma(f_Y) \cap \sigma(t) = Y^c \cap \sigma(t)$. It then follows from $\sigma(t) \subseteq \sigma(e_Y \otimes t) \cup \sigma(f_Y \otimes t)$ that we in fact have equalities $\sigma(e_Y \otimes t) = Y \cap \sigma(t)$ and $\sigma(f_Y \otimes t) = Y^c \cap \sigma(t)$. We thus conclude that $\sigma = \Supp$ by invoking \cref{prop:uniqueness-of-support}. This completes the proof.
\end{proof}

\begin{Rem}\label{rem:uniquenessreinterpretation}
    We may reinterpret \cref{thm:universal-support} as follows. Suppose $(X,\sigma)$ is a stratifying support theory for $\cat T$. If the collection $\sigma(x)^c \subseteq X$ for $x \in \cat T^c$ forms a basis of quasi-compact open sets for the topology of $X$, which is spectral and weakly noetherian, then there is a canonical identification $(X,\sigma) \cong (\Spc(\cat T^c),\Supp)$ and $\cat T$ is stratified in the sense of \cref{def:stratified}. 
\end{Rem}

\begin{Rem}
    The philosophical point of \cref{thm:universal-support} is as follows. Suppose you have a rigidly-compactly generated category $\cat T$ whose spectrum $\Spc(\cat T^c)$ is weakly noetherian, so that the Balmer--Favi notion of support is defined. If the localizing $\otimes$-ideals of $\cat T$ can be stratified by a reasonable notion of support on $\cat T$ in a way that is compatible with the classification of thick $\otimes$-ideals of $\cat T^c$ then the Balmer--Favi notion of support will do the job. Therefore, one might as well use the Balmer--Favi notion of support to begin with. In this sense, the Balmer--Favi notion of support is the ``universal choice'' for stratifying categories whose Balmer spectrum is weakly noetherian.
\end{Rem}

\begin{Rem}\label{rem:BIK}
    We now briefly recall the approach to the classification of localizing $\otimes$-ideals developed by Benson--Iyengar--Krause (BIK) in \cite{BensonIyengarKrause08,BensonIyengarKrause11b}, which in turn is inspired by the work of Hovey--Palmieri--Strickland \cite{HoveyPalmieriStrickland97}. In fact, their theory applies more generally to triangulated categories without monoidal structures, but we will not make use of this here. The starting point of their framework is an auxiliary action of a (graded) commutative noetherian ring $R$ on the category $\cat T$. In parallel to the Balmer--Favi idempotents $\kappaP$, they construct smashing endofunctors $\Gamma_{\frak p}\colon \cat T \to \cat T$ for any (homogeneous) prime ideal $\frak p \in \Spec(R)$. The BIK notion of support is then defined as
        \[
            \SuppBIK(t) \coloneqq \SET{\frak p \in \Spec(R)}{\Gamma_{\frak p}t \neq 0}.
        \]
    We also define $\supp_R(\cat T) \coloneqq \SuppBIK(\unit) \subseteq \Spec(R)$ and note that $\SuppBIK(t) \subseteq \supp_R(\cat T)$ for all $t \in \cat T$. The category $\cat T$ is said to be stratified with respect to the action of $R$ in the sense of BIK if the $\Gamma_{\frak p}\cat T$ are trivial or minimal localizing $\otimes$-ideals of $\cat T$ for any $\frak p$; because $R$ is assumed to be (graded) noetherian, the local-to-global principle automatically holds in their setting. In the language of \cref{def:axiomaticsupp}, the pair $(\supp_R(\cat T),\Supp_{\text{BIK}})$ is a support theory for $\cat T$, and $\cat T$ is stratified by the action of $R$ in the sense of BIK if $(\supp_R(\cat T),\Supp_{\text{BIK}})$ is stratifying in the sense of \cref{def:axiomaticstratifying}.
	Finally, we recall that in \cite{BensonIyengarKrause11b} the category~$\cat T$ is said to be noetherian if the graded $R$-module $\End^*_{\cat T}(c)$ is finitely generated for each compact object $c \in \cat T^c$.
\end{Rem}

\begin{Cor}\label{cor:BIK}
    Let $\cat T$ be a rigidly-compactly generated tensor-triangulated category which, in the terminology of \cite{BensonIyengarKrause11b}, is noetherian and stratified by the action of a graded-noetherian ring $R$. Then the BIK space of supports $\supp_R(\cat T)$ is homeomorphic to $\Spc(\cat T^c)$ and the BIK notion of support in $\supp_R(\cat T)$ coincides with the Balmer--Favi notion of support in $\Spc(\cat T^c)$.
\end{Cor}

\begin{proof}
    It is established in \cite[Sections 5--6]{BensonIyengarKrause08} that the BIK notion of support satisfies axioms (b), (c), and (d) of \cref{def:axiomaticsupp}. Moreover, their restriction to $\supp_R(\cat T)$ rather than using the whole space $\Spec(R)$ ensures that axiom (a) holds. Since the functors $\Gamma_{\frak p}$ are smashing, it follows that axiom (e) holds. Also by \cite[Theorem~5.5]{BensonIyengarKrause08}, the space of supports $\supp_R(\cat T) = \SuppBIK(\unit)$ is a closed subspace of $\Spec(R)$, hence is itself a noetherian spectral space. The claim then follows from \cref{thm:universal-support} and \cite[Theorem~6.1]{BensonIyengarKrause11b}.
\end{proof}

\begin{Exa}[BIK-stratified categories]\label{exa:BIK-stratified}
    Given \cref{cor:BIK}, we can restate all the BIK stratification results in the literature in terms of stratification via the Balmer--Favi notion of support. For example, the stable module category~$\StMod_{kG}$ of a finite group $G$ over a field $k$ is stratified by the Balmer--Favi notion of support \cite{BensonIyengarKrause11a}. More generally, this holds when $G$ is a finite group scheme over~$k$~\cite{BensonIyengarKrausePevtsova18}. There are many more examples, such as module categories of cochain algebras \cite{BarthelCastellanaHeardValenzuela19,BarthelCastellanaHeardValenzuela22} or affine weakly regular tensor-triangulated categories \cite{DellAmbrogioStanley16}.
\end{Exa}

\begin{Rem}
    \Cref{thm:universal-support} and \cref{cor:BIK} are related to  results of Dell'Ambrogio \cite[Theorem 3.1 and Theorem 3.14]{DellAmbrogio10}. From one perspective, our results establish that (under weak hypotheses) a stratifying support theory automatically satisfies the remainder of the hypotheses used in the cited theorems.
\end{Rem}

\begin{Cor}
	Suppose the graded-commutative ring $R^\bullet \coloneqq \Hom_{\cat T}(\unit,\Sigma^{\bullet} \unit)$ is graded-noetherian. If, in the terminology
	of \cite{BensonIyengarKrause11b}, $\cat T$ is noetherian and stratified by
	the canonical action of $R^\bullet$
	then the comparison map $\rho^\bullet\colon  \cat \Spc(\cat T^c) \to \Spec(R^\bullet)$ is a homeomorphism.
\end{Cor}

\begin{proof}
    Let $R\coloneqq R^\bullet$ for notational simplicity 
	and write $\supp_R \coloneqq \SuppBIK$ for the BIK support (\cref{rem:BIK}) for the canonical action of $R$ on $\cat T$. We refer to \cite{Balmer10b} for the construction of the comparison map $\rho^\bullet$. By \cref{cor:BIK}, the canonical continuous map $\supp_R(\cat T) \to \Spc(\cat T^c)$ is a homeomorphism. Moreover, the comparison map $\rho^\bullet\colon \Spc(\cat T^c) \to \Spec(R)$ is surjective by \cite[Theorem~7.3]{Balmer10b}. The composite map $\supp_R(\cat T)\xra{\sim} \Spc(\cat T^c) \twoheadrightarrow \Spec(R)$ sends $\frak p \in \supp_R(\cat T)$ to the homogeneous prime ideal $\SET{f \in R }{\frak p \in \supp_R(\cone(f))} \in \Spec(R)$. By \cite[Lemma~2.6]{BensonIyengarKrause11b}, 
        \[
            \supp_R(\cone(f)) = V(f) \cap \supp_R(\unit) = V(f) \cap \supp_R(\cat T),
        \]
    where $V(f)=\SET{\frak p \in \Spec(R)}{f \in \frak p}$. Hence for any $\frak p \in \supp_R(\cat T)$ we have
        \begin{align*}
            \SET{f \in R }{\frak p \in \supp_R(\cone(f))} &= \SET{f \in R}{\frak p \in V(f) \cap \supp_R(\cat T)} \\
                &= \SET{f \in R}{\frak p \in V(f)} = \frak p.
        \end{align*}
    It follows that the comparison map $\Spc(\cat T^c) \to \Spec(R)$ is injective, hence a continuous bijection.
	The assumption that $\cat T$ is noetherian \cite[p.~642]{BensonIyengarKrause11b}
	means that $\cat T^c$ is End-finite in the terminology of \cite[Definition~2.6]{Lau22pp}. Hence the continuous bijection $\rho^\bullet$ is 
	automatically 
	a homeomorphism
	by \cite[Corollary 2.8]{Lau22pp}.
	Note that we cannot just invoke \cite[Lemma~2.2]{DellAmbrogioStanley16} due to
	\cite{DellAmbrogioStanley16Erratum}.
\end{proof}

\begin{Exa}[Noetherian schemes]
    Recall that the derived category $\Derqc(X)$ of a noetherian scheme $X$ is stratified by the Balmer--Favi notion of support (\cref{cor:noetherian-schemes-strat}). This also follows from \cref{thm:universal-support} and \cite{AlonsoJeremiasSouto04}. Indeed, it is proved in \cite[Corollary~4.13]{AlonsoJeremiasSouto04} that the usual cohomological support
        \[
            \sigma(E) \coloneqq \SET{ x \in X }{ E \otimes_{\cat O_X}^{\bbL} k(x) \neq 0 \text{ in $\Der(k(x))$}}
        \]
    provides a stratifying support theory for $\Derqc(X)$ taking values in $X$, and the remaining hypotheses of \cref{thm:universal-support} are well-known (see \cite[Lemma~3.3 and Lemma~3.4]{Thomason97}). Note that $\Derqc(X)$ is an example of a stratified category which is not, in general, stratified in the sense of BIK.
\end{Exa}

\begin{Rem}
    Throughout the rest of this paper when we say that a category~$\cat T$ is \emph{stratified} we mean in the sense of \cref{def:stratified}. In other words, it means that $\Spc(\cat T^c)$ is weakly noetherian and $\cat T$ is stratified by the Balmer--Favi notion of support. It is a property of the category $\cat T$.
\end{Rem}

\section{The tensor product formula and Bousfield classes}\label{sec:tensor-and-bousfield}

We now mention a couple of known consequences of stratification. Further consequences will be given in \cref{sec:gentelescope} and the companion paper \cite{bhs2}. 

\begin{Rem}\label{rem:half-tensor-to-full-tensor}
    Recall from \cref{lem:half-tensor} that the Balmer--Favi notion of support satisfies a ``half tensor-product theorem'':
        \[
            \Supp(t \otimes x) = \Supp(t) \cap \supp(x)
        \]
    for any $t \in \cat T$ and compact $x \in \cat T^c$. In the presence of stratification, this can be improved to a full tensor-product theorem:
\end{Rem}

\begin{Thm}\label{thm:tensor_product}
    Let $\cat T$ be a rigidly-compactly generated tt-category with $\Spc(\cat T^c)$ weakly noetherian. If $\cat T$ is stratified then
        \[
            \Supp(t_1 \otimes t_2) = \Supp(t_1) \cap \Supp(t_2)
        \]
    for any $t_1,t_2 \in \cat T$.
\end{Thm}

\begin{proof}
    The inclusion $\Supp(t_1 \otimes t_2) \subseteq \Supp(t_1) \cap \Supp(t_2)$ is a basic property of the support (\cref{rem:support-axioms}). For the reverse inclusion, first note that when stratification holds, $\cat P \in \Supp(t)$ implies that $\kappaP \in \Loco{\kappaP\otimes t}$ (\cref{rem:min-consequence}). Thus, if $\cat P \in \Supp(t_1)$ then $\kappaP \in \Loco{\kappaP \otimes t_1}$. Applying \cref{lem:tensor-in} we then have $\kappaP \otimes t_2 \in \Loco{\kappaP \otimes t_1 \otimes t_2}$. Therefore, if $\cat P \in \Supp(t_2)$ as well, then $\kappaP \in \Loco{\kappaP\otimes t_2}$ so that $\kappaP \in \Loco{\kappaP \otimes t_1 \otimes t_2}$. It follows that $\kappaP \otimes t_1 \otimes t_2 \ne 0$, and that $\cat P \in \Supp(t_1 \otimes t_2)$ as required. 
\end{proof}

\begin{Rem}
    Without the stratification condition, the tensor-product theorem does not hold in general, even when $\Spc(\cat T^c)$ is noetherian. For example, in the derived category $\Der(R)$ of the ring considered by Neeman (\cref{rem:Neeman_example}), there exists a nonzero object $I$ with $I \otimes I = 0$ \cite[Theorem C]{DwyerPalmieri08}. Since the detection property holds in this example, the existence of such an $I$ implies that the tensor-product theorem cannot hold in $\Der(R)$. 
\end{Rem}

\begin{Def}[Bousfield classes]
    For a fixed $t \in \cat T$, recall that an object $s \in \cat T$ is \emph{$t$-acyclic} if $s \otimes t = 0$. The \emph{Bousfield class} of $t$ is the collection of all $t$-acyclic objects:
        \[
            A(t) \coloneqq \SET{s \in \cat T}{s \otimes t = 0}.
        \]
    Note that $A(t)$ is a localizing $\otimes$-ideal of $\cat T$.
\end{Def}
    
\begin{Rem}
    By a theorem of Iyengar and Krause \cite[Theorem~3.1]{IyengarKrause13}, if $\cat T$ is compactly generated, then there is only a set of Bousfield classes, which we denote
        \[
            A(\cat T) \coloneqq \SET{ A(t) }{t \in \cat T}.
        \]
\end{Rem}

\begin{Rem}
    There is a partial order on the set $A(\cat T)$ of Bousfield classes given by
        \[
            A(s) \le A(t) \quad \text{when} \quad A(s) \supseteq A(t). 
        \]
    Moreover, for any set of objects $\{t_i\}$, we define the join of the set $\{ A(t_i) \}$ by
        \[
            \bigvee_i A(t_i) \coloneqq A(\coprod_i t_i). 
        \]
    The set $A(\cat T)$ of Bousfield classes is then a partially ordered set in which every subset has a greatest lower bound. By taking the join of all the lower bounds, we can define a meet operation~$\wedge$ on $A(\cat T)$. Under these operations $A(\cat T)$ forms a complete lattice called the \emph{Bousfield lattice} of $\cat T$. It was originally studied (for $\cat T$ the category of spectra) by Bousfield~\cite{Bousfield79_boolean}.
\end{Rem}

\begin{Rem}
    In general, it is difficult to describe the meet operation on the Bousfield lattice $A(\cat T)$. However, if $\cat T$ is stratified, then the following theorem gives a complete description of $A(\cat T)$ and hence also a description of the meet operation. For categories stratified in the sense of BIK (\cref{rem:BIK}), an analogous result was established in \cite[Section~4]{IyengarKrause13}.
\end{Rem}

\begin{Thm}\label{thm:bousfield-lattice}
    Let $\cat T$ be a stratified rigidly-compactly generated tt-category. The Bousfield lattice of $\cat T$ is isomorphic to the lattice of subsets of $\Spc(\cat T^c)$ via the map which sends a Bousfield class $A(t)$ to $\Supp(t)$. In particular,
        \[
            A(s) \wedge A(t) = A(s \otimes t)
        \]
    for any $s,t \in \cat T$.
\end{Thm}

\begin{proof}
    By the local-to-global principle, we have 
        \[
            \Loco{t} = \Loco{t \otimes \kappaP \mid \cat P \in \Supp(t)}.
        \]
    It follows that 
        \[
            A(t) = \bigvee_{\cat P \in \Supp(t)} A(\kappaP \otimes t).
        \]
    Arguing as in \cite[Theorem 4.4]{IyengarKrause13}, the tensor-product theorem for support (\cref{thm:tensor_product}), which holds whenever $\cat T$ is stratified, implies that $A(g(\cat P))$ is a minimal nonzero element of $A(\cat T)$. We deduce that $A(g(\cat P)) = A(\kappaP \otimes t)$ whenever $\cat P \in \Supp(t)$, so that 
        \[
            A(t) = \bigvee_{\cat P \in \Supp(t)} A(\kappaP).
        \]
    It follows that $\Supp(A(t)) = \Supp(t)^c$ for any $t \in \cat T$. This implies that $A(t) \mapsto \Supp(t)$ is a well-defined order-preserving bijection
        \[
            A(\cat T) = \big\{\text{Bousfield classes of } \cat T \big\} \xrightarrow{\simeq} \big\{ \text{subsets of } \Spc(\cat T^c) \big\}
        \]
    with order-preserving inverse $U \mapsto A(\coprod_{\cat P \in U}\kappaP)$.
\end{proof}

\begin{Rem}
    Note that under the bijection of \cref{thm:bousfield-lattice}, a Bousfield class $A\in A(\cat T)$ corresponds to the \emph{complement} of $\Supp(A) \subseteq \Spc(\cat T^c)$.
\end{Rem}

\begin{Rem}
    \Cref{thm:bousfield-lattice} implies that if $\cat T$ is stratified then there is no difference between localizing $\otimes$-ideals and Bousfield classes: Every localizing $\otimes$-ideal is a Bousfield class.
\end{Rem}

\section{The general telescope conjecture}\label{sec:gentelescope}

The goal of this section is to study the telescope conjecture for a rigidly-compactly generated tt-category $\cat T$ in the presence of stratification. We provide a point-set topological condition on the spectrum $\Spc(\cat T^c)$ which ensures that the telescope conjecture holds (in the presence of stratification), generalizing previous work of Stevenson \cite{Stevenson13,Stevenson14}.

\begin{Def}\label{def:telescope-conj}
    Let $\cat T$ be a rigidly-compactly generated tt-category. We say that~$\cat T$ satisfies the \emph{telescope conjecture} if every smashing $\otimes$-ideal of $\cat T$ is compactly generated. In other words, if every smashing localization of $\cat T$ is a finite localization. (Recall \cref{rem:compactly-generated}, \cref{rem:smashing_idempotents}, and \cref{rem:finite-localizations}.)
\end{Def}

\begin{Rem}\label{rem:support_smashing}
    If $\cat T$ is a \emph{stratified} rigidly-compactly generated tt-category, then the telescope conjecture for $\cat T$ is equivalent to the statement that the support $\Supp(\cat L)$ of any smashing $\otimes$-ideal $\cat L \subseteq \cat T$ is a Thomason subset of $\Spc(\cat T^c)$. Indeed, the support of a compactly generated localizing $\otimes$-ideal is Thomason, and the converse holds when we have stratification (cf.~\cref{thm:universal-support} applied to $\sigma = \Supp$).
\end{Rem}

\begin{Rem}
    The telescope conjecture was originally stated for the stable homotopy category $\SH$ by Ravenel \cite{Ravenel84} and is a notorious open problem. There are several equivalent formulations (see \cite{barthel_short_2020} and \cite{Wolcott15}, for example) but the one given by Bousfield \cite[Conjecture~3.4]{Bousfield79} can be stated for an arbitrary tensor-triangulated category~$\cat T$ as given in \cref{def:telescope-conj}. Despite the word ``conjecture,'' it is rather a property that the category $\cat T$ might or might not have, just as the category $\cat T$ might or might not have stratification. The ``telescope'' in the name is also no longer meaningful at this level of generality, and the ``telescope conjecture'' for $\cat T$ is sometimes referred to as the ``smashing conjecture'' for $\cat T$. However, this should not be confused with Ravenel's smashing conjecture introduced in the same paper \cite{Ravenel84} and mentioned in \cref{Rem:e_n}. In any case, we will return to the original $\cat T=\SH$ case of the telescope conjecture in \cref{sec:telescopeconjecture}.
\end{Rem}

\begin{Rem}\label{rem:telescope_conj_model}
    When $\cat T$ is stratified with noetherian spectrum and $\cat T$ arises as the homotopy category of a monoidal model category, the telescope conjecture has been shown to hold by Stevenson \cite[Theorem~7.15]{Stevenson13}. In fact, Stevenson considers more generally the case when $\cat T$ is acted upon by another tensor-triangulated category. We will give a generalization of Stevenson's result (in the case where $\cat T$ acts on itself) in \cref{thm:gennoethtelescopeconj} below. Tracing through Stevenson's proof, one sees that the assumption that $\cat T$ has a model is used (in Lemma 7.10 of \cite{Stevenson13}) to ensure that the local-to-global principle holds. This holds automatically in our situation by \cref{thm:local-to-global}. Therefore, for our purposes below, Stevenson's hypothesis that~$\cat T$ has an underlying model will be unnecessary.
\end{Rem}

\begin{Def}\label{defn:gennoetherian}
    A topological space $X$ is said to be \emph{generically noetherian} if the generalization closure $\gen(x)$ (\cref{rem:set-of-generalizations}) of any point $x \in X$ is noetherian.
\end{Def}

\begin{Exa}\label{ex:noethgennoeth}
    Every noetherian space is generically noetherian, because every subspace of a noetherian space is noetherian.
\end{Exa}

\begin{Exa}\label{ex:profinitegennoeth}
    Every $T_1$ spectral space (a.k.a.~profinite space) $X$ is generically noetherian, because $\gen(x)=\{x\}$ is a singleton for every point $x \in X$.
\end{Exa}

\begin{Rem}\label{rem:localgennoeth}
    It suffices to check the condition of \cref{defn:gennoetherian} for the closed points of $X$. In particular, it follows that a local space is generically noetherian if and only if it is noetherian. 
\end{Rem}

\begin{Lem}\label{lem:gennoethweaklynoeth}
    Let $X$ be a spectral space. If $X$ is generically noetherian, then it is weakly noetherian. 
\end{Lem}

\begin{proof}
    A point $x \in X$ is weakly visible if and only if $x \in \gen(x)$ is weakly visible, see \cref{rem:closed-weakly-visible}. By assumption, $\gen(x)$ is noetherian, so $x \in \gen(x)$ is visible and hence also weakly visible. 
\end{proof}

\begin{Exa}
    The converse to the statement of \cref{lem:gennoethweaklynoeth} does not hold. Indeed, there exists a spectral space $X = (\bbN^+)^{\triangleright}$ determined by the following two properties:
    \begin{enumerate}
        \item The underlying constructible space of $X$ is given by $\bbN^{+} \sqcup \{\omega\}$, i.e., the space obtained by adding a disjoint point $\omega$ to the one-point compactification of the discrete space of natural numbers.
        \item The specialization order $\rightsquigarrow$ on $X$ is given by $x \rightsquigarrow y$ for points $x,y \in X$ if and only if $y=\omega$.
    \end{enumerate}
    On the one hand, this space is not noetherian, as $\bbN^{+}$ is not noetherian. In light of \cref{rem:localgennoeth} and because $X$ is local with unique closed point $\omega$, it is not generically noetherian either. On the other hand, $X$ is weakly noetherian. Indeed, if $x \in X$ is different from $\omega$, we have $\{x\} = \gen(x)$ by construction, so $x \in X$ is weakly visible. If $x = \omega$, then $\omega$ is isolated in $\gen(\omega) = X$ for the constructible topology, hence weakly visible in $X$ by \cite[Corollary 4.5.19]{DickmannSchwartzTressl19}. We refer to \cite{BalchinBarthelGreenlees21ip} for further details on this construction.
\end{Exa}

\begin{Thm}\label{thm:gennoethtelescopeconj}
    Let $\cat T$ be a stratified rigidly-compactly generated $tt$-category with generically noetherian spectrum. The telescope conjecture holds for~$\cat T$. 
\end{Thm}

We start with a preparatory result, which might be of independent interest. 

\begin{Lem}[Stevenson]\label{lem:specclosedtelescopeconj}
    Suppose $\cat T$ is a stratified rigidly-compactly generated $tt$-category. The telescope conjecture holds for $\cat T$ if and only if the support $\Supp(\cat L)$ of any smashing $\otimes$-ideal $\cat L \subseteq \cat T$ is specialization closed.
\end{Lem}

\begin{proof}
    Let $X \coloneqq \Spc(\cat T^c)$. Bearing \cref{rem:visible} in mind, Stevenson proves in \cite[Proposition~2.14]{Stevenson14} that the right orthogonal~${\cat L}^{\bot}$ has proconstructible support; that is, $\Supp(\cat L^{\bot})$ is closed in the constructible topology. Note that if $e \to \unit \to f \to \Sigma e$ is the idempotent triangle associated to the smashing $\otimes$-ideal $\cat L$ then $\Supp(\cat L) = \Supp(e)$ and $\Supp(\cat L^{\bot}) = \Supp(f)$ (\cref{rem:smashing_idempotents}). Since $\cat T$ is stratified, the full tensor-product theorem (\cref{thm:tensor_product}) implies that $\Supp(\cat L^{\bot}) = \Supp(\cat L)^c$ are complements. Thus, if $\Supp(\cat L)$ is specialization closed then $\Supp(\cat L^{\bot})$ is generalization closed. By \cite[Theorem~1.4.3]{DickmannSchwartzTressl19}, this means $\Supp(\cat L^{\bot})$ is specialization closed in the Hochster dual $X^*$. Recall that $X$ and $X^*$ provide the same constructible topology. Hence \cite[Theorem~1.5.4(i)]{DickmannSchwartzTressl19} implies that $\Supp(\cat L^{\bot})$  is a closed subset of $X^*$. Thus the complement $\Supp(\cat L)$ is an open subset of $X^*$, that is, a Thomason subset of $X$. Since $\cat T$ is stratified, the localizing $\otimes$-ideal $\cat L$ must therefore be compactly generated, as desired (\cref{rem:support_smashing}). 
\end{proof}

\begin{proof}[Proof of \cref{thm:gennoethtelescopeconj}]
	Consider a smashing $\otimes$-ideal $\cat L \subseteq \cat T$ with
	associated idempotent triangle $e \to \unit \to f \to \Sigma e$
	(\cref{rem:smashing_idempotents}). By \cref{lem:specclosedtelescopeconj},
	it suffices to prove that $\Supp(\cat L)=\Supp(e)$ is specialization closed
	in $\Spc(\cat T^c)$. In other words, we need to show that for any $\cat Q
	\in \Supp(\cat L)$ and any $\cat P \in \overline{\{\cat Q\}}$ we have $\cat
	P \in \Supp(\cat L)$. To this end, consider the localization $L\colon \cat
	T \to \cat T/\langle \cat P \rangle$ to the local category at~$\cat P$.
	Applying $L$ to the idempotent triangle, we obtain an induced smashing
	localization on $\cat T/\langle\cat P\rangle$ with corresponding smashing
	$\otimes$-ideal $\cat L' = \ker(Lf \otimes -)=\Loco{Le} \subseteq \cat
	T/\langle \cat P \rangle$. Identifying $\Spc((\cat T/\langle\cat
	P\rangle)^c)$ with $\gen(\cat P) \subseteq \Spc(\cat T^c)$ (as in
	\cref{exa:local-categories}), it follows from \cref{prop:ltg-permanence}(c)
	(applied to $d\coloneqq Le$)
	and \cref{lem:support-compat-with-finite-loc} that \[
            \Supp(\cat L') = \Supp(\cat L) \cap \gen(\cat P). 
        \]
    In particular, $\cat Q \in \Supp(\cat L')$. Since $\Spc(\cat T^c)$ is generically noetherian, $\gen(\cat P) = \Spc((\cat T/\langle\cat P\rangle)^{c})$ is noetherian. Moreover, since $\cat T$ is stratified, the localization $\cat T/\langle \cat P \rangle$ is also stratified (\cref{cor:passes-to-finite}). Thus, by the telescope conjecture for stratified rigidly-compactly generated tt-categories having noetherian spectrum verified by Stevenson \cite[Theorem 7.15]{Stevenson13} (recalling \cref{rem:telescope_conj_model}), we conclude that $\Supp(\cat L')$ is specialization closed as a subset of $\gen(\cat P)$. Since $\cat Q \in \Supp(\cat L')$, it follows that $\cat P$ is contained in $\Supp(\cat L')$ and hence also in $\Supp(\cat L)$. We conclude that $\Supp(\cat L)$ is specialization closed, as desired.
\end{proof} 

\begin{Exa}\label{ex:rational-gspectra-gennoeth}
    Let $G$ be a compact Lie group and write $\SH_{G,\bbQ}$ for the stable homotopy category of rational $G$-spectra. It follows from the main theorem of \cite{Greenlees19_rational} that the spectrum $\Spc(\SH_{G,\bbQ}^c)$ is generically noetherian. The details of this argument will be part of forthcoming work \cite{BalchinBarthelGreenlees21ip}. Moreover, the category $\SH_{G,\bbQ}$ is stratified (see \cref{thm:rational_stratification} below) and hence satisfies the telescope conjecture by \cref{thm:gennoethtelescopeconj}. Note that, in general, $\Spc(\SH_{G,\bbQ}^c)$ is neither noetherian nor constructible (see also \cref{rem:o2}), so this result is not covered by the existing results on the telescope conjecture in the literature.
\end{Exa}

\begin{Rem}
    In general, the assumption in \cref{thm:gennoethtelescopeconj} that $\cat T$ is stratified cannot be dropped. Indeed, Keller's example \cite{Keller94b} of a derived category of a commutative ring for which the telescope conjecture fails is rigidly-compactly generated and has spectrum consisting of only two points. See also \cite[Example 5.7]{bhs2}.
\end{Rem}

\part{Applications: Chromatic homotopy theory}\label{part:chromatic}

We now turn to applications and examples of stratification in chromatic homotopy theory. We prove that the category of $E(n)$-local spectra is stratified, establish further connections between stratification and the telescope conjecture, and prove that the category of rational $G$-spectra is stratified for all compact Lie groups $G$.

\section{Stratification of the \texorpdfstring{$E(n)$}{E(n)}-local category}\label{sec:En-local}

\begin{Not}
    Fix a prime number $p$ and let $\cat S \coloneqq \SH_{(p)}$ denote the $p$-local stable homotopy category. For $0 \le n \le \infty$, let $K(n)$ denote the $p$-local Morava $K$-theory and, for $n < \infty$, let $E(n)$ denote the $p$-local Johnson--Wilson spectrum. We assume the reader has some familiarity with the basic properties of these spectra, as discussed in \cite{HoveyStrickland99}, for example. Ravenel's book \cite{Ravenel92} is also instructive.
\end{Not}

\begin{Rem}
    We follow the usual conventions in homotopy theory where $K(0)$ denotes rational homology and $K(\infty)$ denotes mod $p$ homology.  
\end{Rem}

\begin{Rem}\label{Rem:e_n}
    Let $L_n \colon  \cat S \to \cat S$ denote Bousfield localization with respect to $E(n)$. This is a smashing localization, as conjectured by Ravenel \cite[Section 10.6]{Ravenel84} and proved by Hopkins--Ravenel \cite[Theorem~7.5.6]{Ravenel92}. It follows that the $E(n)$-local stable homotopy category $\cat S_{E(n)}$ is rigidly-compactly generated (\cref{ter:tt-category}). In fact, it is generated by its compact monoidal unit $\unit=L_n S^0$ (cf.~\cref{rem:suffice-tensor-compact}).
\end{Rem}

\begin{Rem}
    The localizing $\otimes$-ideals of $\cat S_{E(n)}$ have been classified by Hovey--Strickland \cite{HoveyStrickland99}. Our present goal is to formulate this result in  terms of stratification and the Balmer spectrum. See \cref{thm:hovey-strickland} and \cref{thm:En-stratified} below. These results will serve as the nonequivariant base case for \cref{cor:En-local-mack-strat} in \cref{part:spectral-mackey} which stratifies $E(n$)-local spectral Mackey functors. We start by reviewing the Balmer spectrum of the stable homotopy category of compact $p$-local spectra~$\cat S^c$.
\end{Rem}

\begin{Not}\label{not:Cn}
    For each $0 \le n \le \infty$, let $\cat C_n \coloneqq \SET{ x \in \cat S^c }{K(n)_*(x) = 0}$. Beware that some authors use a slightly different choice of indexing.
\end{Not}

\begin{Rem}\label{rem:definitioncn}
    Since the coefficient ring $K(n)_*$ is a graded field, Morava $K$-theory satisfies a K\"{u}nneth formula and $\cat C_n$ is a \emph{prime} thick $\otimes$-ideal of~$\cat S^c$. In other words,~$\cat C_n$ is a point in the Balmer spectrum $\Spc(\cat S^c)$.  The Thick Subcategory Theorem of Hopkins--Smith \cite{HopkinsSmith98}, as interpreted by Balmer \cite[Corollary~9.5]{Balmer10b}, provides a complete description of this space:
\end{Rem}

\begin{Thm}[Hopkins--Smith]\label{thm:hopkins_smith_thick}
    The spectrum
        \[
            \Spc(\cat S^c) = \cat C_\infty - \dots - \cat C_{n+1} - \cat C_{n} - \cdots - \cat C_{1} - \cat C_0
        \]
	is an infinite tower of connected points, where closure goes to the left:
$\overbar{\{\cat C_n\}} = \SET{\cat C_k}{n \le k \le \infty}$.  In particular,
the space is local and irreducible with generic point $\cat C_0 =
\SH_{(p),\mathrm{tor}}^c$ and unique closed point $\cat C_\infty = (0)$. \end{Thm}

\begin{Rem}\label{rem:type}
    A compact $p$-local spectrum $x \in \cat S^c$ has \emph{type} $n$ if $K(n)_*(x) \neq 0$ but $K(m)_*(x) = 0$ for $m < n$. Every nonzero compact $p$-local spectrum has a well-defined type. Indeed $x$ has type $n$ if and only if $x \in \cat C_{n-1} \setminus \cat C_n$ if and only if $\supp(x) = \overbar{\{\cat C_n\}}$ if and only if the codimension of $\supp(x)$ is $n$. (Here we interpret $\cat C_{-1} \coloneqq \cat S^c$ to be the whole category of compact $p$-local spectra.)
\end{Rem}

\begin{Not}\label{not:En-primes}
    For each $0 \le h \le n$, we define
        \[ 
            \cat P_h \coloneqq \SET{x \in \cat S_{E(n)}^c }{K(h)_*(x)=0},
        \]
    a prime $\otimes$-ideal of $\cat S_{E(n)}^c$.  Equivalently, by \cite[Proposition~6.8]{HoveyStrickland99},
        \[ 
            \cat P_h = \SET{x \in \cat S_{E(n)}^c}{K(i)_*(x) =0 \text{ for all } 0 \le i \le h}.
        \]
\end{Not}

\begin{Thm}[Hovey--Strickland] \label{thm:hovey-strickland}
    The spectrum
        \[
            \Spc(\cat S_{E(n)}^c) = \cat P_{n} - \cdots - \cat P_{1} - \cat P_0
        \]
    is a local irreducible space consisting of $n+1$ points, where closure goes to the left: $\overbar{\{\cat P_h\}} = \SET{\cat P_k}{h \le k \le n}$.
\end{Thm}

\begin{proof}
    This was established in \cite[Theorem~6.9]{HoveyStrickland99} and reinterpreted in tt-geometric language in \cite[Proposition~3.5]{BarthelHeardNaumann22}.
\end{proof}

\begin{Rem}\label{rem:compute-gP}
    Our goal is to establish stratification for the $E(n)$-local category. First observe that (in the notation of \cref{exa:local-categories} and \cref{rem:set-of-generalizations})
        \[
            Y_{\cat P_k}  = \SET{\cat P_i \in \Spc(\cat S_{E(n)}^{c})}{\cat P_k \not \subseteq \cat P_i } = \big\{ \cat P_{k+1},\cat P_{k+2},\ldots,\cat P_{n}\big\}
        \]
    for any $0 \le k \le n$. Moreover, finite localization on $\cat S_{E(n)}$ away from $L_nF(k+1)$ is equivalent to $L_k$. This follows from \cite[Corollary~6.10]{HoveyStrickland99}, which only uses the thick subcategory theorem. (Here $F(k+1) \in \cat S^c$ denotes a type $k+1$ spectrum; recall \cref{rem:type}.) This implies that $f_{Y_{\cat P_k}} \simeq L_kS^0$. It follows that $e_{Y_{\cat P_k}} \simeq M_{k,n}S^0$ where $M_{k,n}S^0$ denotes the fiber of the natural map $L_nS^0 \to L_kS^0$. For $k \ge 1$, we have $\overbar{\{\cat P_{k}\}} = Y_{\cat P_{k-1}}$ and hence $e_{\overbar{\{\cat P_k\}}} \simeq M_{k-1,n}S^0$. This is also true for $k=0$ with the interpretation that $L_{-1}S^0 = 0$ and $M_{-1,n}S^0 = L_n S^0$.

    Recall that for each $\cat P_k \in \Spc(\cat S_{E(n)}^c)$, we have the Balmer--Favi idempotent
        \[ 
            g(\cat P_k) \coloneqq e_{\overbar{\{\cat P_k\}}}\otimes f_{Y_{\cat P_k}}.
        \]
    The calculations above show that $g(\cat P_k) \simeq M_{k-1,n}S^0 \otimes L_kS^0 \simeq L_kM_{k-1,n}S^0$ where the last equivalence uses the fact that both functors are smashing. The fiber sequence
        \[
            L_kM_{k-1,n}S^0 \to L_kL_nS^0 \to L_kL_{k-1}S^0
        \] 
    along with the fact that $L_iL_j \simeq L_{\min(i,j)}$ shows that $L_kM_{k-1,n}S^0 \simeq M_{k-1,k}S^0$ which is usually denoted in the literature as $M_kS^0$. To summarize, we have established the following:
\end{Rem}

\begin{Prop}\label{prop:kappa-of-En}
    For each $0 \le k \le n$, the Balmer--Favi idempotent $g(\cat P_k)$ is given by $g(\cat P_k) \simeq M_kS^0$. The support of $X \in \cat S_{E(n)}$ is 
        \[
            \begin{split}
				\Supp(X) &= \SET{\cat P_k \in \Spc(\cat S_{E(n)}^{c})}{g(\cat P_k) \otimes X \ne 0} \\
				&= \SET{\makebox[\widthof{$\cat P_k \in \Spc(\cat S_{E(n)}^{c})$}][c]{$k \in \{0,\ldots, n\}$}}{\makebox[\widthof{$g(\cat P_k) \otimes X \ne 0$}][c]{$M_kS^0 \otimes X \ne 0$}}.
            \end{split}
        \]
\end{Prop}

\begin{Rem}\label{rem:chromatic_support}
    By \cite[Proposition~5.3]{HoveyStrickland99}, $K(k) \otimes X \neq 0$ if and only if ${M_kS^0 \otimes X \neq 0}$. The Balmer--Favi notion of support thus coincides with the usual notion of chromatic support defined using Morava $K$-theories.
\end{Rem}

\begin{Thm}[Hovey--Strickland]\label{thm:En-stratified}
    The $E(n)$-local stable homotopy category~$\cat S_{E(n)}$ is stratified.
\end{Thm}

\begin{proof}
    The space $\Spc(\cat S^c_{E(n)})$ was described in \cref{thm:hovey-strickland}. It is noetherian since it has finitely many points. It therefore satisfies the local-to-global principle by \cref{thm:local-to-global}. In fact, it is easy to give a direct verification: For any $X \in \cat S_{E(n)}$ we have a tower of cofiber sequences
        \[
            \begin{tikzcd}[column sep=small]
            M_nS^0 \otimes X \arrow[d] & M_{n-1}S^0 \otimes X \arrow[d] &                  & M_1S^0 \otimes X \arrow[d] & M_0S^0 \otimes X \arrow[d] &   \\
            X \arrow[r]    & L_{n-1}X \arrow[r] & \cdots \arrow[r] & L_1X \arrow[r] & L_0X \arrow[r] & 0
            \end{tikzcd}
        \]
    and hence $X \in \Loco{M_kS^0 \otimes X \mid k\in \Supp(X)}$, so that 
        \begin{equation}\label{eq:ltg-En}
            \Loco{X} \subseteq \Loco{M_kS^0 \otimes X \mid k \in \Supp(X)}.
        \end{equation}
    The converse holds automatically, so \eqref{eq:ltg-En} is an equality. We note in passing that a closely related result is obtained in \cite[Proposition~6.18]{HoveyStrickland99}.

    To establish stratification it suffices by \cref{thm:equiv-strat} and \cref{prop:kappa-of-En} to show that $\Loco{g(\cat P_k)} = \Loco{M_kS^0}$ is a minimal localizing subcategory of~$\cat S_{E(n)}$. Note that this is the monochromatic category $\cat M_k$ studied in \cite[Section~6.3]{HoveyStrickland99}. Using the equivalence of categories $\cat M_k \cong \cat S_{K(k)}$ between the monochromatic category and the $K(k)$-local category \cite[Theorem~6.19]{HoveyStrickland99} the claim follows from \cite[Theorem~7.5]{HoveyStrickland99}.
\end{proof}

\begin{Rem}
    Because $\Spc(\cat S_{E(n)}^c)$ is noetherian, \cref{thm:En-stratified} and \cref{thm:gennoethtelescopeconj} imply that the telescope conjecture holds in $\cat S_{E(n)}$. This was previously shown by Hovey and Strickland \cite[Corollary~6.10]{HoveyStrickland99} using only the thick subcategory theorem (equivalently, the computation of $\Spc(\cat S_{E(n)}^c))$. In fact, we have already used this result to determine $g(\cat P_k)$ in \cref{rem:compute-gP}.
\end{Rem}

\begin{Rem}
    The stratification of the category of $E(n)$-local spectra provided by \cref{thm:En-stratified} is of particular interest because it is an example of a tt-category which is not canonically BIK-stratified (recall \cref{rem:BIK}).
\end{Rem}

\section{The chromatic telescope conjecture}\label{sec:telescopeconjecture}

In this section, we'll reformulate the classical telescope conjecture in stable homotopy theory in terms of stratification.

\begin{Not}\label{not:enfn}
    For each $0 \le n <\infty$, let
        \[ 
            e_n \to \unit \to f_n \to \Sigma e_n 
        \]
    denote the idempotent triangle in $\cat S\coloneqq \SHp$ corresponding to the thick $\otimes$-ideal $\cat C_n$ (\cref{rem:finite-localizations} and \cref{rem:definitioncn}), so that $\Ker(f_n \otimes -) = \Loco{e_n} = \Loc\langle \cat C_n \rangle$.  Thus $f_n \otimes -\colon \cat S \to \cat S$ is the finite localization whose subcategory of compact acyclics is precisely $\cat C_n$.  It corresponds to the Thomason subset
        \[ 
            Y_n \coloneqq \bigcup_{x \in \cat C_n} \supp(x) = \overbar{\{\cat C_{n+1}\}}
        \]
    whose complement is the open subset
        \[ 
            U_n \coloneqq \Spc(\cat S^c)\setminus Y_n = \{ \cat C_n, \cat C_{n-1}, \cdots, \cat C_0\}.
        \]
    We could also consider the $n=\infty$ case, but this is just the trivial localization ($e_\infty = 0$) corresponding to the empty Thomason subset.
\end{Not}

\begin{Rem}\label{rem:spec-of-fn}
    As explained in \cref{rem:finite-localizations}, the finite localization 
        \[
            \cat S \to \cat S_n \coloneqq \cat S(U_n) \coloneqq \cat S/\Loco{e_n}=\cat S/\Loc\langle \cat C_n \rangle\cong f_n\otimes \cat S
        \]
    excises the closed set $\overbar{\{\cat C_{n+1}\}}$:
        \begin{equation}\label{eq:open-piece}
            \Spc(\cat S_n^c) \cong U_n = \{ \cat C_n, \cat C_{n-1}, \ldots , \cat C_0 \} \subset \Spc(\cat S^c).
        \end{equation}
\end{Rem}

\begin{Rem}
    A result due to Ravenel \cite[Theorem~2.1]{Ravenel84} provides a decomposition of the Bousfield class of $E(n)$ in terms of the Morava $K$-theories: 
        \[
            \langle E(n)\rangle = \langle K(0) \rangle \vee \langle K(1) \rangle \vee \cdots \vee \langle K(n) \rangle.
        \]
    Moreover, since $L_n$ is smashing (recall \cref{Rem:e_n}), we also have $\langle E(n) \rangle = \langle L_n S^0 \rangle$. These results of Ravenel and Hopkins--Ravenel establish that a $p$-local spectrum $t \in \cat S$ is $E(n)$-acyclic if and only if it is $K(i)$-acyclic for each $0 \le i \le n$. However, if a \emph{compact} $p$-local spectrum is $K(i+1)$-acyclic then it is necessarily $K(i)$-acyclic. This amounts to the inclusion $\cat C_{i+1} \subseteq \cat C_i$ (cf.~\cref{thm:hopkins_smith_thick}) and was proved in \cite[Theorem~2.11]{Ravenel84}. Thus
        \[
            \SET{x \in \cat S^c}{E(n)_*(x) = 0} = \SET{x \in \cat S^c}{K(n)_*(x) = 0} = \cat C_n.
        \]
    In other words, $L_n$ has the same finite acyclics as $L_{K(n)}$, namely $\cat C_n$. Thus the finite localization $L_n^f$ associated to $L_n$ is the localization $f_n \otimes -$ from \cref{not:enfn}. We use the notation $\cat S_n$ for the finitely localized category (as has already appeared in \cref{rem:spec-of-fn}). We can now state Ravenel's Telescope Conjecture:
\end{Rem}

\begin{Conj}[The Telescope Conjecture]
    Let $p$ be a prime number. For any integer $n \ge 0$, $E(n)$-localization on the $p$-local stable homotopy category coincides with its associated finite localization: $L_n^f \xra{\sim} L_n$.
\end{Conj}

\begin{Rem}
    As explained in \cite[Section 4]{barthel_short_2020}, this is equivalent to the telescope conjecture of \cref{def:telescope-conj} for the $p$-local stable homotopy category $\cat S = \SHp$. In other words, the statement that every smashing localization on $\SHp$ is finite reduces to the case of the particular smashing localizations~$L_n$, $n \ge 0$.
\end{Rem}

\begin{Rem}\label{rem:induced-homeo}
    For each $0 \le n <\infty$, the smashing localization $\cat S \to \cat S_{E(n)}$ induces a smashing localization 
        \begin{equation}\label{eq:smashing-comparison}
            \cat S_n \to \cat S_{E(n)}
        \end{equation}
    which has no finite acyclics, and which the Telescope Conjecture asserts is an equivalence. From \cref{thm:hovey-strickland} and \cref{rem:spec-of-fn}, we see that the smashing localization~\eqref{eq:smashing-comparison} induces a homeomorphism 
        \[
            \Spc(\cat S_{E(n)}^c) \xra{\cong} \Spc(\cat S_n^c)
        \]
    between Balmer spectra.
\end{Rem}

\begin{Prop}\label{prop:homeo-smashing}
    Let $\cat T$ be a stratified rigidly-compactly generated tt-category. If a smashing localization $\cat T \to \cat L$ induces a homeomorphism
        \[
            \varphi\colon \Spc(\cat L^c) \xra{\cong} \Spc(\cat T^c)
        \]
    then the localization must be trivial.
\end{Prop}

\begin{proof}
    A smashing localization $\cat L$ of a rigidly-compactly generated category $\cat T$ is again rigidly-compactly generated. Thus we can apply \cref{prop:ltg-permanence} to the localization functor $L\colon \cat T \to \cat L$. We claim that $\Supp(\Ker(L)) \cap \im(\varphi) =\emptyset$. Indeed, if $\varphi(\cat Q) \in \Supp(t)$ then $g(\varphi(\cat Q)) \in \Loco{t}$ (\cref{rem:min-consequence}). Hence $g(\cat Q) \cong L(g(\varphi(\cat Q))) \in \Loco{L(t)}$ which is a contradiction if $L(t)=0$. In particular, the left idempotent~$e$ of the smashing localization (\cref{rem:smashing_idempotents}) is in the kernel. Since $\varphi$ is surjective by assumption, we conclude that $\Supp(e)=\emptyset$ so that $e=0$.
\end{proof}

\begin{Cor}\label{cor:telescope-strat}
    Let $n \ge 0$ be an integer and $p$ a prime number. The Telescope Conjecture holds at height~$n$ and prime $p$ if and only if the finite localization $\cat S_n\coloneqq \cat S/\langle \cat C_n\rangle$ of the $p$-local stable homotopy category $\cat S\coloneqq \SHp$ is stratified.
\end{Cor}

\begin{proof}
    Recall from \cref{rem:induced-homeo} that the induced smashing localization $\cat S_n \to \cat S_{E(n)}$ induces a homeomorphism $\Spc(\cat S_{E(n)}^c) \cong \Spc(\cat S_n^c)$. Thus, if $\cat S_n$ is stratified then \cref{prop:homeo-smashing} implies that this induced smashing localization is trivial. In other words, $E(n)$-localization coincides with its associated finite localization, which is the Telescope Conjecture at height~$n$. Alternatively, \cref{thm:gennoethtelescopeconj} implies that if $\cat S_n$ is stratified then every smashing localization on $\cat S_n$ is finite. The induced smashing localization $\cat S_n \to \cat S_{E(n)}$ has no nonzero compact acyclics by construction hence must be trivial. Conversely, if the Telescope Conjecture holds at height~$n$ then $\cat S_{n} \to \cat S_{E(n)}$ is an equivalence of tensor-triangulated categories, and \cref{thm:En-stratified} establishes that $\cat S_{E(n)}$ is stratified.
\end{proof}

\begin{Rem}
    The spectrum of $\cat S\coloneqq \SH_{(p)}$ is unfortunately not weakly noetherian. The problem is that the closed point $\cat C_\infty$ is not ``weakly visible'' in the sense of \cref{rem:visible}. All the other primes are visible, however, and one can still consider the Balmer--Favi support 
        \[
            \Supp_{<\infty}(t) \coloneqq \SET{ \cat C_i}{ 0\le i<\infty, g(\cat C_i) \otimes t \neq 0}.
        \]
	As shown in \cite[Remark 5.13]{bhs2}, we have 
		\[
            \Supp_{<\infty}(t) = \SET{ i \in \mathbb{N}}{ T(i) \otimes t \neq 0}
        \]
    where $T(i)$ is the telescope of a finite type $i$ spectrum. However, this notion of support doesn't satisfy the local-to-global principle nor even the detection property. For example, dissonant spectra will have empty support. In fact, even if we extend this notion of support to the closed point $\cat C_\infty$ by defining $\cat C_\infty \in \Supp(t)$ iff $\HFp \otimes t \neq 0$, there are still nonzero objects in $\SH_{(p)}$ with empty support. Indeed, an example is given by the $p$-local Brown--Comenetz dual of the sphere spectrum (see \cite[Corollary~B.12]{HoveyStrickland99} and \cite[Lemma 7.1(d)]{HoveyPalmieri99}). Nevertheless, $\Supp_{<\infty}$ allows us to make the following statement:
\end{Rem}

\begin{Cor}\label{cor:support-approach-to-tc} 
    Let $p$ be a prime number. The telescope conjecture for $\SH_{(p)}$ is equivalent to the statement that the surjective map
        \[
             \left\{
             \begin{array}{c}
             \text{localizing $\otimes$-ideals of $\SH_{(p)}$ which}\\
             \text{contain a nonzero compact object} 
             \end{array}
             \right\} 
             \to
             \left\{
             \begin{array}{c}
             \text{subsets of $\Spc(\SH_{(p)}^c)$ which}\\
             \text{contain a nonempty closed subset}
             \end{array}
             \right\}
        \]
    defined by $\cat L \mapsto \big\{\cat C_{\infty}\big\} \cup \Supp_{<\infty}(\cat L)$ is injective (and hence a bijection).
\end{Cor}

\begin{proof}
    For any $t \in \cat S \coloneqq \SHp$ and Thomason subset $Y \subseteq \Spc(\cat S^c)$, we have $\Supp_{<\infty}(t \otimes e_Y) = \Supp_{<\infty}(t) \cap Y$ as in the proof of \cref{lem:support-compat-with-finite-loc}. In particular, $\Supp_{<\infty}(e_Y) = Y \setminus \{\cat C_\infty\}$ and, as in the proof of \cref{lem:half-tensor}, $\Supp_{<\infty}(x) = \supp(x) \setminus \{ \cat C_\infty \}$ for any compact $x \in \cat S^c$. Define $\Supp(t) \coloneqq \Supp_{<\infty}(t) \cup \{ \cat C_{\infty}\}$. Note that a subset of $\Spc(\cat S^c)$ contains a nonempty closed subset if and only if it contains the closed point $\cat C_\infty$. The map in the statement of the Corollary is thus surjective by (the proof of) \cref{lem:surjectivity}. Now, for any $0 \le n <\infty$, let $L_n^f\colon\cat S \to \cat S_n$ be the corresponding finite localization of \cref{rem:spec-of-fn} and let $\varphi_n\colon \Spc(\cat S_n^c) \hookrightarrow \Spc(\cat S^c)$ be the induced map on spectra. We claim that the following diagram commutes:
        \begin{equation}\label{eq:mini-telescope}
            \begin{tikzcd}
             \bigg\{
             \begin{array}{c}
             \makebox[\widthof{\text{localizing $\otimes$-ideals of $\cat S_n$}}][l]{\text{localizing $\otimes$-ideals of $\cat S$}}\\
             \text{which contain $e_n$} 
             \end{array}
             \bigg\} 
             \ar[r,"\Supp","(\dagger 1)"'] &
             \bigg\{
             \begin{array}{c}
             \text{subsets of $\Spc(\cat S^c)$}\\
             \text{which contain $\overbar{\{\cat C_{n+1}\}}$}
             \end{array}
             \bigg\} \ar[d,"\varphi_n^{-1}","\cong"']\\
             \bigg\{
             \begin{array}{c}
             \text{localizing $\otimes$-ideals of $\cat S_n$}
             \end{array}
             \bigg\} \ar[u,"(L_n^f)^{-1}","\cong"']
             \ar[r,"\Supp","(\dagger 2)"'] &
             \bigg\{
             \begin{array}{c}
             \makebox[\widthof{\text{which contain $\overbar{\{\cat C_{n+1}\}}$}}][c]{\text{subsets of $\Spc(\cat S_n^c)$}}
             \end{array}
             \bigg\}.
             \end{tikzcd}
        \end{equation}
    Although $\cat S$ does not have weakly noetherian spectrum, the proof of \cref{prop:ltg-permanence}(c) still establishes that $\Supp_{<\infty}(U_n(t)) = \varphi_n(\Supp(t))$ for any $t \in \cat S_n$. Here $U_n\colon\cat S_n \to \cat S$ is the right adjoint of $L_n^f\colon\cat S \to \cat S_n$. Since $\varphi_n$ is injective and $\cat C_{\infty} \not\in \im\varphi_n$, it follows that $\varphi_n^{-1}(\Supp(U_n(t))) = \Supp(t)$ for all $t\in \cat S_n$. Moreover, note that $\varphi_n^{-1}(\Supp(e_n)) = \varphi_n^{-1}(\overbar{\{\cat C_{n+1}\}})=\emptyset$. The commutativity of \eqref{eq:mini-telescope} then follows from the definitions, \cref{lem:pullback-localizing}, and \cref{rem:spec-of-fn}. The vertical maps are bijections. Hence the injectivity of $(\dagger 1)$ is equivalent to the injectivity of $(\dagger 2)$, which by \cref{thm:equiv-strat} and \cref{cor:telescope-strat} is equivalent to the Telescope Conjecture at height $n$. The injectivity of the map in the corollary	evidently implies the injectivity of $(\dagger 1)$ for each $n$ and hence implies the Telescope Conjecture. Conversely, suppose $\cat L_1$ and $\cat L_2$ are localizing $\otimes$-ideals of $\cat S$ which contain nonzero compact objects, say $x_1$ and $x_2$, respectively. Then $x_1 \otimes x_2$ is a nonzero compact object (since $\cat S$ is local) which is common to both localizing $\otimes$-ideals, and we have $\supp(x_1 \otimes x_2) = \overbar{\{\cat C_{n+1}\}}$ for some $0 \le n <\infty$. Hence both $\cat L_1$ and $\cat L_2$ contain~$e_n$. Thus, if $\Supp(\cat L_1) = \Supp(\cat L_2)$, then the injectivity of $(\dagger 1)$ for this particular $n$ (which we have established is equivalent to the Telescope Conjecture at height $n$) implies that $\cat L_1 = \cat L_2$.
\end{proof}

\begin{Rem}
    There is some evidence that the Telescope Conjecture at height~$n$ may be false (see \cite{ravenel92a, MahowaldRavenelShick01}) for $n>1$. Nevertheless, \cref{cor:support-approach-to-tc} clarifies that in order to \mbox{disprove} the Telescope Conjecture (say at height $n=2$), it would suffice to find any two distinct localizing $\otimes$-ideals which contain nonzero compact objects, and demonstrate that they have the same Balmer--Favi supports.
\end{Rem}

\section{Rational equivariant spectra}\label{sec:rational-G-spectra}

Let $G$ be a compact Lie group and let $\SH_{G,\bbQ}$ denote the category of rational equivariant $G$-spectra. In this section, we review the results of Greenlees \cite{Greenlees19_rational} and show that his work can be interpreted as the statement that $\SH_{G,\bbQ}$ is stratified. As we shall see, this is an example of stratification where the Balmer spectrum $\Spc(\SH_{G,\bbQ}^c)$ is weakly noetherian, but not noetherian in general.

\begin{Rem}
    The category $\SH_{G,\bbQ}$ is a rigidly-compactly generated tt-category, with a set of compact-rigid generators $\{ \Sigma^{\infty}_G G/H_+ \otimes \bbQ \} $ with $H$ running through the set of closed subgroups of $G$. For legibility, we will omit the rationalization and $G$-suspension spectrum from our notation.
\end{Rem}

\begin{Rem}
    The Balmer spectrum of $\SH_{G,\bbQ}^c$ is determined, for any compact Lie group $G$, by Greenlees \cite[Theorem~1.3]{Greenlees19_rational}. In fact, he first determines the localizing $\otimes$-ideals of $\SH_{G,\bbQ}$ and uses this to deduce the classification of thick $\otimes$-ideals of~$\SH_{G,\bbQ}^c$. This is the opposite of how we will proceed. We also note that an independent calculation of the spectrum of $\SH_{G,\bbQ}^c$ can be deduced from the description of the Balmer spectrum of the whole $G$-equivariant stable homotopy category $\SH_G^c$ obtained by Barthel, Greenlees, and Hausmann; see \cite[Remark~1.5]{bgh_balmer}.
\end{Rem}

\begin{Rem}\label{rem:geometric-isotropy}
    We recall that there exist tt-functors $\Phi^H \colon \SH_{G,\bbQ} \to \SH_{\mathbb{Q}}$, one for each closed subgroup $H \le G$, and that these ``geometric fixed point'' functors are jointly conservative (see \cite[Proposition~3.3.10]{Schwede18_global}). The \emph{geometric isotropy} of $X \in \SH_{G,\bbQ}$ is defined as 
        \[
            \SET{ H \le G }{ \PhiHX \neq 0 }.
        \]
\end{Rem}

\begin{Def}
    For each closed subgroup $H \le G$, let $\mathfrak{p}_H$ denote the prime ideal
        \[
            \mathfrak{p}_H \coloneqq \SET{ x \in \SH_{G,\bbQ}^c }{ \PhiHx = 0 } \in \Spc(\SH_{G,\bbQ}^c).
        \]
    Note that it only depends on the $G$-conjugacy class of $H$.
\end{Def}

\begin{Rem}
    For the following result, recall that, by definition, $L$ is \emph{cotoral} in $K$ if $L$ is a normal subgroup of $K$ and $K/L$ is a torus. 
\end{Rem}

\begin{Thm}[Greenlees]\label{thm:spectrumrationalgspectra}
    Let $G$ be a compact Lie group. Then 
        \[
            \Spc(\SH_{G,\bbQ}^c) = \SET{ \mathfrak{p}_H}{H \le G \text{ is a closed subgroup}}.
        \]
    The specialization order is determined by cotoral inclusions:
        \[
            \frak{p}_K \subseteq \mathfrak{p}_H  \text{ if and only if $K$ is conjugate to a subgroup cotoral in $H$.}
        \]
    The topology on $\Spc(\SH_{G,\bbQ}^c)$ is the ``Zariski topology on the $f$-topology'' of \cite{Greenlees98_rational}. Under this identification, the support of an object $x \in \SH_{G,\bbQ}^c$ corresponds to its geometric isotropy:
        \[
            \supp(x) = \SET{\mathfrak{p}_H}{\PhiHx \neq 0}.
        \]
\end{Thm}

\begin{Rem}\label{rem:o2}
    If $G$ is finite, then $\Spc(\SH_{G,\bbQ}^c)$ is a finite space and is thus noetherian, but this does not hold in general for a compact Lie group. For example, it follows from \cite{Greenlees98_O2, Barnes17_O2} that the Balmer spectrum of the dihedral part of $\SH_{O(2),\bbQ}$ is homeomorphic to the one-point compactification of the set of natural numbers. In particular, $\Spc(\SH_{O(2),\bbQ}^c)$ is not noetherian. However, $\Spc(\SH_{G,\bbQ}^c)$ is always weakly noetherian, as we show in \cref{lem:rational_weakly_noetherian}.
\end{Rem}

\begin{Rem} \label{rem:families}
    Let $\calF$ be a family of closed subgroups of $G$, that is, a collection of closed subgroups which is closed under subconjugation~$\le_G$. As in \cite[Example~5.14]{BalmerSanders17}, we have a thick $\otimes$-ideal of  $\SH_{G,\bbQ}^c$ defined by $\cat {C}_{\calF} \coloneqq \thickt{G/H_+ \mid H \in \calF}$ with associated Thomason subset $ Y_{\calF} \coloneqq \bigcup_{H \in \calF}\supp(G/H_+)$. It follows from \cref{thm:spectrumrationalgspectra} that
        \begin{equation}\label{eq:YF}
            Y_{\calF} = \SET{\mathfrak{p}_H}{H \in \calF}
        \end{equation}
    since $\Phi^K(G/H_+) \neq 0$ if and only if $K \le_G H$. For this, see the proofs of \cite[Lemma~4.10]{BalmerSanders17} and \cite[Lemma~3.12]{bgh_balmer}. Finite localization with respect to the thick $\otimes$-ideal $\cat C_{\calF}$ has an associated idempotent triangle
        \[
            e_{\calF} \to \unit \to f_{\calF} \to \Sigma e_{\calF}
        \]
    which recovers the usual isotropy cofiber sequence: $e_{\calF} \simeq E\calF_+$ and $f_{\calF} \simeq \widetilde E\calF$.
\end{Rem}

\begin{Rem}\label{rem:pH-weakly-visible}
    For any closed subgroup $H \le G$, consider the family $\calF_{\le H}$ consisting of all closed subgroups conjugate to a subgroup of $H$ and the family $\calF_{< H}$ consisting of all closed subgroups conjugate to a proper subgroup of $H$. By \eqref{eq:YF}, we have
        \begin{equation}\label{pH-weak-vis}
            \{\mathfrak p_H\} = Y_{\calF_{\le H}} \cap (Y_{\calF_{< H}})^c
        \end{equation}
	which shows that $\frak p_H$ is weakly visible.
\end{Rem}

\begin{Lem}\label{lem:rational_weakly_noetherian}
    For any compact Lie group $G$, $\Spc(\SH_{G,\bbQ}^c)$ is weakly noetherian. 
\end{Lem}

\begin{proof}
    By \cref{thm:spectrumrationalgspectra}, every prime is of the form $\frak p_H$ for some closed subgroup $H \le G$, and these are weakly visible by \cref{rem:pH-weakly-visible}.
\end{proof}

\begin{Lem}\label{lem:idempotent-for-pH}
    The Balmer--Favi idempotent for the prime $\frak p_H \in \Spc(\SH_{G,\bbQ}^c)$ is
        \[
            g(\frak p_H) = (E\calF_{\le H})_+ \wedge \widetilde E\calF_{< H}.
        \] 
\end{Lem}

\begin{proof}
    This is immediate from \eqref{pH-weak-vis}.
\end{proof}

\begin{Rem}\label{rem:char_kappa}
    The idempotent $g(\frak p_H)$ has the property that 
        \begin{equation}\label{eq:isotropy_dirac}
            \Phi^K(g(\frak p_H) ) =
                \begin{cases}
                S^0 & \text{if } K \sim_G H\\
                0 & \text{otherwise}.
                \end{cases}
        \end{equation}
    In other words, the geometric isotropy of $g(\frak p_H)$ consists precisely of those subgroups of $G$ conjugate to $H$. This follows from \cref{lem:idempotent-for-pH} because
        \[
            \Phi^K(E\calF_+) = 
                \begin{cases}
                0 & \text{ if } K \not \in \calF \\
                S^0 & \text{ if } K \in \calF
                \end{cases}
        \]
    and
        \[
            \Phi^K(\widetilde E\calF) = 
                \begin{cases}
                0 & \text{ if } K  \in \calF \\
                S^0 & \text{ if } K \not \in \calF.
                \end{cases}
        \]
\end{Rem}

\begin{Rem}\label{rem:greenlees_compare}
    By smashing the cofiber sequence 
        \[
            (E\calF_{< H})_+  \to \unit \to \widetilde E\calF_{< H} \to \Sigma (E\calF_{< H})_+
        \]
    with the idempotent $( E\calF_{\le H})_+$, we see that $g(\frak p_H)$ fits in a cofiber sequence
        \begin{equation}\label{eq:equivariant_kappa}
            (E\calF_{< H})_+  \to (E\calF_{\le H})_+ \to g(\frak p_H) \to \Sigma (E\calF_{< H})_+.
        \end{equation}
    This shows that $g(\frak p_H)$ is the $G$-suspension spectrum of the space denoted $E\langle (H)_G \rangle$ in \cite[Definition~2.1]{Greenlees19_rational}.
\end{Rem}

\begin{Rem}
    Since $\Spc(\SH_{G,\bbQ}^c)$ is not noetherian in general, the local-to-global principle for $\SH_{G,\bbQ}$ is not a consequence of \cref{thm:local-to-global}. Nevertheless, it does always hold:
\end{Rem}

\begin{Prop}[Greenlees]\label{prop:local-to-global-rational}
    The local-to-global principle holds for $\SH_{G,\bbQ}$. 
\end{Prop}

\begin{proof}
    By \cite[Lemma~2.2]{Greenlees19_rational}, we have 
        \[
            \SH_{G,\bbQ} = \Loco{E\langle(H)_G \rangle \mid H \le G} = \Loco{g(\frak p_H) \mid H \le G}.
        \]
    In particular, $\unit \in \Loco{g(\frak p_H) \mid H \le G}$. The local-to-global principle then follows from \cref{lem:tensor-in}.
\end{proof}

\begin{Def}
    We let $\SH_{G,\langle H \rangle,\bbQ}$ denote the full subcategory of $\SH_{G,\bbQ}$ consisting of those rational $G$-spectra which are either contractible or have geometric isotropy conjugate to $H$. This category is denoted \textbf{$G$-spectra$\langle H \rangle$} in \cite{Greenlees19_rational}. 
\end{Def}

\begin{Prop}\label{prop:equivalence_gath}
    For any closed subgroup $H \le G$, we have
        \[
            \Loco{g(\frak p_H)} = \SH_{G,\langle H \rangle,\bbQ}
        \]
    as full subcategories of $\SH_{G,\bbQ}$.
\end{Prop}

\begin{proof}
    Because geometric fixed points form a coproduct-preserving tt-functor, it follows from \eqref{eq:isotropy_dirac} that $\Loco{g(\frak p_H)} \subseteq \Sp_{G,\langle H \rangle,\bbQ}$. Conversely, suppose that $X \in \SH_{G,\langle H \rangle,\bbQ}$. By \eqref{eq:equivariant_kappa}, there is a natural map $(E\calF_{\le H})_+ \wedge X \to g(\frak p_H) \wedge X$ and there is also a natural map $(E\calF_{\le H})_+ \wedge X \to X$. Using the joint conservativity of geometric fixed points (\cref{rem:geometric-isotropy}),  both these maps are equivalences:
        \[
            X \xleftarrow{\sim}  (E\calF_{\le H})_+ \wedge X \xrightarrow{\sim} g(\frak p_H) \wedge X. 
        \]
    In particular, $X \in \Loco{g(\frak p_H)}$, as required. 
\end{proof}

\begin{Thm}[Greenlees]\label{thm:rational_stratification}
    For any compact Lie group $G$, the category of rational $G$-equivariant spectra $\SH_{G,\bbQ}$ is stratified.
\end{Thm}

\begin{proof}
    By \cref{thm:equiv-strat}, along with \cref{lem:rational_weakly_noetherian,prop:local-to-global-rational}, we are reduced to showing that $\Loco{g(\frak p_H)}$ is a minimal localizing $\otimes$-ideal of $\SH_{G,\bbQ}$. By \cref{prop:equivalence_gath}, $\Loco{g(\frak p_H)} = \SH_{G,\langle H \rangle,\bbQ}$ and minimality then follows from \cite[Corollary~4.7]{Greenlees19_rational}.
\end{proof}

\begin{Rem}
    As noted in \cref{ex:rational-gspectra-gennoeth}, the space $\Spc(\SH_{G,\bbQ}^c)$ is generically noetherian, and so \cref{thm:rational_stratification} along with \cref{thm:gennoethtelescopeconj} show that the telescope conjecture holds for $\SH_{G,\bbQ}$. 
\end{Rem}

\begin{Rem}
    One can also use \cref{thm:universal-support} to show that the notion of support used by Greenlees coincides with the Balmer--Favi support. Indeed, the notion of support used by Greenlees is the geometric isotropy of $X \in \SH_{G,\bbQ}$ (\cref{rem:geometric-isotropy}). It follows directly from the definitions and properties of geometric fixed points that axioms (a) to (e) of \cref{def:axiomaticsupp} hold, and Greenlees shows that this notion of support stratifies $\SH_{G,\bbQ}$. Moreover, \cite[Remark~1.4 and Theorem~8.4]{Greenlees19_rational} verify that the equivalent conditions of \cref{thm:universal-support} are satisfied. Therefore, \cref{thm:universal-support} implies that the support defined by geometric isotropy agrees with the support defined by Balmer--Favi.
\end{Rem}

\part{Applications: Spectral Mackey functors}\label{part:spectral-mackey}

Our final objective is to study spectral Mackey functors \cite{Barwick17} with coefficients in a commutative ring spectrum. We will build on the work of \cite{PatchkoriaSandersWimmer22} which studies $\HZ$-valued spectral Mackey functors (a.k.a.~derived Mackey functors).  Our main goal is to prove that the category of spectral Mackey functors valued in \mbox{$\bbE$-modules} is stratified whenever the category of $\bbE$-modules is stratified.  This is the content of \cref{thm:mackeystrat}.  For example, it establishes stratification for $\HZ$-valued spectral Mackey functors (\cref{cor:dermackstratified}), as well as $E(n)$-local spectral Mackey functors (\cref{cor:En-local-mack-strat}). Throughout, $G$ is a finite group.

\section{Derived categories of ring spectra with trivial $G$-action}\label{sec:spectrum-E-mackey}

Spectral Mackey functors can be identified with modules over certain equivariant ring spectra. In this section, we study the Balmer spectrum of these module categories. Their connection with spectral Mackey functors will be explained in \cref{sec:equivalence-spectral-Mackey}. Our present goal is to show that the method used in \cite{PatchkoriaSandersWimmer22} to compute the Balmer spectrum of the category of $\HZ$-valued spectral Mackey functors~$\Der(\HZ_G)$ actually holds in much greater generality, replacing $\HZ_G\coloneqq \triv_G(\HZ)$ with $\triv_G(\bbE)$ for a general commutative ring spectrum $\bbE \in \CAlg(\Sp)$.

\begin{Rem}
	Here $\triv_G \colon \Sp \to \Sp_G$ denotes the essentially unique colimit-preserving symmetric monoidal functor from the $\infty$-category of spectra to the $\infty$-category of $G$-spectra. It can be identified with the inflation functor $\infl_{G/G}^G$ under an equivalence $\Sp \simeq \Sp_{G/G}$. See \cite[Section 3]{PatchkoriaSandersWimmer22}.
\end{Rem}

\begin{Def}
    Let $\bbE \in \CAlg(\Sp)$ be a commutative algebra in the symmetric monoidal stable $\infty$-category of spectra $\Sp$. For any finite group $G$, define $\bbE_G \coloneqq \triv_G(\bbE) \in \CAlg(\SpG)$ and let $\DEG \coloneqq \Ho(\EG\MMod_{\SpG})$.
\end{Def}

\begin{Exa}
    Taking $\bbE \coloneqq S^0$ to be the sphere spectrum, we obtain the equivariant stable homotopy category $\SHG \coloneqq \Ho(\SpG) = \DEG$.
\end{Exa}

\begin{Rem}\label{rem:gfp_finite_localization}
    With the symbol $\HZ$ replaced with $\bbE$, the arguments in \cite[Section~3]{PatchkoriaSandersWimmer22} go through verbatim. In particular, the tensor-triangulated categories~$\DEG$ satisfy Properties (A) through (I) listed at the start of \cite[Section~2]{PatchkoriaSandersWimmer22}. For example, there is an adjunction
        \[
            F_G : \SHG\leftrightarrows \Der(\bbE_G) : U_G 
        \]
    where the left adjoint (given by extension-of-scalars) is a tt-functor, and $\DEG$ is rigidly-compactly generated by $\SET{F_G(G/H_+)}{H \le G}$. The categories $\DEG$ also satisfy the key property (F) that for a normal subgroup $N \lenormal G$, the composite
        \[
            \Der(\bbE_{G/N}) \xra{\infl_{G/N}^G} \Der(\bbE_G) \to \Der(\bbE_G)/\Loco{F_G(G/H_+) \mid H \not\supseteq N}
        \]
    is an equivalence (see \cite[\S 2(F) and Proposition~3.22]{PatchkoriaSandersWimmer22}) and we define
        \[
            \tilde \Phi^{N,G} \colon \Der(\bbE_G) \to \Der(\bbE_{G/N})
        \]
    to be the finite localization
        \[
            \Der(\bbE_G) \to \Der(\bbE_G)/\Loco{F_G(G/H_+) \mid H \not\supseteq N} \cong \Der(\bbE_{G/N}).
        \]
    For any subgroup $H \le G$, we then define the (absolute) geometric $H$-fixed point functor $\Phi^{H,G} \colon \Der(\bbE_G) \to \Der(\bbE)$ as the composite
        \begin{equation}\label{eq:absolute-geom-def}
            \Der(\bbE_G) \xrightarrow{\res^G_H} \Der(\bbE_H) \xrightarrow{\tilde \Phi^{H,H}} \Der(\bbE_{H/H}) \cong \Der(\bbE). 
        \end{equation}
    It only depends (up to natural isomorphism) on the conjugacy class of $H$; see \cite[\S2(H)]{PatchkoriaSandersWimmer22}. Note that the base category which serves as the target for the absolute geometric fixed point functors is the nonequivariant derived category $\Der(\bbE) \coloneqq \Ho(\bbE\MMod_{\Sp})$. It follows from the definitions and property (F) that $\Phi^{H,G}$ is split by $\triv_G \colon \Der(\bbE) \to \Der(\bbE_G)$:
        \begin{equation}\label{eq:triv-splits}
            \Phi^{H,G} \circ \triv_G \simeq \Id_{\Der(\bbE)}.
        \end{equation}
    Finally observe that for $H = G$, the restriction functor in \eqref{eq:absolute-geom-def} is the identity and the geometric fixed point functor $\Phi^G\coloneqq \Phi^{G,G}\colon \DEG \to \DE$ is just a finite localization. 
\end{Rem}

\begin{Def}
    Because $\Phi^{H,G}\colon\DEG \to \DE$ is a tt-functor between rigidly-compactly generated tt-categories there is an induced map
        \[
            \Spc(\Phi^{H,G}) \colon \Spc(\Der(\bbE)^c) \to \Spc(\Der(\bbE_G)^c)
        \] 
    on Balmer spectra. For each prime $\frak  p \in \Spc(\DE^c)$ and subgroup $H \le G$ we define the prime $\cat P_G(H,\frak p) \in \Spc(\DEG^c)$ by
        \[
            \cat P_G(H,\frak p) \coloneqq \Spc(\Phi^{H,G})(\frak p) = \SET{x \in \DEG^c}{\Phi^{H,G}(x) \in \frak p} \in \Spc(\DEG^c).
        \]
\end{Def}

\begin{Prop}\label{prop:all-primes}
    Every prime of $\Spc(\DEG^c)$ is of the form $\cat P_G(H,\frak p)$ for some subgroup $H \le G$ and $\frak p \in \Spc(\DE^c)$.
\end{Prop}

\begin{proof}
    For any subgroup $H\le G$, the restriction functor $\DEG \to \Der(\bbE_H)$ has a left adjoint \cite[Remark~3.16]{PatchkoriaSandersWimmer22}, hence the proof of \cite[Proposition~2.13]{PatchkoriaSandersWimmer22} also holds, showing that the image on Balmer spectra of restriction is precisely $\supp(F_G(G/H_+))$. With this in hand, and the fact that $\Phi^{G}$ is the usual finite localization (\cref{rem:gfp_finite_localization}), the inductive argument in \cite[Theorem~2.22]{PatchkoriaSandersWimmer22} goes through verbatim.
\end{proof}

\begin{Lem}\label{lem:varphiG}
    Let $\varphi_G \colon \Spc(\DEG^c) \to \Spc(\SHGc)$ be the map induced by the extension-of-scalars functor $F_G\colon \SHG \to \DEG$. Also write $\varphi \colon \Spc(\DE^c) \to \Spc(\SHc)$ for the $G=1$ case. Then 
		\[
			\varphi_G(\cat P_G(H,\frak p)) = \cat P_G(H,\varphi(\frak p)) \in \Spc(\SHGc)
		\]
    for any $H \le G$ and $\frak p \in \Spc(\DE^c)$.
\end{Lem}

\begin{proof}
    The commutative diagram of tt-functors
        \[
            \begin{tikzcd}
            \SHG \arrow[d, "F_G"'] \arrow[r, "{\Phi^{H,G}}"] & \SH \arrow[d, "F=F_1"] \\
            \Der(\bbE_G) \arrow[r, "{\Phi^{H,G}}"']           & \Der(\bbE)       
            \end{tikzcd}
        \] 
    (see \cite[Lemma 2.4 and Remark 2.6]{PatchkoriaSandersWimmer22}) shows that $\varphi_G(\cat P_G(H,\frak p)) = \cat P_G(H,\varphi(\frak p))$ in $\Spc(\SHGc)$.
\end{proof}

\begin{Prop}\label{prop:prime-overlap}
    Let $H,K \le G$ and $\frak p,\frak q \in \Spc(\DE^c)$. Then $\cat P_G(H,\frak p) = \cat P_G(K,\frak q)$ in $\Spc(\DEG^c)$ if and only if $H \sim_G K$ and $\frak p = \frak q$ as primes of $\DE^c$.
\end{Prop}

\begin{proof}
    If $H \sim_G K$ are conjugate subgroups then the geometric fixed point functors $\Phi^{H,G},\Phi^{K,G}\colon \SHG \to \SH$ are naturally isomorphic. The same is true of the induced functors $\Phi^{H,G},\Phi^{K,G}\colon \DEG \to \DE$; see \cite[\S2(H)]{PatchkoriaSandersWimmer22}. Thus, if ${H \sim_G K}$ then $\cat P_G(H,\frak p) = \cat P_G(K,\frak p)$ in $\Spc(\DEG^c)$ for each $\frak p \in \Spc(\DE^c)$. This gives the $(\Leftarrow)$ direction. On the other hand, if $\cat P_G(H,\frak p) =\cat P_G(K,\frak q)$ in $\Spc(\DEG^c)$ then $\cat P_G(H,\varphi(\frak p)) = \cat P_G(K,\varphi(\frak q))$ in $\Spc(\SH(G)^c)$ by \cref{lem:varphiG}. Hence, by \cite[Theorem 4.14]{BalmerSanders17}, $H \sim_G K$ are $G$-conjugate. Thus, $\cat P_G(H,\frak p) = \cat P_G(K,\frak q) = \cat P_G(H,\frak q)$ in $\Spc(\DEG^c)$. This implies $\frak p=\frak q$ since $\Phi^{H,G}\circ \triv_G \simeq \Id_{\DE}$ (\cref{rem:gfp_finite_localization}).
\end{proof}

\begin{Thm}\label{thm:spcDEG}
    Let $G$ be a finite group and let $\bbE \in \CAlg(\Sp)$ be a commutative ring spectrum. As a set $\Spc(\DEG^c)$ consists of primes $\cat P_G(H,\frak p)$ for $H \le G$ and $\frak p \in \Spc(\DE^c)$. Moreover:
    \begin{enumerate}
        \item Two such primes $\cat P_G(H,\frak p)$ and $\cat P_G(K,\frak q)$ coincide precisely when the subgroups $H$ and $K$ are $G$-conjugate and $\frak p=\frak q$.
        \item For each conjugacy class of subgroups $H\le G$, the geometric fixed point functor $\Phi^{H,G}\colon \Der(\bbE_G) \to \Der(\bbE)$ induces a homeomorphism 
                \[
                    \Spc(\Der(\bbE)^c) \xra{\cong} \SET{\cat P_G(H,\frak p)}{\frak p \in \Spc(\Der(\bbE)^c)} \subset \Spc(\Der(\bbE_G)^c).
                \]
        \item If $\Spc(\DE^c)$ is a noetherian space, then so is $\Spc(\DEG^c)$.
    \end{enumerate}
\end{Thm}

\begin{proof}
    \Cref{prop:all-primes} establishes that every prime is of the form $\cat P_G(H,\frak p)$ and \cref{prop:prime-overlap} characterizes when two such primes coincide. Since $\Phi^{H,G} \circ \triv_G \simeq \Id_{\DE}$ (\cref{rem:gfp_finite_localization}), the composite
        \[
            \Spc(\DE^c) \xra{\Spc(\Phi^{H,G})} \Spc(\DEG^c) \xra{\Spc(\triv_G)} \Spc(\DE^c)
        \]
    is the identity. Thus, $\Spc(\Phi^{H,G})$ embeds $\Spc(\DE^c)$ homeomorphically onto its image in $\Spc(\DEG^c)$. This proves part (b). If $\Spc(\DE^c)$ is noetherian, then the images of $\Spc(\DE^c)$ under the various geometric fixed point functors provide a finite cover of $\Spc(\DEG^c)$ by noetherian spaces. Hence part (c) follows as in the proof of \cite[Proposition~2.38]{PatchkoriaSandersWimmer22}.
\end{proof}

\begin{Rem}
    We thus understand $\Spc(\DEG^c)$ completely as a set: It is covered by disjoint copies of $\Spc(\DE^c)$, one for each conjugacy class of subgroups $H \le G$. However, there could be topological interaction between the copies of $\Spc(\DE^c)$ for different conjugacy classes of subgroups. For the case $\bbE=\HZ$, the topology has been completely determined in \cite[Remark~2.39]{PatchkoriaSandersWimmer22}. Note also that the above statements specialize to basic results about $\Spc(\SHGc)$ from \cite[Section~4]{BalmerSanders17} by taking $\bbE=S^0$ to be the sphere spectrum. In this case, the topology is understood for finite abelian groups \cite{BHNNNS19} and extraspecial 2-groups \cite{KuhnLloyd20pp}, while we have only partial information for general finite groups $G$.
\end{Rem}

\begin{Cor}
    Let $\varphi_G\colon \Spc(\Der(\bbE_G)^c) \to \Spc(\SHGc)$ be the map induced by extension-of-scalars. Also write $\varphi \colon  \Spc(\Der(\bbE)^c) \to \Spc(\SHc)$ for the $G=1$ case.
    \begin{enumerate}
        \item If $\varphi$ is injective then $\varphi_G$ is injective.
        \item If $\varphi$ is surjective then $\varphi_G$ is surjective.
    \end{enumerate}
\end{Cor}

\begin{proof}
    This is a consequence of \cref{lem:varphiG} by applying \cref{thm:spcDEG} both to the given $\bbE$ and to the sphere spectrum~$S^0$. Indeed, every prime of $\DEG^c$ is of the form $\cat P_G(H,\frak p)$ and $\varphi_G(\cat P_G(H,\frak p)) = \cat P_G(H,\varphi(\frak p))$ by \cref{lem:varphiG}. Thus, if $\varphi_G(\cat P_G(H,\frak p)) = \varphi_G(\cat P_G(K,\frak q))$ then we have an equality of primes $\cat P_G(H,\varphi(\frak p)) = \cat P_G(K,\varphi(\frak q))$ in $\Spc(\SHGc)$. Hence $\varphi(\frak p)=\varphi(\frak q)$ in $\Spc(\SHc)$ and the subgroups~$H$ and $K$ are $G$-conjugate. Thus, if $\varphi$ is injective, then $\cat P_G(H,\frak p) = \cat P_G(K,\frak q)$ in $\Spc(\DEG^c)$. This establishes part (a). Part (b) is an immediate corollary of \cref{thm:spcDEG} and \cref{lem:varphiG}.
\end{proof}

\begin{Rem}\label{rem:continuousbijection}
    If $\varphi\colon\Spc(\DE^c) \to \Spc(\SHc)$ is injective then $\varphi_G\colon\Spc(\DEG^c) \to \Spc(\SHGc)$ is a continuous bijection onto its image 
        \[
            \im(\varphi_G)=\SET{ \cat P_G(H,\varphi(\frak p))}{H\le G, \frak p \in \Spc(\DE^c)} \subset \Spc(\SHGc).
        \]
    However, we have not established that if $\varphi$ is a homeomorphism onto its image, then the same is true for $\varphi_G$. This point will be relevant for the following:
\end{Rem}

\begin{Exa}\label{exa:en_spectral_mackey_1}
    Let $\bbE =L_nS^0$, so that $\DE \coloneqq \Ho(\bbE\MMod_{\Sp}) \cong \cat S_{E(n)}$ is the $E(n)$-local stable homotopy category (see \cref{Rem:e_n}). The spectrum of $\cat S^c_{E(n)}$ was described in \cref{thm:hovey-strickland}:
        \[
            \Spc(\cat S_{E(n)}^c) = \cat P_{n} - \cdots - \cat P_{1} - \cat P_0.
        \]
    We let $\SpGEn \coloneqq \bbE_G\MMod_{\Sp_G}$ and $\SHGEn \coloneqq \Ho(\SpGEn) = \DEG$ for $\bbE = L_nS^0$. Applying \cref{thm:spcDEG} we deduce that, as a set, $\Spc(\SHGEn^c)$ consists of primes $\cat P_G(H,k)$ for $H \le G$ and $0 \le k \le n$. Two such primes $\cat P_G(H,k)$ and $\cat P_G(K,\ell)$ coincide precisely when the subgroups $H$ and $K$ are $G$-conjugate and $k=\ell$. Moreover, as a consequence of \cref{thm:spcDEG} and \cref{rem:continuousbijection}, the following hold:
    \begin{enumerate}
        \item For each conjugacy class of subgroups $H \le G$, the geometric fixed point functor $\Phi^{H,G}\colon \SHGEn^c \to \cat S_{E(n)}^c$ induces a homeomorphism of $\Spc(\cat S_{E(n)}^c)$ onto the subspace 
                \[
                    \SET{\cat P_G(H,k)}{ 0 \le k \le n} \subset \Spc(\SHGEn^c).
                \]
        \item The map $\varphi_G\colon \Spc(\SHGEn^c) \to \Spc(\SHG^c)$ induced by extension-of-scalars is injective. It provides a continuous bijection onto
                \[
                    \SET{\cat P_G(H,p,k)}{H \le G, 0\le k \le n} \subset \Spc(\SHG^c)
                \]
              where $p$ is the prime number implicit in our choice of $E(n)$.
    \end{enumerate}
\end{Exa}

\begin{Rem}
    In other words, by part (b), the spectrum $\Spc(\SHGEn^c)$ bijects onto the first $n$ chromatic layers of the spectrum $\Spc(\SHGp^c)\subset \Spc(\SHG^c)$ of the $p$-local equivariant stable homotopy category (recall \cite[Corollary~5.25]{BalmerSanders17} and \cref{not:Cn}). However, we are unable to prove that this continuous bijection is a homeomorphism.
\end{Rem}

\begin{Conj}\label{conj:homeomorphism}
    For any finite group $G$, the continuous injective map
        \[
            \varphi_G\colon \Spc(\SHGEn^c) \to \Spc(\SHG^c)
        \]
    is a homeomorphism onto its image.
\end{Conj}

\section{Equivalence of \texorpdfstring{$\Der(\bbE_G)$}{D(E_G)} with spectral Mackey functors}\label{sec:equivalence-spectral-Mackey}

We now connect the equivariant module categories of the previous section with categories of spectral Mackey functors. These results follow \cite[Section 4]{PatchkoriaSandersWimmer22} where we refer the reader for further details and discussion.

\begin{Not}
    We write $\AGeff$ for Barwick's effective Burnside $\infty$-category \cite{Barwick17} whose objects are finite $G$-sets and whose $n$-simplicies are $n$-fold spans of finite \mbox{$G$-sets.} For an additive $\infty$-category $\cat C$, we write $\Fun_{\add}(\AGeff,\cat C)$ for the $\infty$-category of additive functors from $\AGeff$ to $\cat C$. This is a smashing localization of the functor category $\Fun(\AGeff, \cat C)$ and can be equipped with the localized Day convolution product \cite[Section 3]{BarwickGlasmanShah20}.
\end{Not}

\begin{Def}
    Let $\bbE \in \CAlg(\Sp)$ be a commutative algebra in the symmetric monoidal stable $\infty$-category of spectra $\Sp$. For any finite group $G$, we have the $\infty$-category
        \[
            \Mack{G}{\bbE} \coloneqq \Fun_{\add}(\AGeff,\bbE\MMod_{\Sp})
        \]
    of spectral $G$-Mackey functors valued in $\bbE$-module spectra (or with $\bbE$-module coefficients).
\end{Def}

\begin{Prop}\label{prop:mackmorita}
    For any finite group $G$ and commutative algebra $\bbE \in \CAlg(\Sp)$, there is an equivalence
        \[
            \Mack{G}{\bbE}\simeq \bbE_G\MMod_{\Sp_G}
        \]
    of presentably symmetric monoidal stable $\infty$-categories.
\end{Prop}

\begin{proof}
    As explained in \cite[Remark~4.5]{PatchkoriaSandersWimmer22}, there are functors
        \[
            L \colon \Fun(\AGeff,\Sp) \to \Fun_{\add}(\AGeff,\Sp)
        \]
    and 
        \[
            F_1 \colon \Sp \to \Fun(\AGeff,\Sp)
        \]
    given respectively as the localization and the adjoint to evaluation at $G/G \in \AGeff$. As in \cite[Remark~4.5 and Corollary~4.6]{PatchkoriaSandersWimmer22}, we have an equivalence of  commutative algebras $\bbE_G \simeq LF_1(\bbE)$ which gives rise to an equivalence of symmetric monoidal $\infty$-categories
        \[
            \bbE_G\MMod_{\SpG} \simeq LF_1(\bbE)\MMod_{\SpG}. 
        \]
    Replacing all instances of $\HZ$ with $\bbE$, the argument in \cite[Proposition~4.9]{PatchkoriaSandersWimmer22} establishes an equivalence of symmetric monoidal $\infty$-categories 
        \[
            LF_1(\bbE)\MMod_{\Sp_G} \simeq \Fun_{\add}(\AGeff,\bbE\MMod_{\Sp}) 
        \]
    and the latter is $\Mack{G}{\bbE}$.
\end{proof}

\begin{Exa}\label{exa:en_spectral_mackey_2}
    We return to \cref{exa:en_spectral_mackey_1}. By \cref{prop:mackmorita}, there is an equivalence of symmetric monoidal stable $\infty$-categories 
        \[
            \begin{split}
            \SpGEn= \triv_G(L_nS^0)\MMod_{\SpG} &\simeq \Fun_{\add}(\AGeff,{L_nS^0}\MMod_{\Sp}) \\
            & \simeq \Fun_{\add}(\AGeff,\cat S_{E(n)}). 
            \end{split}
        \]
    Therefore, $\SpGEn$ can alternatively be regarded as the category of spectral Mackey functors with coefficients in $E(n)$-local spectra.
\end{Exa}

\begin{Rem}
    Since Bousfield localization at $E(n)$ is smashing on $\Sp$, it follows from \cite[Proposition 3.12]{Carrick22} that Bousfield localization at $E(n)_G \coloneqq \triv_G(E(n))$ is smashing on $\SpG$. Thus
        \[
            L_{E(n)_G}\SpG \simeq L_{E(n)_G}(S^0)\MMod_{\SpG}.
        \]
    Moreover, $L_{E(n)_G}(S^0) \simeq \triv_G(L_nS^0)$ by \cite[Proposition 3.2]{Carrick22}, so that
        \[
            L_{E(n)_G}\SpG \simeq \SpGEn. 
        \]
    This gives a third interpretation of the category $\SpGEn$ as the Bousfield localization of $\SpG$ at $E(n)_G$. This justifies calling $\SpGEn$ the category of $E(n)$-local spectral Mackey functors rather than the exhausting ``spectral Mackey functors valued in $E(n)$-local spectra.'' In any case, one can see using \cite[Proposition~3.3]{PatchkoriaRoitzheim20} that when $n = 1$, this category agrees with the category studied by Patchkoria and Roitzheim in \cite{PatchkoriaRoitzheim20}. 
\end{Rem}

\section{Stratification for spectral Mackey functors}\label{sec:strat-for-spectral-mack}

The goal of this section is to prove the following theorem, which shows that the category of spectral Mackey functors with coefficients in a commutative ring spectrum $\bbE$ is stratified provided the non-equivariant category of coefficients $\DE$ is stratified and has a noetherian spectrum $\Spc(\DE^c)$.

\begin{Thm}\label{thm:mackeystrat}
    Let $\bbE \in \CAlg(\Sp)$ be a commutative ring spectrum such that $\Spc(\DE^c)$ is noetherian. If $\DE$ is stratified, then $\DEG$ is stratified for any finite group $G$.
\end{Thm}

\begin{Rem}
    In light of \cref{prop:mackmorita}, we have stated the theorem with $\DEG$ in place of the equivalent tt-category $\Ho(\Mack{G}{\bbE})$. Since the spectrum of $\DEG$ is noetherian by \cref{thm:spcDEG}(c), it remains (by \cref{thm:local-to-global} and \cref{thm:equiv-strat}) to verify the minimality of the local categories for each prime ideal. These prime ideals are of the form $\cat P_G(H,\frak p)$ for a (conjugacy class of) subgroup $H \le G$ and a prime ideal $\frak p \in \Spc(\DE^c)$. The idea of the proof may then be outlined via the following commutative diagram:
        \begin{equation}\label{eq:proofofmackeystrat}
            \begin{tikzcd}
                \DEG \ar[r,"\res^G_H"'] \ar[rr,bend left=15,"\Phi^{H,G}"] & \DEH \ar[r,"\Phi^H"'] \ar[d] & \DE \ar[d] \\
                & \DEH/\langle\cat P_H(H,\frak p)\rangle \ar[r,"\sim"'] & \DE/\langle\frak p\rangle,
            \end{tikzcd}
        \end{equation}
    in which the vertical maps are given by the localization functors. Minimality at $\cat P_G(H,\frak p)$ will then be a consequence of the following three claims:
        \begin{enumerate}
            \item The functor $\res^G_H$ is finite \'etale (\cref{def:finite-etale}) and the fiber over $\cat P_G(H,\frak p)$ of the induced map on spectra  is precisely $\{\cat P_H(H,\frak p)\}$. 
            \item The geometric fixed point functor $\Phi^H$ induces an equivalence on localizations $\DEH/\langle\cat P_H(H,\frak p)\rangle \xrightarrow{\sim}  \DE/\langle\frak p\rangle$.
            \item The vertical functors preserve and reflect minimality at $\cat P_H(H,\frak p)$ and $\frak p$, respectively. 
        \end{enumerate}
\end{Rem}

\begin{Lem}\label{lem:resfiniteetale}
    For any subgroup $H \le G$, the restriction functor
        \[ 
            \res^G_H\colon  \DEG \longrightarrow \DEH
        \]
    is finite \'{e}tale. Moreover, if $\varphi\colon \Spc(\DEH^c) \to \Spc(\DEG^c)$ denotes the induced map on spectra, we have 
        \[
            \varphi^{-1}(\varphi(\{\cat P_H(H,\frak p)\})) = \{\cat P_H(H,\frak p)\}
        \]
    for any $\frak p \in \Spc(\DE^c)$.
\end{Lem}

\begin{proof}
    The restriction functor is finite \'etale by \cite[Example~5.8]{Sanders22}. As in \cite[Lemma~2.12]{PatchkoriaSandersWimmer22}, we can verify that $\varphi(\cat P_H(K,\frak q)) = \cat P_G(K,\frak q)$ for any subgroup $K \le H$ and any prime $\frak q \in \Spc(\DE^c)$. In particular, we compute
        \begin{align*}
            \varphi^{-1}(\{\varphi(\cat P_H(H,\frak p))\}) & = \varphi^{-1}(\{\cat P_G(H,\frak p)\}) \\
            & = \SET{\cat P_H(K,\frak q)}{\cat P_G(K,\frak q) = \cat P_G(H,\frak p)}. 
        \end{align*}
    \Cref{thm:spcDEG} implies that the last set consists only of the single prime $\cat P_H(H,\frak p)$ since if $K \le H$ is $G$-conjugate to $H$ then $K=H$.
\end{proof}

\begin{Lem}\label{lem:transferquotient}
    For any prime $\frak p \in \Spc(\DE^c)$, there is a commutative diagram
        \[
            \begin{tikzcd}
                \DEH \ar[r,"\Phi^H"] \ar[d] & \DE \ar[d] \\
                \DEH/\langle\cat P_H(H,\frak p)\rangle \ar[r,"\sim"] & \DE/\langle\frak p\rangle
            \end{tikzcd}
        \]
    where the vertical arrows are the localization functors to the corresponding local categories (\cref{def:primeverdierquotient}) and the bottom functor is an equivalence.
\end{Lem}

\begin{proof}
    By construction, the geometric fixed point functor $\Phi^H\colon \DEH \to \DE$ is a finite localization (\cref{rem:gfp_finite_localization}). Since $\cat P_H(H, \frak p) = (\Phi^H)^{-1}(\frak p)$ by definition, the result is a special case of \cref{prop:localfinitelocalization}.
\end{proof}

\begin{proof}[Proof of \cref{thm:mackeystrat}]
    By assumption, the spectrum $\Spc(\DE^c)$ is noetherian, hence so is $\Spc(\DEG^c)$ by \cref{thm:spcDEG}(c). Therefore, the local-to-global principle holds for $\DEG$ by \cref{thm:local-to-global}. By \cref{thm:equiv-strat}, it remains to verify minimality (in the sense of \cref{def:stratified}) at all the primes of $\Spc(\DEG^c)$.

    By \cref{thm:spcDEG}, an arbitrary prime of $\Spc(\DEG^c)$ is of the form $\cat P_G(H,\frak p) = (\Phi^{H,G})^{-1}(\frak p) \in \Spc(\DEG^c)$ for a subgroup $H \le G$ and a prime $\frak p \in \Spc(\DE^c)$. Recall diagram \eqref{eq:proofofmackeystrat}. Since $\DE$ is stratified, it has minimality at $\frak p$. It follows from \cref{prop:single-prime-local} that the local category $\DE/\langle \frak p\rangle$ satisfies minimality at its unique closed point. By \cref{lem:transferquotient}, this is equivalent to the local category $\DEH/\langle \cat P_H(H,\frak p) \rangle$ having minimality at its closed point. Using \cref{prop:single-prime-local} again, we obtain minimality of $\DEH$ at the prime $\cat P_H(H,\frak p) \in \Spc(\Der(\bbE_H)^c)$.

    Finally, by \cref{lem:resfiniteetale} we can invoke \cref{thm:finite-etale} to conclude that $\DEG$ has minimality at $\cat P_G(H,\frak p) \in \Spc(\DEG^c)$.
\end{proof}

\begin{Rem}\label{rem:mackeystratop}
    The following converse to \cref{thm:mackeystrat} also holds: If $G$ is a finite group and $\bbE$ is a commutative ring spectrum such that $\DEG$ is stratified with (weakly) noetherian spectrum, then $\DE$ is also stratified with (weakly) noetherian spectrum. Indeed, just apply \cref{cor:ltg-passes-to-finite}, \cref{prop:single-prime-local} and \cref{thm:equiv-strat} to the finite localization $\Phi^G\colon\DEG \to \DE$.
\end{Rem}

\begin{Rem}
    If $\DE$ is a local category with unique closed point $\frak p \in \Spc(\DE^c)$ then the finite localization $\Phi^G\colon\DEG\to\DE$ is nothing but localization at the point $\cat P_G(G,\frak p) \in \Spc(\DEG^c)$. Indeed, this follows from \cref{lem:transferquotient}. In particular, $\DE\cong \DEG/\langle \cat P_G(G,\frak p) \rangle$ is the local category of $\DEG$ at the point $\cat P_G(G,\frak p)$.
\end{Rem}

\begin{Rem}\label{rem:topology_ideals}
    It may be interesting to note that we have obtained the stratification result of \cref{thm:mackeystrat} for $\DEG$ without using (or even having) complete knowledge of the topology of the Balmer spectrum $\Spc(\DEG^c)$. We note, moreover, that understanding the (radical) thick tensor-ideals of compact objects in a category $\cat T$ is equivalent to understanding the topological space $\Spc(\cat T^c)$ by Balmer's classification theorem \cite[Theorem 4.10]{Balmer05a}. On the other hand, when $\cat T$ is stratified the classification of localizing $\otimes$-ideals only depends on $\Spc(\cat C^c)$ as a set. In particular, \cref{thm:mackeystrat} allows us to find examples (see \cref{rem:no-thick-classification} below) where we understand the localizing $\otimes$-ideals but not the thick $\otimes$-ideals of compact objects. 
\end{Rem}

\begin{Cor}\label{cor:dermackstratified}
    For any finite group $G$, the category of derived Mackey functors $\DHZG$ is stratified.
\end{Cor}

\begin{proof}
    This is the special case of \cref{thm:mackeystrat} corresponding to $\bbE = \HZ$. The hypotheses on the base category $\Der(\HZ)\cong\Der(\bbZ)$ are provided by Neeman \cite{Neeman92a}; see \cref{thm:ring-strat} and \cref{exa:local-for-DR}, for example.
\end{proof}

\begin{Rem}
    More generally, the previous result applies to $\bbE = \HR$ where $R$ is any discrete commutative noetherian ring.
\end{Rem}

\begin{Cor}\label{cor:En-local-mack-strat}
    For any finite group $G$, prime number $p$, and $0 \le n < \infty$, the category of $E(n)$-local spectral Mackey functors $\SHGEn$ is stratified.
\end{Cor}

\begin{proof}
    Invoking \cref{thm:hovey-strickland} and \cref{thm:En-stratified}, we can apply \cref{thm:mackeystrat} with $\bbE=L_n S^0$. See \cref{exa:en_spectral_mackey_1,exa:en_spectral_mackey_2}.
\end{proof}

\begin{Rem}\label{rem:no-thick-classification}
    When $n > 1$ we do not fully understand the topology on $\Spc(\SHGEn^c)$ and so this gives a concrete example of a category for which we understand the localizing $\otimes$-ideals, but not the thick $\otimes$-ideals of compact objects.
\end{Rem}

\begin{Rem}
    The previous examples do not exhaust the cases for which \cref{thm:mackeystrat} applies. For example, if $\bbE$ is a commutative ring spectrum such that $\pi_*(\bbE)$ is a regular noetherian ring concentrated in even degrees, then $\DE$ is also stratified~\cite{DellAmbrogioStanley16}. In particular, this applies to $\mathbb{E} = E_n$, the $n$-th Morava $E$-theory.  
\end{Rem}

\addtocontents{toc}{\vspace{\normalbaselineskip}}
\printbibliography
	
\end{document}